\documentclass{article}

\usepackage{arxiv}

\usepackage[utf8]{inputenc} 
\usepackage[T1]{fontenc}    
\usepackage{hyperref}       
\usepackage{url}            
\usepackage{booktabs}       
\usepackage{amsfonts}       
\usepackage{nicefrac}       
\usepackage{microtype}      
\usepackage{lipsum}
\usepackage{mathptmx}
\usepackage{bigints}

\usepackage{sectsty}
\usepackage{graphicx}
\usepackage{amsmath}
\usepackage{hyperref}
\usepackage{subcaption}
\usepackage{placeins}
\usepackage{bm}

\usepackage{pdflscape,rotating}

\usepackage{algorithmic, algorithm2e, threeparttable, multirow}
\RestyleAlgo{boxruled}

\hypersetup{colorlinks,breaklinks,
           linkcolor=blue,urlcolor=blue,
           anchorcolor=blue,citecolor=blue}

\newcommand{\Am}{\mathbf{A}}
\newcommand{\Cm}{\mathbf{C}}
\newcommand{\Cb}{\mathbb{C}}
\newcommand{\Chom}{\mathbb{C}_{\mathrm{hom}}}
\newcommand{\Cbm}{\mathbb{C}_{\mathrm{m}}}
\newcommand{\Cbf}{\mathbb{C}_{\mathrm{f}}}

\newcommand{\ee}{\mathbf{e}}
\newcommand{\fm}{\mathbf{f}}

\newcommand{\gm}{\mathbf{g}}

\newcommand{\hb}{\mathbf{h}}

\newcommand{\mm}{\mathbf{m}}
\newcommand{\ngg}{n_{\mathbf{g}}}
\newcommand{\np}{n_{\mathbf{\ppm}}}

\newcommand{\ny}{n_{\mathbf{\Ym}}}

\newcommand{\ppm}{{\bm\theta}}
\newcommand{\rr}{\mathbf{r}}
\newcommand{\Sm}{\mathbf{S}}
\newcommand{\uu}{\mathbf{u}}
\newcommand{\vm}{\mathbf{v}}

\newcommand{\xm}{\mathbf{x}}
\newcommand{\Ym}{\pmb{\xi}}
\newcommand{\ym}{\mathbf{y}}

\newcommand{\Exp}{\mathbb{E}}

\newcommand{\kappaa}{\bm\kappa}
\newcommand{\thetaa}{\bm\theta}
\newcommand{\sigmaa}{\bm{\sigma}}
\newcommand{\epsl}{\bm{\epsilon}}
\newcommand{\xii}{\bm\xi}
\newcommand{\epsb}{\bm\epsilon}

\newcommand\Tstrut{\rule{0pt}{2.6ex}}         
\newcommand\Bstrut{\rule[-0.9ex]{0pt}{0pt}}   

\usepackage{pgfplots}
\usepackage{mathptmx}
\DeclareMathAlphabet{\mathcal}{OMS}{cmsy}{m}{n}

\usepackage{tikz}
\usetikzlibrary{matrix} 
\usetikzlibrary{arrows} 
\usetikzlibrary{arrows.meta}
\usetikzlibrary{calc} 
\usetikzlibrary{shapes}
\usetikzlibrary{fit}
\usetikzlibrary{positioning}
\usetikzlibrary{intersections,patterns,pgfplots.fillbetween}
\usetikzlibrary{calc}
\usetikzlibrary{decorations.pathreplacing}

\usepackage{ifpdf}
\ifpdf
  \usepackage{epstopdf}
\fi

\tikzstyle{block} = [draw,rectangle,thick,minimum height=2em,minimum width=2em]
\tikzstyle{sum} = [draw,circle,inner sep=0mm,minimum size=2mm]
\tikzstyle{connector} = [->,thick]
\tikzstyle{line} = [thick]
\tikzstyle{branch} = [circle,inner sep=0pt,minimum size=1mm,fill=black,draw=black]
\tikzstyle{guide} = []
\pgfplotsset{compat=1.16}

\pgfkeys{
  /pgf/arrow keys/.cd,
  pitch/.code={%
    \pgfmathsetmacro\pgfarrowpitch{#1}
    \pgfmathsetmacro\pgfarrowsinpitch{abs(sin(\pgfarrowpitch))}
    \pgfmathsetmacro\pgfarrowcospitch{abs(cos(\pgfarrowpitch))}
  },
}

\pgfdeclarearrow{
  name = Cone,
  defaults = {       
    length     = +3.6pt +5.4,
    width'     = +0pt +0.5,
    line width = +0pt 1 1,
    pitch      = +0, 
  },
  cache = false,     
  setup code = {},   
  drawing code = {   
    \pgfmathsetmacro\pgfarrowhalfwidth{.5\pgfarrowwidth}
    \pgfmathsetmacro\pgfarrowhalfwidthsin{\pgfarrowhalfwidth*\pgfarrowsinpitch}
    \pgfpathellipse{\pgfpointorigin}{\pgfqpoint{\pgfarrowhalfwidthsin pt}{0pt}}{\pgfqpoint{0pt}{\pgfarrowhalfwidth pt}}
    \pgfusepath{fill}
    \pgfmathsetmacro\pgfarrowlengthcos{\pgfarrowlength*\pgfarrowcospitch}
    \pgfmathparse{\pgfarrowlengthcos>\pgfarrowhalfwidthsin}
    \ifnum\pgfmathresult=1
      \pgfmathsetmacro\pgfarrowlengthtemp{\pgfarrowhalfwidthsin*\pgfarrowhalfwidthsin/\pgfarrowlengthcos}
      \pgfmathsetmacro\pgfarrowwidthtemp{\pgfarrowhalfwidth/\pgfarrowlengthcos*sqrt(\pgfarrowlengthcos*\pgfarrowlengthcos-\pgfarrowhalfwidthsin*\pgfarrowhalfwidthsin)}
      \pgfpathmoveto{\pgfqpoint{\pgfarrowlengthcos pt}{0pt}}
      \pgfpathlineto{\pgfqpoint{\pgfarrowlengthtemp pt}{ \pgfarrowwidthtemp pt}}
      \pgfpathlineto{\pgfqpoint{\pgfarrowlengthtemp pt}{-\pgfarrowwidthtemp pt}}
      \pgfpathclose
      \pgfusepath{fill}
    \fi
    \pgfpathmoveto{\pgfpointorigin}
  }
}

\title{Topology Optimization under Microscale Uncertainty using Stochastic Gradients}

\author{
  Subhayan De \\
  Aerospace Engineering Sciences\\
  University of Colorado\\
  Boulder, CO 80303 \\
  \texttt{Subhayan.De@colorado.edu} \\
   \And
 Kurt Maute \\
  Aerospace Engineering Sciences\\
  University of Colorado\\
  Boulder, CO 80303 \\
  \texttt{maute@colorado.edu} \\
  \And
 Alireza Doostan \\
  Aerospace Engineering Sciences\\
  University of Colorado\\
  Boulder, CO 80303 \\
  \texttt{doostan@colorado.edu} \\
}

\begin{document}
\maketitle

\begin{abstract}
This paper considers the design of structures made of engineered materials, accounting for uncertainty in material properties. We present a topology optimization approach {that optimizes the structural shape and topology at the macroscale assuming design-independent uncertain microstructures.} {}
The structural geometry at the macroscale is described by an explicit level set approach, and the macroscopic structural response is predicted by the eXtended Finite Element Method (XFEM). We describe the microscopic layout by either an analytic geometric model with uncertain parameters or {a level cut from} a Gaussian random field. The macroscale properties of the microstructured material are predicted by homogenization. 
Considering the large number of possible microscale configurations, one of the main challenges of solving such topology optimization problems is the computational cost of estimating the statistical moments of the cost and constraint functions and their gradients with respect to the design variables. 
Methods for predicting these moments, such as Monte Carlo sampling, and Taylor series and polynomial chaos expansions often require a large number of random samples resulting in an impractical computation. 
To reduce this cost, we propose an approach wherein, at every design iteration, we only use a small number of microstructure configurations to generate an independent, stochastic approximation of the gradients. 
These gradients are then used either with a gradient descent algorithm, namely Adaptive Moment (Adam), or the globally convergent method of moving asymptotes (GCMMA). 
Three numerical examples from structural mechanics are used to show that the proposed approach provides a computationally efficient way for macroscale topology optimization in the presence of microstructural uncertainty {and enables the designers to consider a new class of problems that are out of reach today with conventional tools}. 
\end{abstract}

\keywords{Topology optimization \and Microscale uncertainty \and Stochastic gradients}

\section{Introduction}

The ubiquitous presence of uncertainty in geometry, material properties, and loading conditions of a structure must be considered in the design process in order to achieve a robust and reliable performance. 
Most commonly, in reliability-based design optimization, a probabilistic failure criterion estimated by first- or second-order Taylor series expansion \cite{haldar2000probability} is considered through a design constraint \cite{bae2002reliability,maute2003reliability,kharmanda2004reliability,jung2004reliability,moon2004reliability,kim2006reliability,mogami2006reliability,eom2011reliability}. 
On the other hand, in {design optimization under uncertainty}, the effect of uncertainty is taken into account by optimizing the mean value of the structural performance \cite{beyer2007robust,de2017efficient,Diwekar2020}. Often a contribution from the variability of the performance is added to the objective to generate designs that are less sensitive to uncertainty \cite{alvarez2005minimization,conti2009shape,guest2008structural,chen2010rtso,chen2011new,asadpoure2011robust,tootkaboni2012topology,maute2014touu,keshavarzzadeh2017topology,de2019topology}. 


Topology optimization (TO) considers how  one or more materials can be optimally placed within a design domain to achieve a desired mechanical performance while satisfying design constraints. TO has found applications in several fields, such as structural mechanics, fluid flow, optics, and acoustics. The readers are referred to Sigmund and Maute \cite{sigmund2013topology},  Deaton and Grandhi \cite{deaton2014survey}, and references therein for details. 
In most of these applications, the structure is designed exclusively at the macroscale, assuming homogeneous materials \cite{bendsoe1988generating,suzuki1991homogenization,diaz1992shape,xia2017recent}. 
However, an increasing number of studies consider the use of heterogeneous materials, \textit{i.e.}, materials with spatially varying properties, such as engineered composites. 
The design of these heterogeneous materials at the microscale was performed, for example, in Sigmund \cite{sigmund1994materials,sigmund1995tailoring}, Lipton and Stuebner \cite{lipton2007optimal}, No\"{e}l and Duysinx \cite{noel2017shape}, Collet et al. \cite{collet2018topology}, and Chatterjee et al. \cite{Chatterjee2021} to achieve prescribed effective properties at the macroscale. An overview of the design of layered microstructure is given in Eschenauer and Olhoff \cite{eschenauer2001topology}. 
Concurrent multiscale TO \cite{bendsoe1988generating,xia2016multiscale,xia2017recent} seeks to optimize the structure at the macro as well as in the microscale (see references in Xia and Breitkopf \cite{xia2017recent} for a comprehensive list). 
{In Rodrigues et al. \cite{rodrigues2002hierarchical} and Coelho et al. \cite{coelho2008hierarchical}, a hierarchical approach is used, where a microstructure is designed for each element in the finite element mesh, which is then used to estimate the macroscale objective. 
Schury et al. \cite{schury2012efficient} employed free material optimization method that designs the stiffness tensor values for each finite element with appropriate constraints from the microscale problem. 
Xia and Breitkopf \cite{xia2014concurrent} solved the concurrent design problem at macro and for every element in the finite element mesh in the microscale by using computational homogenization that estimates the macroscale responses due to microscale inhomogeneity. Sivapuram et al. \cite{sivapuram2016simultaneous} considered the microstructure to remain the same inside a sub-region in the structure and used linearization to decompose the concurrent design problem. The connection between different adjacent microstructures was addressed in Du et al. \cite{du2018connecting}.} 

These highly optimized structures are, however, sensitive to defects \cite{pasini2019imperfect}, which can be introduced at the microscale during manufacturing. 
For example, melting process parameters, such as cooling rate in selective laser melting based additive manufacturing, can affect the pore and grain sizes \cite{beuth2001role,aboulkhair2014reducing,parry2016understanding}. This, in turn, affects the structural properties in the macroscale \cite{liu2019additive} and limits the application of TO 
\cite{dong2017survey,marmarelis2020data}. 
Hence, the design optimization process needs to incorporate the effects of these random defects at the microscale in order to produce robust structures. 

In this paper, we describe the variability in geometry at the microscale and impurities in the constituent materials of the microstructure using uncertainties characterized by known probability distributions. 
Hence, the TO design problem at the macroscale needs to address the uncertainties in geometry and material properties of the microstructure. 
In gradient-based approaches for solving the resulting TO problem, the evaluation of a large number of objective, constraints, and their gradients for many possible microstructure configurations may be needed using a standard Monte Carlo approach if variances of the gradients are large. As the number of possible microstructure scenarios {combining many realizations of the microstructure} can be extremely large, the optimization process becomes computationally burdensome. 
{To quantify the uncertainty, polynomial chaos expansion and its sparse variation \cite{ghanem2003stochastic,doostan2007stochastic,blatman2008sparse,doostan2011non} can be used to construct surrogate models of the performance, but the number of expansion coefficients and hence the number of objective, constraints, and gradient evaluations rapidly increases as the stochastic dimension of the problem increases. Perturbation methods, such as Taylor series expansion, with respect to the uncertain parameters can be used efficiently to estimate many gradients, but the expansion acuracy deteriorates for nonlinear objectives and constraints.}

To alleviate the computational burden of gradient-based approaches to solve TO under microscale uncertainty, we herein propose a method based on stochastic estimates of these quantities. 
We construct stochastic approximations of the objective, constraints, and gradients, using only a handful of random configurations of the microstructure {generated independently at every} design optimization iteration. To solve the optimization problem, these stochastic estimates are then used with Adam, a popular variant of the stochastic gradient descent, and globally convergent method of moving asymptotes (GCMMA). 
{To the best of our knowledge, this is the first approach that can address TO of structures in the presence of high-dimensional microstructure uncertainty.} 
We illustrate the efficacy of the proposed approach using two- and three-dimensional structures with two types for microstructural materials: (a) randomly dispersed inclusions in host matrix and (b) chopped fiber composites. 
The microstructural properties, such as the shape, size, and distribution of the inclusions, elastic moduli of the fiber and matrix, as well as the orientation of the fibers, are assumed uncertain and design independent. 
The results from these numerical examples show that the proposed approach produces an average design that has {computational cost only a small factor larger the cost of the corresponding deterministic TO.} 


The rest of the paper is organized as follows: In the next section, we define the TO problem under microscale uncertainty and briefly discuss the level set method used to describe the geometry of the structure at the macroscale. 
The subsequent section discusses the use of stochastic gradients with random microstructures and homogenization techniques to solve the TO problem. We illustrate the proposed approach using three numerical examples in Section \ref{sec:ex} before concluding the paper with a discussion on the future directions of this approach. 

\section{Optimization Problem}
\label{sec:methodology} 
\begin{figure}[!htb]
    \centering
    \begin{tikzpicture}[scale=1]  
    \node[inner sep=0pt] (3dmat) at (4,-0.05)
    {\includegraphics[width=.09\textwidth]{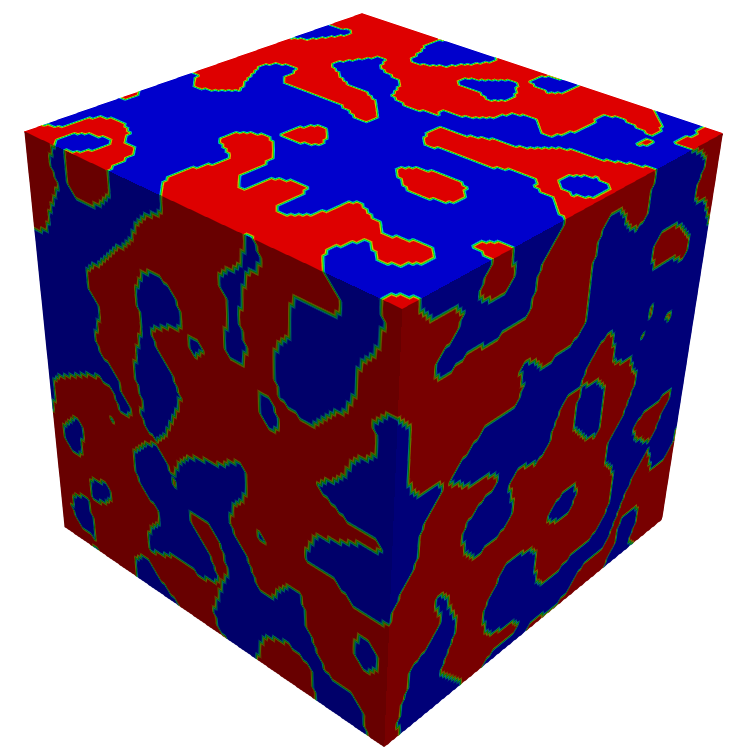}};
    \fill[red!30!black,opacity=0.6] (-3,-1,0) -- (-1,0,0) -- (-1,2,0) -- (-3,1,0) -- (-3,-1,0);
\pgfmathsetseed{3}
\shade [ball color=green!30] plot [smooth cycle, samples=8,domain={1:8}]
     (\x*360/8+5*rnd:1cm+2cm*rnd);
\draw[fill=white,draw = none, thick] (1,0.5) circle (0.5 cm);
\draw[fill=white,draw = none, thick] (0,-1.25) ellipse (0.75cm and 0.4cm);
\draw[very thick,-{Cone[pitch=30]}, color=black!20!red](1.8,1,0)--(2.7,1.4,0);
\draw[very thick,-{Cone[pitch=30]}, color=black!20!red](1.8,1.2,0)--(2.7,1.6,0);
\draw[very thick,-{Cone[pitch=30]}, color=black!20!red](1.8,1.4,0)--(2.7,1.8,0);
\draw[very thick,-{Cone[pitch=30]}, color=black!20!red](1.8,1.6,0)--(2.7,2.0,0);
\draw[draw=gray!50,very thick] (1.6,-0.5) circle (0.1cm); 
\draw[draw=gray!50,very thick] (4,0) circle (1cm); 
\draw[draw=gray!50,very thick] (1.6,-0.4) -- (3.35,0.75); 
\draw[draw=gray!50,very thick] (1.6,-0.6) -- (3.8,-0.975); 
\node[color=white] at (1.55,1.2) {$\Gamma_t$}; 
\node[color=black] at (-2.2,0.8) {$\Gamma_u$}; 
\node[rotate = 30] at (-2.2,1.7) {$u=\bar{u}$}; 
\node[] at (3.4,1.8) {$t=\bar{t}$}; 
\node[] at (4.1,-1.4) {Microstructure}; 
\node[] at (-0.5,0) {$\Omega$}; 
\node[color=white] at (1.5,-0.8) {$A$}; 
\end{tikzpicture}
    \caption{A schematic showing the structure at the macroscale with design domain $\Omega$, Dirichlet boundary condition $u=\bar{u}$ applied at the boundary $\Gamma_u$, and Neumann boundary condition $t=\bar{t}$ applied at the boundary $\Gamma_t$. A representative volume element (RVE) of the microstructure at point $A$ is shown as the inset figure. }
    \label{fig:prob_schem}
\end{figure}
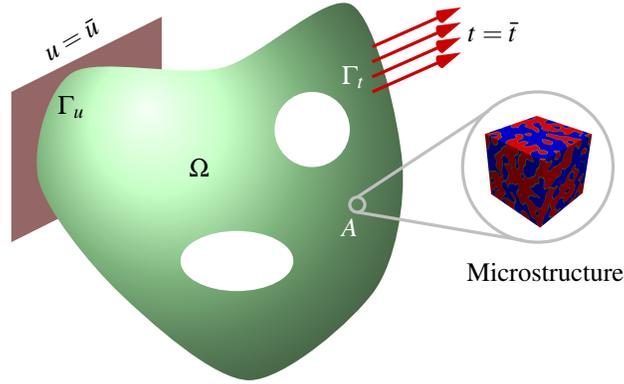

In this paper, we seek to optimize the structural shape and topology at the macroscale. The structure is made of a composite material whose effective properties are stochastic due to uncertainty at the microscale. We assume spatially varying uncertainty in geometry and material properties. 
Figure \ref{fig:prob_schem} shows a schematic of the problem, where the macrostructure is optimized over a design domain $\Omega$. 
Dirichlet boundary condition $u=\bar{u}$ is applied at the boundary $\Gamma_u$, and Neumann boundary condition $t=\bar{t}$ is applied at the boundary $\Gamma_t$. 
The figure also shows one realization of the representative volume element (RVE) of the microstructure at macroscopic point A. 
In this section, we formulate the optimization problem and then discuss the analysis model used to solve this problem. 


\subsection{Problem Formulation}

In macroscale optimization of structures under microscale uncertainty, the design optimization is performed at the macroscale while accounting for uncertainties in the microstructures. The cost function $f(\ppm;\xii): \mathbb{R}^{\np} \times \mathbb{R}^{\ny} \rightarrow \mathbb{R}$ and the constraint $\gm(\ppm;\xii): \mathbb{R}^{\np} \times \mathbb{R}^{\ny} \rightarrow \mathbb{R}^{\ngg}$ depend on the macroscale optimization variables $\ppm\in \mathbb{R}^{\np}$ as well as on the random variables $\xii\in \mathbb{R}^{\ny}$ with known probability distributions associated with the microstructure uncertainty. 
Accordingly, {the design optimization problem under uncertainty using average values of the cost function and constraints} 
can be defined 
as 
\begin{equation} \label{eq:opt_def}
\begin{split}
    &\mathop{\min~}\limits_{\ppm}R(\ppm):= \Exp_{\Ym}[f(\ppm;\Ym)]\\
    &\text{subject to } C_i(\ppm):= \Exp_{\Ym}\left[\gm(\ppm;\Ym)\right] \leq 0,\quad i=1,\dots,\ngg,
\end{split}
\end{equation}
where $\Exp_{\xii}[\cdot]$ denotes the expectation of its argument with respect to the probability distribution of $\Ym$; $R(\thetaa)$ is the objective and known as the \textit{expected risk}; and $C(\thetaa)$ is the \textit{expected constraint violation}. 
In the next subsection, we briefly discuss the model used to describe the geometry of the structure and estimate the sensitivity of the design with respect to the optimization variables. 

\subsection{Geometry and Analysis Model} \label{sec:level_set}

\begin{figure}[!htb]
    \centering
    \begin{subfigure}[b]{0.5\textwidth}
		\centering
\begin{tikzpicture}
    \begin{axis}[grid=none,view={20}{40},z buffer=sort, data cs=polar, axis line style={draw opacity=0},tick style={draw=none},yticklabels={,,},xticklabels={,,},zticklabels={,,}]
      \addplot3 [surf, domain=0:360, domain y=5:10,samples=60, samples y=20,shader=interp,colormap/jet]
      {-y+5};
      \addplot3 [data cs=cart,surf,domain=-10:10,samples=2, opacity=0.5,shader=interp,colormap/hot]
      {0};
      \addplot3 [domain=0:360, samples y=0, samples=50, ultra thick, z buffer=auto]
      (x,5.1,0);
      \addplot3 [surf,domain=0:360, domain y=0:5,samples=60, samples y=20,shader=interp,colormap/jet]
      {-y+5};
      \node[draw=none,color=white] at (0, -7)   (a) {$\phi(\xm)<0$};
      \node[draw=none] at (-4, 12)   (b) {$\phi(\xm)>0$};
      \node[draw=none] at (2, 15)   (c) {$\phi(\xm)=0$};

    \end{axis}
  \end{tikzpicture} 
  \caption{Level set field $\phi(\xm)$.}  
  \end{subfigure}
  \hfill
\begin{subfigure}[b]{0.45\textwidth}
		\centering
		\begin{tikzpicture}
  \filldraw [fill=gray, draw=black, thick] (0,0) rectangle (4,3);
  \filldraw [fill=white, draw=black,thick] (2,1.5) circle (0.75);
  \node[draw=none,color=white] at (0.75, 0.75)   (a) {$\Omega_I$};
  \node[draw=none] at (2, 1.5)   (b) {$\Omega_{II}$};
  \node[draw=none,color=white] at (3,1)   (c) {$\Gamma_{I,II}$};
  \end{tikzpicture}
  \vspace{20pt}
  \caption{Subdomains $\Omega_I$ and $\Omega_{II}$ with interface $\Gamma_{I,II}$.}
  \end{subfigure}
    \caption{An example of level set field $\phi(\xm)$ with corresponding subdomains $\Omega_I$ and $\Omega_{II}$, and interface $\Gamma_{I,II}$.}
    \label{fig:level_set}
\end{figure}
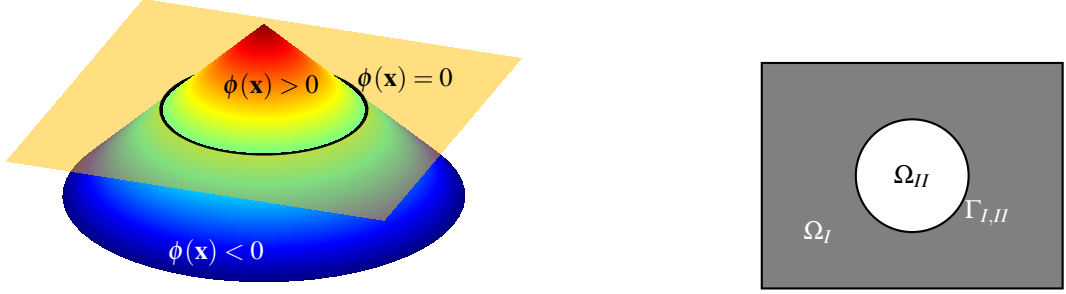
In this paper, we describe the geometry of the structures at the macroscale using a level set method following specifically the approach of Villanueva and Maute \cite{villanueva2014density}, and Sharma et al. \cite{sharma2017shape}. 
Considering a two-phase problem, the level set field $\phi(\xm)$ decomposes the design domain $\Omega_D$ into two distinct subdomains $\Omega_I$ and $\Omega_{II}$ as follows
\begin{equation}
\xm \in \begin{cases}
\Omega_I      \ \forall \ \phi(\xm)<0;\\
\Omega_{II}   \ \forall \ \phi(\xm)>0; \\
\Gamma_{I,II} \ \forall \ \phi(\xm)=0,
\end{cases}
\end{equation} 
where $\xm$ denotes the spatial coordinates and $\Gamma_{I,II}$ represents the interface between $\Omega_I$ and $\Omega_{II}$ as shown in Figure \ref{fig:level_set}. 
In this paper, the level set design field, $\phi(\xm)$, is discretized by
\begin{equation}
    \phi(\xm) = \sum_{i=1}^{N_n}\bar{\phi}_i \mathcal{N}_i(\xm), 
\end{equation}
where $\mathcal{N}_i(\xm)$ are bi-linear or tri-linear shape functions in 2D and 3D, respectively, and $N_n$ is the total number of nodes in the finite element mesh. 
Following Kreissl and Maute \cite{kreissl2012levelset}, we use a linear filter for the nodal level set values as follows
\begin{equation}
    \bar{\phi}_i =  \frac{\sum_{j=1}^{N_{rf}}w_{ij}\phi_i}{\sum_{j=1}^{N_{rf}}w_{ij}},
\end{equation}
where $N_{rf}$ is the number of nodes within the filter radius $r_f$; the weights are defined as $w_{ij}=\max(0,r_f-|\xm_i-\xm_j|)$; and $\phi_i$ is the level set function value at node $i$ before filtering. 
These unfiltered nodal level set values, $\phi_i$, are used as the design parameters. 
This approach allows for solving
the optimization problem with the help of mathematical programming methods.

We use the eXtended Finite Element Method (XFEM) to discretize the governing equations and use a generalized Heaviside enrichment strategy to approximate the displacement field in the solid domain. We use the unsymmetric Nitsche’s method to weakly implement the Dirichlet boundary conditions \cite{nitsche1971variational,burman2012fictitious}. The XFEM formulation is stabilized by the face-oriented ghost penalty method \cite{burman2014fictitious,schott2015face} to avoid ill-conditioning arising from small intersections of elements. 
Structural springs are also added to the disconnected solid subdomains to prevent rigid body motion \cite{geiss2019combined}. 
The shape sensitivities of the cost function $f(\thetaa;\xii)$ and the constraints $g_i(\thetaa;\xii),i=1,\dots,n_{\gm}$, with respect to the optimization variables are computed by the adjoint method as detailed in Sharma et al. \cite{sharma2017shape}. We use these gradients to determine the search direction during the optimization process, as described in the next section.

\section{Use of Stochastic Gradients} \label{sec:proposed_method}
In the standard Monte Carlo approach, $R(\ppm)$ and $\Cm(\ppm)$ in \eqref{eq:opt_def}, and their gradients with respect to the design optimization variables $\thetaa$ are estimated using $N_s$, \textit{e.g.}, $N_s\sim\mathcal{O}(10^{3})$, evaluations of the structural response for a given vector of optimization variables $\thetaa$ as follows 
\begin{equation}\label{eq:exp_risk}
    \begin{split}
        \widehat{R}(\ppm) & = \frac{1}{N_s}\sum_{i=1}^{N_s}f(\ppm;\xii_i);\\
        \widehat{C}_j(\ppm) & = \frac{1}{N_s}\sum_{i=1}^{N_s}g_j(\ppm;\xii_i);
        \quad j = 1,\dots,\ngg;\\
        { \widehat{\nabla R}(\ppm)} & = {  \frac{1}{N_s}\sum_{i=1}^{N_s}\nabla f(\ppm;\xii_i);}\\
        { \widehat{\nabla C_j}(\ppm)} & = {  \frac{1}{N_s}\sum_{i=1}^{N_s}\nabla g_j(\ppm;\xii_i); \quad j = 1,\dots,\ngg.}\\
    \end{split}
\end{equation}
Here, a hat notation is used for an estimate of a quantity. We do not explicitly state the dependence of these estimates on $\xii$. Note that, in general, to achieve a small estimation error, $N_s$ needs to be large, leading to a high computational cost. 
{Motivated by the success of the stochastic gradient descent method and its different variants for solving nonlinear non-convex optimization problems, \textit{e.g.}, in deep learning and macroscale topology optimization \cite{de2019topology,de2019bifidelity,de2021reliability,li2020momentum}, we propose a stochastic gradient based approach to alleviate the computational burden of topology optimization under microscale uncertainty.}
In this approach, instead of calculating the objective, constraints, and their gradients for large number of  microstructures realizations at every integration point, {a small number of random samples $n_s\ll N_s$, \textit{e.g.}, $n_s\sim\mathcal{O}(1)$, are randomly chosen to give small-sample unbiased estimates of the mean values in \eqref{eq:exp_risk}} as 
\begin{equation}\label{eq:grad}
    \begin{split}
    {\widehat{R}^{(n_s)}(\ppm)} & ={ \frac{1}{n_s}\sum_{i=1}^{n_s}f(\ppm;\xii_i);}\\
        {\widehat{C}_j^{(n_s)}(\ppm)} & ={ \frac{1}{n_s}\sum_{i=1}^{n_s}g_j(\ppm;\xii_i);}
        \quad j = 1,\dots,\ngg;\\
        \widehat{\nabla R}^{(n_s)}(\ppm) & = \frac{1}{n_s}\sum_{i=1}^{n_s}\nabla f(\ppm;\xii_i);\\
        \widehat{\nabla C_j}^{(n_s)}(\ppm) & = \frac{1}{n_s}\sum_{i=1}^{n_s}\nabla g_j(\ppm;\xii_i); \quad j = 1,\dots,\ngg,\\
    \end{split} 
\end{equation} 
{where the superscript $(n_s)$ is used to specify that the approximation uses $n_s\ll N_s$ number of random samples. Note that at every iteration different $n_s$ number of independent and identically distributed random samples are used to evaluate \eqref{eq:grad}.}
These coarse approximations of objective, constraints, and their design sensitivities are used in optimization algorithms, which are described next. 
In particular, following De et al. \cite{de2019topology}, two algorithms are investigated in this paper, namely the Globally Convergent Method of Moving Asymptotes (GCMMA) \cite{Svanberg:02} popular in TO, and Adaptive Moment (Adam) \cite{kingma2014adam}, a variant of the stochastic gradient descent algorithm. 




\subsection{Stochastic Gradient Descent (SGD) Method} 
\label{sec:sgd}
In the standard SGD method \cite{bottou2018optimization}, a single realization of $\xii$ is used at every optimization iteration to estimate the gradients in \eqref{eq:grad}. Here, a straightforward extension of the standard SGD method using a small batch of $n_s\geq 1$ random samples to estimate the gradients in \eqref{eq:grad} is used in the numerical examples, {which reduces the variance of the gradients used in each iteration resulting in a faster convergence}. This version is known as the \textit{mini-batch gradient descent} \cite{ruder2016overview,bottou2018optimization}. The parameters are updated at $k$th iteration as follows
%
\begin{equation}\label{eq:sgd}
\begin{split}
&\hb_k=\widehat{\nabla{R}}^{(n_s)}(\ppm_k) + \sum_{j=1}^{\ngg}{\kappa_j}{\widehat{\nabla C_j^+}}^{(n_s)}(\ppm_k);\\
&\ppm_{k+1}=\ppm_k-\eta \hb_{k},\\
\end{split}
\end{equation}
where $\eta $ is the step size, also known as the \textit{learning rate}; $\kappa_j,j=1,\dots,\ngg$ are positive hyperparameters to enforce the constraints. 
As SGD methods are restricted to unconstrained problem, we use a penalty formulation to account for the inequality constraints in \eqref{eq:opt_def} and define $C_j^+(\ppm)$ as 
\begin{equation}
    C_j^+(\ppm) = \Exp_{\xii}\left[\left(g_j^+(\ppm;\xii)\right)^2\right]; \quad j=1,\dots,\ngg, 
\end{equation} 
where $g_j^+(\ppm;\xii) = 0$ for $g_j(\ppm;\xii)\leq 0$ and $g_j^+(\ppm;\xii) = g_j(\ppm;\xii)$ otherwise. 
The computational cost of SGD is small {compared to that of standard Monte Carlo methods}. However, the convergence of the standard SGD method can be slow since 
the descent is only achieved in expectation {(\textit{i.e.}, the expectation of the stochastic gradients is the same as the gradients of objective and constraints in \eqref{eq:exp_risk})}. 
%


Recently, several variants of the standard SGD method with improved convergence have been proposed for training of neural networks \cite{ruder2016overview}. 
In this paper, one such variant, namely, the Adaptive Moment (Adam) \cite{kingma2014adam} is used in the numerical examples. 
This algorithm seeks to reduce the variability in $\hb$ over the iterations by accumulating historical gradient and squared gradient information using two exponential decay rates, $\beta_m$ and $\beta_v$. 
At $k$th iteration, the gradients are updated as follows 
%
\begin{equation}\label{eq:adamdecay}
\begin{split}
\mm_k &= \beta_m \mm_{k-1} + (1-\beta_m)\hb_k;\\
v_{k,j} &= \beta_v v_{k-1,j} + (1-\beta_v)h_{k,j}^2,\qquad j=1,\dots,\np.\\
\end{split}
\end{equation}
We use $\beta_m = 0.9$ and $\beta_v = 0.999$ herein as suggested in Kingma and Ba \cite{kingma2014adam}. 
An initialization bias correction is applied to $\mm_k$ and $v_{k,j}$ as follows
%
\begin{equation}\label{eq:adambias}
\begin{split}
\widehat{\mm}_k &= \frac{\mm_k}{1-\beta_m^k};\\
\widehat{v}_{k,j} &= \frac{v_{k,j}}{1-\beta_v^k}, \quad j=1,\dots,\np.\\
\end{split}
\end{equation}
Using these quantities, the parameters are updated as
\begin{equation}\label{eq:adam}
\theta_{k+1,j} = \theta_{k,j} - \eta\frac{\widehat{m}_{k,j}}{\sqrt{\widehat{v}_{k,j}}+{\epsilon}}\qquad j=1,2,\dots,\np.
\end{equation}
Algorithm \ref{alg:adam} summarizes these steps. 
%
\begin{algorithm}[htb]
	\begin{algorithmic}
		\STATE Given $\eta$, $\beta_m$, $\beta_v$, and $\epsilon$.
		\STATE Initialize $\ppm_1$.
		\STATE Initialize $\mm = \bm{0}$.
		\STATE Initialize $\vm = \bm{0}$.
		\FOR {$k=1,2,\dots,$}
		\STATE Compute $\hb_k := \hb(\ppm_{k})$.
		\STATE Set $\mm \leftarrow \beta_m\mm + (1-\beta_m)\hb_k$. {[see Eqn. (\ref{eq:adamdecay})]}
		\STATE Set $v_j \leftarrow \beta_vv_j + (1-\beta_v)h_{k,j}^2\qquad j=1,2,\dots,\np$. {[see Eqn. (\ref{eq:adamdecay})]}
		\STATE Set $\widehat\mm \leftarrow \mm/(1-\beta_m^k)$. {[see Eqn. (\ref{eq:adambias})]}
		\STATE Set $\widehat\vm \leftarrow \vm/(1-\beta_v^k)$. {[see Eqn. (\ref{eq:adambias})]}
		\STATE Set $\theta_{k+1,j} \leftarrow \theta_{k,j} - \eta\frac{\widehat{m}_j}{\sqrt{\widehat{v}_j}+{\epsilon}}\qquad j=1,2,\dots,\np$. {[see Eqn. (\ref{eq:adam})]}
		\ENDFOR
	\end{algorithmic}
	\caption{\textit{Adam} \cite{kingma2014adam}}
	\label{alg:adam}
\end{algorithm} 

\subsection{Globally Convergent Method of Moving Asymptotes (GCMMA) with Stochastic Gradients}\label{sec:gcmma} 
We also study GCMMA \cite{Svanberg:02} with stochastic gradients in this paper. In this algorithm, conservative approximations of the objective and constraint functions around the current design are used to formulate $n_{\thetaa}$ decoupled convex subproblems, which are solved by a primal-dual solution strategy.  
This algorithm gained popularity in solving TO problems with a large number of design parameters as the subproblems are separable, i.e., the primal problem can be decomposed in $n_{\thetaa}$ single-variable problems. {In this paper, we use GCMMA with no inner iteration and shape sensitivities are computed using the adjoint method \cite{sharma2017shape}.}

\subsection{Generation of Random Microstructures} \label{sec:micro_gen}

At every iteration of the algorithms described above, we use $n_s$ number of microstructural combinations per iteration. Each of these combinations assign randomly generated microstructures  to each element in the finite element discretization of the structure. 
{Note that, without loss of generality, we assume the microstructures at all integration points inside an element are the same. Next, we randomly assign these microstructures in each element. Another way to assign the microstructures is to assign them using spatial correlation. However, this is only an extension of the currently used random assignment. Also, with random assignment, the microstructure properties can vary significantly in the neighboring elements compared to the correlated assignment.} 
In this subsection, we describe how these random microstructures are generated and utilized in the numerical examples in Section \ref{sec:ex}. 

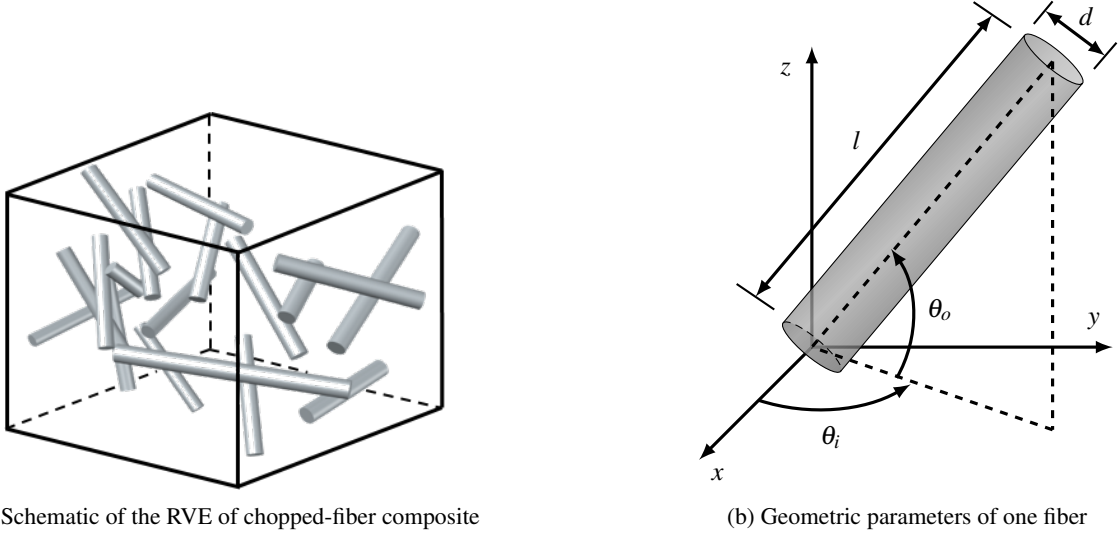
\begin{figure}[!htb]
    \centering
    \begin{subfigure}[t]{0.45\textwidth}
		\centering
    \includegraphics{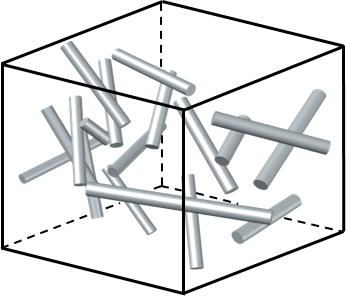} 
    \caption{Schematic of the RVE of chopped-fiber composite} \label{fig:chopped}
\end{subfigure}
  \hfill
\begin{subfigure}[t]{0.45\textwidth}
		\centering
\begin{tikzpicture}    
    \draw[-latex,very thick] (0,0) -- (0,4); 
    \draw[-latex,very thick] (0,0) -- (4,0); 
    \draw[-latex,very thick] (0,0) -- (-1.5,-1.5); 
\begin{scope}[rotate=-40]
\fill[top color=gray!50!black,bottom color=gray!10,middle color=gray,shading angle=45,opacity=0.25] (0,0) circle (0.5cm and 0.125cm);
\fill[left color=gray!30!black,right color=gray!30!black,middle color=gray!1,shading angle=45,opacity=0.5] (0.5,0) -- (0.5,5) arc (360:180:0.5cm and 0.125cm) -- (-0.5,0) arc (180:360:0.5cm and 0.125cm);
\fill[top color=gray!90!,bottom color=gray!2,middle color=gray!30,shading angle=45,opacity=0.25] (0,5) circle (0.5cm and 0.125cm);
\draw (-0.5,5) -- (-0.5,0) arc (180:360:0.5cm and 0.125cm) -- (0.5,5) ++ (-0.5,0) circle (0.5cm and 0.125cm);
\draw[densely dashed] (-0.5,0) arc (180:0:0.5cm and 0.125cm); 
\end{scope} 
\draw[dashed,very thick] (3.2,3.8) -- (3.2,-1.1); 
\draw[dashed,very thick] (3.2,-1.1) -- (0,0); 
\draw[dashed,very thick] (3.2,3.8) -- (0,0); 
\draw[-latex,>=stealth',very thick] (-135:1.0cm) arc[radius=2.5, start angle=252, end angle=300]; 
\draw[-latex,>=stealth',very thick] (-20:1.2cm) arc[radius=1.6, start angle=-30, end angle=35]; 
\node[] at (0.25,-1.2) {$\theta_i$}; 
\node[] at (1.7,0.5) {$\theta_o$}; 
\node[] at (-1.25,-1.7) {$x$}; 
\node[] at (3.75,0.35) {$y$}; 
\node[] at (-0.35,3.7) {$z$}; 
\draw[thick] (-0.5,0.5) -- (-1,0.84); 
\draw[thick] (2.65,4.25) -- (2.15,4.59); 
\draw[latex-latex,very thick] (-0.75,0.67) -- (2.4,4.42); 
\node[] at (0.6,2.75) {$l$}; 

\draw[thick] (2.95,4.3) -- (3.2,4.6);
\draw[thick] (3.75,3.7) -- (4.05,4.05);
\draw[latex-latex, very thick] (3.1,4.475) -- (3.9,3.875); 
\node[] at (3.65,4.4) {$d$};

\end{tikzpicture}
\caption{Geometric parameters of one fiber} \label{fig:fiber}
\end{subfigure}
    \caption{A schematic of the chopped-fiber composite that consists of multiple short fibers suspended in a matrix is shown in (a). The uncertain geometric parameters of these fibers, namely the length $l$, diameter $d$, in-plane angle $\theta_i$, and out-of-plane angle $\theta_o$ are shown in (b). }
    \label{fig:fib_comp} 
\end{figure}

We investigate two types of microstructures --- chopped fiber composite and randomly distributed two-phase material. In the chopped fiber composite, short fibers are suspended in a matrix as shown in Figure \ref{fig:chopped}. These fibers are parameterized using the aspect ratio $l/d$, in-plane angle $\theta_i$, and out-of-plane angle $\theta_o$; $l$ and $d$ are length and diameter of the fiber, respectively (see Figure \ref{fig:fiber}). These geometric parameters of the fiber, as well as the elastic moduli of the fiber and matrix, are assumed uncertain with known probability distributions as specified in Section \ref{sec:ex}. 

The second class of microstructures, {two-phase material with a stiff phase and a compliant phase}, is generated following the procedure in Roberts and Teubner \cite{roberts1995transport}. A $T$-periodic zero-mean Gaussian random field is considered as follows
\begin{equation} \label{eq:random_field}
    \psi(\rr) = \sum_{l=-N}^N\sum_{m=-N}^N\sum_{n=-N}^N c_{l,m,n}e^{\hat{i}\fm_{l,m,n}. \rr},
\end{equation}
where $ \rr$ is the position vector; $\Exp[\psi(\rr)]=0$; $\hat{i}=\sqrt{-1}$; and the vector $\fm_{l,m,n} = \frac{2\pi}{T}\left( l\ee_x + m\ee_y + n\ee_z \right)$ with $\ee_x$, $\ee_y$, and $\ee_z$ being the unit vectors along $x$, $y$, and $z$ directions, respectively. {Here, a cube with $2N\times2N\times2N$ resolution is used to simulate an RVE of the microstructure, which fixes the maximum wavenumber $T$ since $N=KT/2\pi$. Hence, the expansion in \eqref{eq:random_field} is truncated at $f_{l,m,n}=|\fm_{l,m,n}|\geq K$. 
Note that for a fixed $N$, the period $T$ can be increased to generate RVE of microstructure with higher levels of dispersion and vice versa.} 
To obtain strictly real-valued $\psi(\rr)$, we require $c_{l,m,n} = c^*_{-l,-m,-n}$, where $(\cdot)^*$ denotes the complex conjugate. We set $c_{0,0,0} = 0$ as the Gaussian random field is zero-mean. 
\begin{figure}[!htb]
    \centering
    \begin{subfigure}[t]{0.32\textwidth}
		\centering
    \includegraphics[scale=0.24]{Figures/3d_mat_v2_1.png} 
    \caption{Level-cut of the Gaussian random field at zero} \label{fig:3dmat1}
\end{subfigure}
  \hfill
\begin{subfigure}[t]{0.32\textwidth}
			\centering
    \includegraphics[scale=0.24]{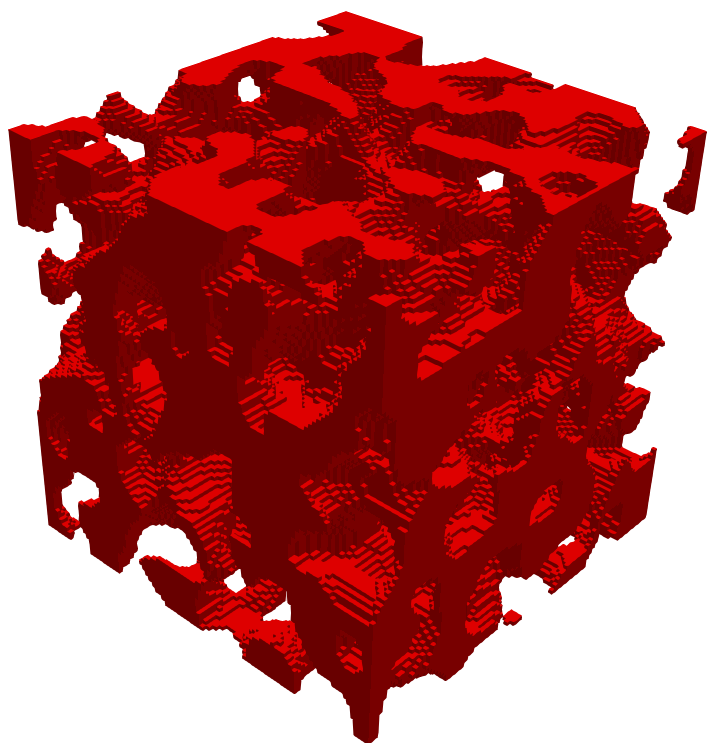} 
    \caption{Stiff phase (red) without the compliant phase (blue)} \label{fig:3dmat2}
\end{subfigure}
\hfill
\begin{subfigure}[t]{0.32\textwidth}
			\centering
    \includegraphics[scale=0.24]{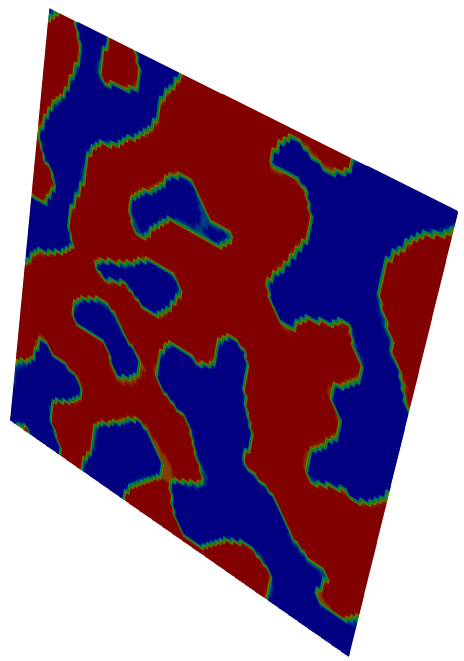} 
    \caption{A slice of the level-cut showing the stiff phase (red) and the compliant phase (blue)} \label{fig:3dmat3}
\end{subfigure}
    \caption{Randomly distributed red phase floating in the blue phase generated using a level-cut with a Heaviside function at zero of a $T=4\pi$ periodic Gaussian random field with maximum wavenumber $K=25$ and an exponential correlation function (following Roberts and Teubner \cite{roberts1995transport}) is shown in (a). In (b), only the red phase is shown by removing the blue phase. A slice of the cube, shown in (c), is used for two-dimensional problems. } 
    \label{fig:random_mat}
\end{figure}
Further, we assume $c_{l,m,n} = a_{l,m,n}+ib_{l,m,n}$, where $a_{l,m,n}$ and $b_{l,m,n}$ are independent and identically distributed zero-mean Gaussian random variables with variance 
\begin{equation} \label{eq:ab} 
    \Exp[a^2_{l,m,n}] = \Exp[b^2_{l,m,n}] = \frac{1}{2} {S}_{\psi\psi,K}(f_{l,m,n})\left(\frac{2\pi}{T}\right)^3. 
\end{equation}
{Here, $S_{\psi\psi,K}(f)$ is the spectral density of the random field $\psi(\rr)$ at $f$ (\textit{i.e.}, Fourier transform of the correlation function of $\psi(\rr)$ at $f$) truncated at $f_{l,m,n}\geq K$, which is related to the non-truncated spectral density $S_{\psi\psi}(f)$ of the random field as follows }
\begin{equation}
    S_{\psi\psi,K}(f) = \frac{S_{\psi\psi}(f)}{\bigintss_0^K 4\pi f^2 S_{\psi\psi}(f) \mathrm{d}f}.
\end{equation}
We assume the correlation function of the random field and the corresponding spectral density are respectively given by 
\begin{equation}
    \begin{split}
        \Exp\left[ \psi(\rr_1)\psi(\rr_2) \right] &= e^{-\lVert \rr_1-\rr_2 \rVert^2_2};\\
        {S_{\psi\psi}(f)} &{= \frac{e^{-{f^2}/{4}}}{(4\pi)^{3/2}}},\\
    \end{split}
\end{equation}
where $\lVert\cdot\rVert_2$ denotes the $\ell_2$-norm of its argument. {Using \eqref{eq:ab}, realizations of $a_{l,m,n}$ and $b_{l,m,n}$ are generated first, which are then used to generate realizations of the Gaussian random field $\psi(\rr)$.} Once we generate the random field, we perform a level-cut at zero using the Heaviside function {assigning $\psi(\rr)>0$ to the stiff phase (red) and $\psi(\rr)\leq0$ to the compliant phase (blue)}, which gives a three-dimensional RVE of the microstructure as shown in Figure \ref{fig:random_mat} for $T=4\pi$ and $K=25$. For two-dimensional problems, we use a slice through the 3D RVE {at the left end} (see Figure \ref{fig:3dmat3}). 

\subsection{Homogenization Methods} 
Homogenization methods are used to estimate the effective material properties at the macroscale considering the microstructural layout. 
In the following, we consider a finite element setting, where the weak form of the governing equations are integrated at the macroscale. 

For linear elastic problems, the material constitutive relation at the macroscale position $\xm$ is given by
\begin{equation}\label{eq:hom_eqn}
    \overline{\sigmaa}(\xm) = \Chom:\overline{\epsl}(\xm),
\end{equation}
where $\overline{\sigmaa}(\cdot)$ is the macroscopic stress tensor; $\overline{\epsl}(\cdot)$ is the macroscopic strain tensor; and $\Chom$ is the fourth order homogenized constitutive tensor. 
The Mori-Tanaka homogenization method \cite{mori1973average,benveniste1987new} is used in this paper for chopped-fiber composite and is discussed next. For random two-phase materials, we use a first-order computational homogenization method, which is also briefly discussed in this section. 

\begin{figure}[!htb]
    \centering
    \begin{tikzpicture}[scale=1]  
    \fill[red!30!black,opacity=0.6] (-3,-1,0) -- (-1,0,0) -- (-1,2,0) -- (-3,1,0) -- (-3,-1,0);
\pgfmathsetseed{3}
\shade [ball color=green!30] plot [smooth cycle, samples=8,domain={1:8}]
     (\x*360/8+5*rnd:1cm+2cm*rnd);
\draw[fill=white,draw = none, thick] (1,0.5) circle (0.5 cm);
\draw[fill=white,draw = none, thick] (0,-1.25) ellipse (0.75cm and 0.4cm);
\draw[very thick,-{Cone[pitch=30]}, color=black!30!red](1.8,1,0)--(2.7,1.4,0);
\draw[very thick,-{Cone[pitch=30]}, color=black!30!red](1.8,1.2,0)--(2.7,1.6,0);
\draw[very thick,-{Cone[pitch=30]}, color=black!30!red](1.8,1.4,0)--(2.7,1.8,0);
\draw[very thick,-{Cone[pitch=30]}, color=black!30!red](1.8,1.6,0)--(2.7,2.0,0);
\node[color=white] at (1.55,1.2) {$\Gamma_t$}; 
\node[color=black] at (-2.2,0.8) {$\Gamma_u$}; 
\node[rotate = 30] at (-2.2,1.7) {$u=\bar{u}$}; 
\node[] at (3.4,1.8) {$t=\bar{t}$}; 
\node[] at (-0.5,0) {$\Omega$}; 

\fill[green!50,opacity=0.75,draw=black,xshift=-1cm,yshift=-1cm,thick] (-0.3,-0.1,0) -- (-0.1,-0.1,0) -- (-0.1,0.1,0) -- (-0.3,0.1,0) -- (-0.3,-0.1,0); 
\draw[thick,black] (-4.75,-1.175) -- (-1.3,-0.9); 
\draw[thick,black] (-3.5,-2.8) -- (-1.1,-1.1);


	\coordinate (P1) at (-12cm,0cm); 
	\coordinate (P2) at (0cm,0cm); 

	\coordinate (A1) at (-5,-2cm); 
	\coordinate (A2) at (-5,-4cm); 

	\coordinate (A3) at ($(P1)!.8!(A2)$); 
	\coordinate (A4) at ($(P1)!.8!(A1)$);

	\coordinate (A7) at ($(P2)!.7!(A2)$);
	\coordinate (A8) at ($(P2)!.7!(A1)$);

	\coordinate (A5) at
	  (intersection cs: first line={(A8) -- (P1)},
			    second line={(A4) -- (P2)});
	\coordinate (A6) at
	  (intersection cs: first line={(A7) -- (P1)}, 
			    second line={(A3) -- (P2)});



	\fill[green!35] (A2) -- (A3) -- (A6) -- (A7) -- cycle; 
	
	\fill[green!25] (A3) -- (A4) -- (A5) -- (A6) -- cycle; 
	
	\fill[green!10] (A5) -- (A6) -- (A7) -- (A8) -- cycle; 
	
	\draw[thick,dashed] (A5) -- (A6);
	\draw[thick,dashed] (A3) -- (A6);
	\draw[thick,dashed] (A7) -- (A6);


	\fill[gray!50,opacity=0.2] (A1) -- (A2) -- (A3) -- (A4) -- cycle; 
	\fill[gray!90,opacity=0.2] (A1) -- (A4) -- (A5) -- (A8) -- cycle; 

	\draw[thick] (A1) -- (A2);
	\draw[thick] (A3) -- (A4);
	\draw[thick] (A7) -- (A8);
	\draw[thick] (A1) -- (A4);
	\draw[thick] (A1) -- (A8);
	\draw[thick] (A2) -- (A3);
	\draw[thick] (A2) -- (A7);
	\draw[thick] (A4) -- (A5);
	\draw[thick] (A8) -- (A5);
	
	

\node[inner sep=0pt] (3dmat) at (-1,-9.25)
    {\includegraphics[width=.09\textwidth]{Figures/3d_mat_v2_1.png}}; 
\node[inner sep=0pt] (3dmat2) at (-2,-5)
    {\includegraphics[width=.09\textwidth]{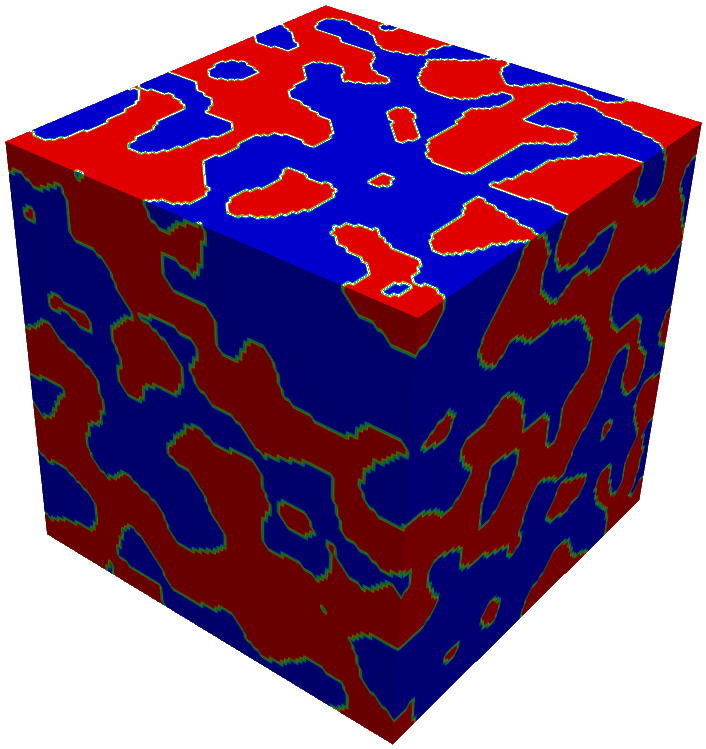}}; 
    \node[inner sep=0pt] (3dmat2) at (-1.5,-6.5)
    {\includegraphics[width=.09\textwidth]{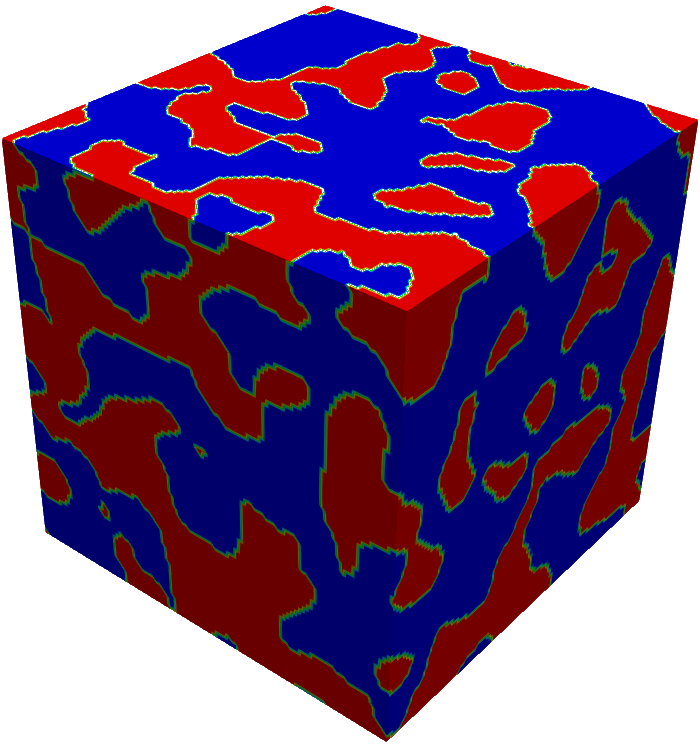}}; 
\node[rotate = 20,scale=3] at (-1.5,-7.65) {$\vdots$}; 

\node[inner sep=0pt] (3dmat) at (-8,-9.25)
    {\includegraphics[width=.09\textwidth]{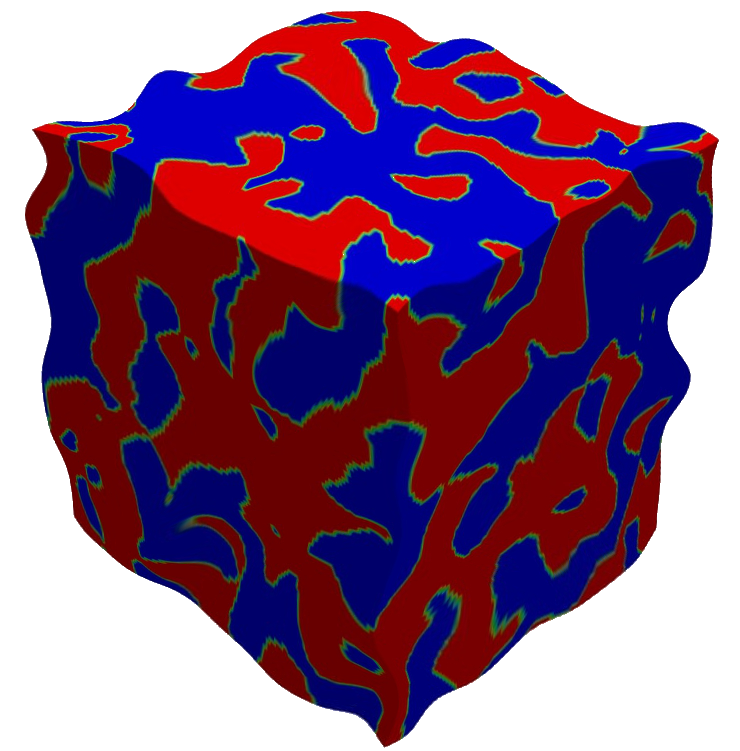}}; 
\node[inner sep=0pt] (3dmat2) at (-7,-5)
    {\includegraphics[width=.09\textwidth]{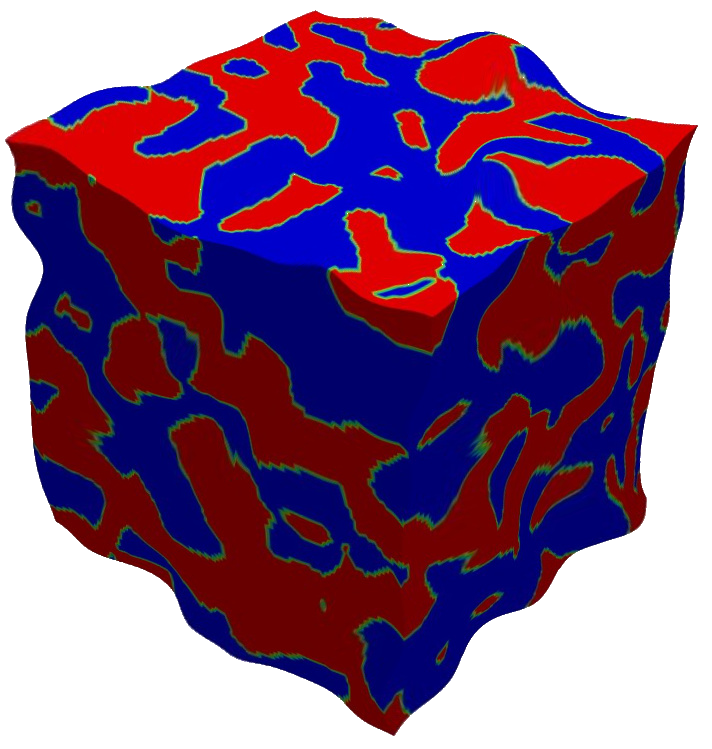}}; 
    \node[inner sep=0pt] (3dmat2) at (-7.5,-6.5)
    {\includegraphics[width=.09\textwidth]{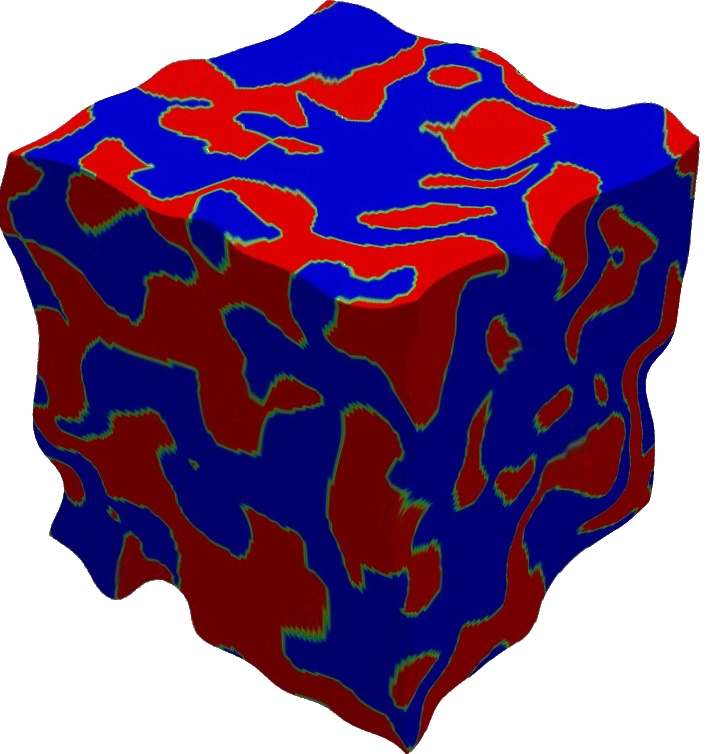}}; 
\node[rotate = -20,scale=3] at (-7.5,-7.65) {$\vdots$};


\shade[ball color = gray!20!red, opacity = 0.75] (-5.2,-3.4) circle (0.075cm); 
\shade[ball color = gray!20!red, opacity = 0.75] (-6,-3) circle (0.075cm); 
\shade[ball color = gray!20!red, opacity = 0.75] (-4.2,-2.8) circle (0.075cm); \shade[ball color = gray!20!red, opacity = 0.75] (-4.8,-2.5) circle (0.075cm); 
\shade[ball color = gray!20!red, opacity = 0.75] (-4.2,-1.8) circle (0.075cm); \shade[ball color = gray!20!red, opacity = 0.75] (-5.2,-2.4) circle (0.075cm); 
\shade[ball color = gray!20!red, opacity = 0.75] (-6,-2) circle (0.075cm); 
\shade[ball color = gray!20!red, opacity = 0.75] (-4.8,-1.5) circle (0.075cm); 

\draw[-{Triangle[width=12pt,length=10pt]}, line width=5pt,color=gray] (-4.1,-2.9) -- (-2.5,-4.1); 
\draw[-{Triangle[width=12pt,length=10pt]}, line width=5pt,color=gray] (-6.5,-4.1) -- (-4.3,-2.9); 
\draw[-{Triangle[width=12pt,length=10pt]}, line width=5pt,color=gray] (-3,-8) -- (-6.4,-8); 

\node[] at (-4.7,-7.7) {$(\sigmaa,\epsb)$}; 
\node[] at (-4.7,-8.3) {PDE solve}; 
\node[] at (-7.5,-3.5) {$\overline{\sigmaa}(\xm)= \langle \sigmaa(\xm,\xm^\prime) \rangle$}; 
\node[] at (-1.75,-3.5) {$\overline{\epsb}(\xm)= \langle \epsb(\xm,\xm^\prime) \rangle$}; 
\node[] at (-4.7,-0.7) {$\left(\overline{\sigmaa}(\xm),\overline{\epsb}(\xm)\right)$}; 

\draw[line width=3.5pt,color=gray,dotted] (-9.5,-4) -- (-7,-4); 

\node[rotate = 90] at (-9.25,-2.5) {Macroscale}; 
\node[rotate = 90] at (-9.25,-5.5) {Microscale}; 
\end{tikzpicture}
\caption{A schematic showing the use of homogenization method for multiple realizations of the random microstructure for an integration point in an element of the finite element mesh to estimate the macroscale average stress and strain. }
    \label{fig:fe2}
\end{figure}

\subsubsection{Mori-Tanaka Method}\label{sec:mori_tanaka}

{The Mori-Tanaka method \cite{mori1973average,benveniste1987new}, a mean-field homogenization method based on Eshelby's inclusion problem \cite{eshelby1957determination}, 
estimates the effective elasticity tensor for a fiber-reinforced composite material by \cite{benveniste1987new,zaoui2002continuum,bohm2004short,tran2018mori,budarapu2019multiscale}}
\begin{equation}
    \Chom = \Cbf+v_\mathrm{f}\langle(\Cbm-\Cbf):\Am\rangle:\left( (1-v_\mathrm{f})\mathbb{I} + v_\mathrm{f} \langle\Am\rangle \right)^{-1}. 
\end{equation}
Here, $v_\mathrm{f}$ is the fiber volume fraction of the composite material; $\Cbm$ and $\Cbf$ are fourth order elasticity tensors of the matrix and fiber, respectively; $\mathbb{I}$ is the fourth order identity matrix; $\Am$ is the strain concentration tensor that relates the strain $\epsl_f$ in the fiber with the strain $\epsl_\mathrm{m}$ in the matrix as $\epsl_\mathrm{f} = \Am : \epsl_\mathrm{m}$; and $\langle \cdot \rangle$ denotes the average over all fiber orientations. The stress concentration tensor $\Am$ is computed by 
\begin{equation}
    \Am = \left[ \mathbb{I}+\Sm:\Cbm^{-1}:(\Cbf-\Cbm) \right]^{-1},
\end{equation}
where $\Sm$ is the fourth order Eshelby's tensor \cite{eshelby1957determination}.

\subsubsection{Computational Homogenization Method} 

This approach assumes separation of scales and periodicity of the microstructure \cite{xia2016multiscale,xia2017recent}. 
The homogenized stress tensor and the macroscopic constitutive tensor are computed from the solutions of a boundary value problem defined over a representative volume element  for different boundary conditions.
Consider one such integration point with macroscale position $\xm$ and microscale position $\ym$. In this approach, the macroscale stress $\overline{\boldsymbol{\sigma}}(\mathbf{x})$ is related to the microscale  stress ${\boldsymbol{\sigma}}(\xm,\xm^\prime)$ over the RVE domain $\Omega_\mathrm{mi}$ as follows
\begin{equation}
\overline{{\sigmaa}}(\xm) = \langle \sigmaa(\xm,\xm^\prime) \rangle = \frac{1}{|\Omega_\mathrm{mi}|}\int_{\Omega_\mathrm{mi}}{{\sigmaa}}(\xm,\xm^\prime) \mathrm{d} \xm^\prime,
\end{equation}
where $\langle\cdot\rangle$ denotes the volume average over the RVE domain. 
The microscale stress ${\sigmaa}(\xm,\xm^\prime)$ is estimated from the boundary value problem associated with the RVE with a constraint on the macroscale strain $\overline{\epsl}(\xm)$ given by
\begin{equation}
\overline{\epsl}(\xm) = \langle \epsb(\xm,\xm^\prime) \rangle =  \frac{1}{|\Omega_\mathrm{mi}|}\int_{\Omega_\mathrm{mi}}{\epsl}(\xm,\xm^\prime) \mathrm{d} \xm^\prime,
\end{equation}
where ${\epsl}(\xm,\xm^\prime)$ is the microscale strain. 
{
To define the RVE boundary value problem that needs to be solved for each of these microstructures, let us consider the first-order local displacement field $\uu(\xm,\xm^\prime)$ and the strain field ${\epsl}(\xm,\xm^\prime)$ as, \cite{michel1999effective,xia2014concurrent},
\begin{equation} \label{eq:u_period}
\begin{split}
        \uu(\xm,\xm^\prime) &= \overline{\epsl}(\xm)\cdot\xm^\prime + \uu^*(\xm^\prime); \\
        {\epsl}(\xm,\xm^\prime) &= \overline{\epsl}(\xm) + {\epsl^*}(\xm^\prime),\\
\end{split}
\end{equation}
where $\uu^*(\xm^\prime)$ is periodic up to a rigid body motion and 
${\epsl^*}(\xm^\prime)$ with $\langle \epsl^*(\xm^\prime)\rangle=\mathbf{0}$ is due to the periodic displacement field $\uu^*(\xm^\prime)$. The local stress field $\sigmaa(\xm,\xm^\prime)$ is periodic as well, and the boundary value problem at the microscale becomes 
\begin{equation} \label{eq:micro_prob} 
    \begin{split}
        \sigmaa(\xm,\xm^\prime) &= \Cb(\xm^\prime) : \left( \overline{\epsl}(\xm) + {\epsl}^*(\xm^\prime) \right);\\
        \mathrm{div}\left(\sigmaa(\xm,\xm^\prime)\right) &= \mathbf{0} ~~\mathrm{in}~~\Omega_\mathrm{mi};\\
        \uu^*(\xm^\prime)&\text{ is periodic};\\
        \sigmaa(\xm,\xm^\prime)\cdot\mathbf{n}&\text{ is anti-periodic}.\\
        \end{split} 
\end{equation} 
Here, each of the phases in the microstructure is assumed linear elastic with constitutive tensor $\Cb(\xm^\prime)$; $\mathrm{div}(\cdot)$ denotes the divergence of its vector argument; the anti-periodicity imposes the condition that $\sigmaa(\xm,\xm^\prime)\cdot\mathbf{n}$ has opposite values at the opposite boundaries of $\Omega_\mathrm{mi}$, which is due to the periodicity of $\sigmaa(\xm,\xm^\prime)$ and the unit normal vector $\mathbf{n}$ being opposite on the opposite boundaries. The implementation of the periodicity in \eqref{eq:u_period} in the finite element approach is performed by specifying $\overline{\epsl}(\xm)$ and applying nodal displacement constraints \cite{michel1999effective,xia2014concurrent}. Once the problem in \eqref{eq:micro_prob} is solved for all six independent components of the strain, \eqref{eq:hom_eqn} can be used to estimate the homogenized constitutive tensor.} 

Figure \ref{fig:fe2} depicts a schematic showing that in the presence of uncertainty in the microscale, many realizations of the random microstructure at an integration point in an element are needed to estimate the macroscale average stress tensor for a given macroscopic strain state. 
{Instead, in the proposed stochastic gradient based approach, for every element in the finite element mesh, we consider one random realization of the uncertain microstructure to generate a microstructure layout configuration. A small number, $n_s\sim\mathcal{O}(1)$, of these microstructure layouts are then used at every iteration of the optimization process resulting in a computational cost that is only a few times larger than the cost of the corresponding deterministic optimization.} 

\section{Numerical examples} \label{sec:ex} 
In this section, we use three numerical examples to illustrate the utility of stochastic gradients in designing the macrostructure under microscale uncertainty. The first two examples consider the design of a two- and a three-dimensional beam, respectively. In the third example, we design a bracket to support a payload box. 
{
The geometry of the structures in these examples are described using the level set method and the governing equations are discretized using the XFEM approach as described in Section \ref{sec:level_set}. The common parameters used during analysis of these examples are listed in Table \ref{tab:ex_params}.  
\begin{table}[!htb] 
\centering
\begin{threeparttable}
\caption{{Parameters used in the description of geometry and analysis in the numerical examples. }} 
\centering 
\begin{tabular}{c | c } 
\hline 
\Tstrut Parameter & Value \\ [0.5ex] 
\hline  
\Tstrut Filter radius, $r_f$ & $1.6h$ (2D), $1.8h$ (3D)$^\dagger$ \\ 
Nitsche penalty factor & $100$ \\
Ghost penalty factor & $0.01$ \\
Spring stiffness factor & $10^{-6}$ \Bstrut \\  
\hline 
\end{tabular}
\label{tab:ex_params} 
\begin{tablenotes}
\item [$\dagger$] $h$ is the element size 
\end{tablenotes}
\end{threeparttable}
\end{table}
}

\subsection{Example I: Design of a Two-dimensional Beam}

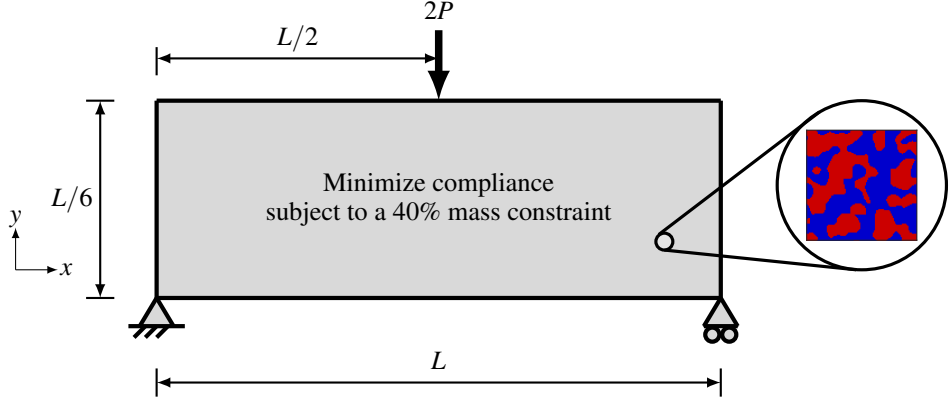
\begin{figure}[!htb]
    	\centering
	\begin{tikzpicture}[scale=0.75,every node/.style={minimum size=1cm},on grid]
	\draw [ultra thick,fill=gray!30,draw=none] (0,-3.5) rectangle (10,0);
	\draw [ultra thick] (0,0) -- (10,0);
	\draw [ultra thick] (10,0) -- (10,-3.5);
	\draw [ultra thick] (10,-3.5) -- (0,-3.5);
	\draw [ultra thick] (0,-3.5) -- (0,0);
	
	\draw[ultra thick,fill=gray!30] (9.85,-4.15) circle (0.125);
	\draw[ultra thick,fill=gray!30] (10.15,-4.15) circle (0.125);
	\draw[fill=gray!30]    (10,-3.5) -- ++(0.3,-0.5) -- ++(-0.6,0) -- ++(0.3,0.5);
	\draw [ultra thick] (10,-3.5) -- (10.3,-4);
	\draw [ultra thick] (10.3,-4) -- (9.7,-4);
	\draw [ultra thick] (9.7,-4) -- (10,-3.5);
	
	\draw[fill=gray!30]    (0,-3.5) -- ++(0.3,-0.5) -- ++(-0.6,0) -- ++(0.3,0.5);
	\draw [ultra thick] (0,-3.5) -- (0.3,-4);
	\draw [ultra thick] (0.5,-4) -- (-0.5,-4);
	\draw [ultra thick] (-0.3,-4) -- (0,-3.5);
	\draw [ultra thick] (0,-4) -- (-0.2,-4.2);
	\draw [ultra thick] (-0.2,-4) -- (-0.4,-4.2);
	\draw [ultra thick] (0.2,-4) -- (0,-4.2);
	
	\node[draw=none] at (5, -1.5)  (c)     {Minimize compliance};
	\node[draw=none] at (5, -2)  (c)     {subject to a 40\% mass constraint};
	
	\draw[-latex, line width=1mm] (5,1.25) -- (5,0);
	\node[draw = none] at (5,1.6) () {$2P$};
	
	
	
	\draw[thick,latex-latex] (0,-5) -- (10,-5);
	\node[draw = none] at (5,-4.6) () {$L$};
	\draw[thick] (0,-4.75) -- (0,-5.25);
	\draw[thick] (10,-4.75) -- (10,-5.25);
	
	\draw[thick,latex-latex] (0,0.75) -- (5,0.75);
	\node[draw = none] at (2.5,1.1) () {$L/2$};
	\draw[thick] (0,0.95) -- (0,0.45);
	\draw[thick] (5,0.95) -- (5,0.45);
	
	\draw[thick,latex-latex] (-1,0) -- (-1,-3.5);
	\node[draw = none] at (-1.5,-1.75) () {$L/6$};
	\draw[thick] (-1.25,0) -- (-0.75,0); 
	\draw[thick] (-1.25,-3.5) -- (-0.75,-3.5); 
	
	\draw[very thick,fill=white!90!black] (9,-2.5) circle (0.15cm); 
	\draw[very thick] (8.9,-2.375) -- (11.75,-0.2); 
	\draw[very thick] (9,-2.65) -- (12.4,-3); 
	\draw[very thick] (12.5,-1.5) circle (1.5cm); 
	\node[inner sep=0pt] (inset) at (12.5,-1.5)
    {\includegraphics[width=.09\textwidth]{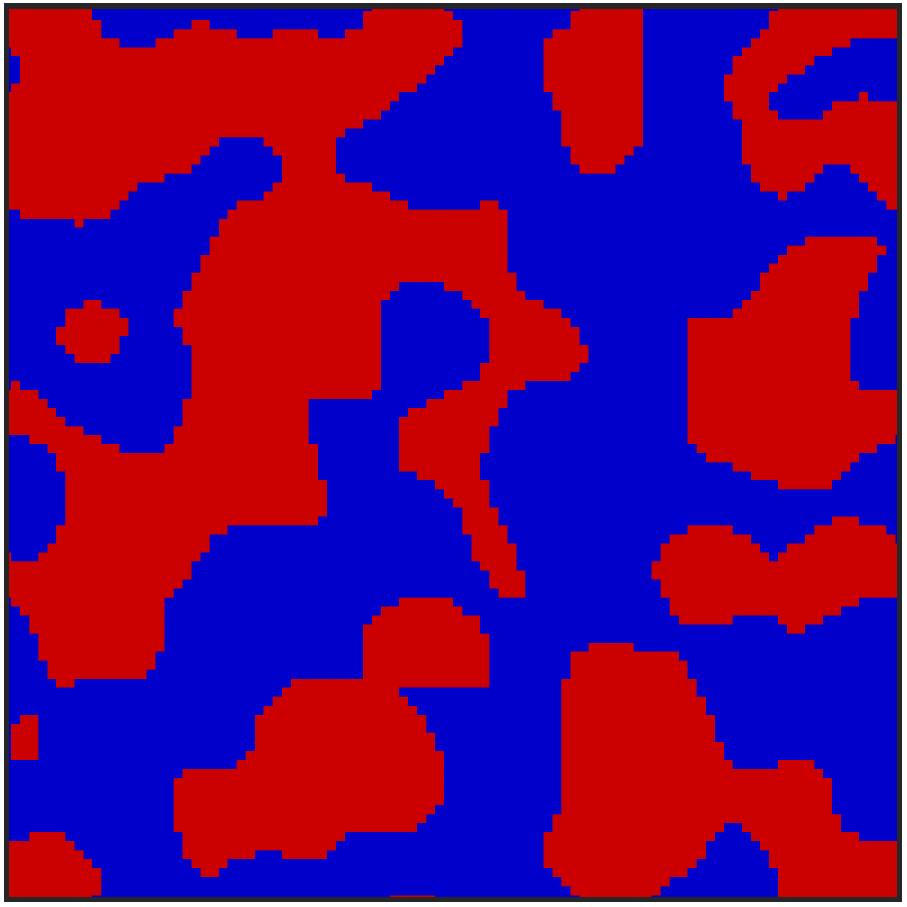}}; 
    \draw[-latex] (-2.5,-3) -- (-1.75,-3); 
    \draw[-latex] (-2.5,-3) -- (-2.5,-2.25); 
    \node[] at (-1.6,-3) {$x$}; 
    \node[] at (-2.5,-2.1) {$y$}; 
	\end{tikzpicture}
    \caption{The strain energy of a simply supported beam is minimized in Example I considering a random two-phase composite; a point load is applied at the center of the upper edge, and the optimization problem is subject to a 40\% mass constraint. } 
    \label{fig:mbb2d_prob}
\end{figure}


\subsubsection{Problem Description}
We first consider a two-dimensional simply supported beam subjected to a point load at the midspan (Figure \ref{fig:mbb2d_prob}). At the macroscale, the shape and topology of the beam is optimized to minimize the strain energy subject to a 40\% mass constraint. Using symmetry, only one-half of the domain is used for optimization. 
The optimization problem is defined as follows 
\begin{equation}
\begin{split}
    &\min_{\thetaa} ~~ R(\thetaa) = {w_{\!_\Psi}}\frac{\Exp_{\xii}\left[ {\Psi(\thetaa;\xii)} \right]}{\Psi_0} + w_{\mathrm{per}} P_\mathrm{per}(\thetaa) + w_{\mathrm{reg}} P_\mathrm{reg}(\thetaa);\\
    &\text{subject to } 
    C(\ppm) = g(\ppm) = 
    \frac{\int_{\Omega} \rho(\ppm)\mathrm{d}\xm}{\int_{\Omega}\mathrm{d}\xm}
    - \gamma_\mathrm{req}\leq 0,\\
    \end{split}
\end{equation}
where $w_{\!_\Psi}$ is the weight for the strain energy in the objective; $\Psi_0$ is the initial strain energy of the structure; {and $\Omega$ denotes design domain in the macroscale}. The two penalty terms, $P_\mathrm{per}$ and $P_\mathrm{reg}$, are defined following Barrera et al. \cite{barrera2020hole} as 
\begin{equation}
    \begin{split}
        P_\mathrm{per} &= \frac{\int_{\Gamma_{I,II}(\thetaa)}\mathrm{d}A}{\int_{\Gamma}\mathrm{d}A}; \quad
        P_\mathrm{reg} = w_\phi\frac{\int_{\Omega}(\phi(\xm;\thetaa)-\tilde{\phi}(\xm;\thetaa))^2 \mathrm{d}\xm}{\int_{\Omega}(\tilde{\phi}_{\max}-\tilde{\phi}_{\min})^2 \mathrm{d}\xm} + w_{\nabla \phi}\frac{\int_{\Omega}|\nabla\phi(\xm;\thetaa) - \nabla \tilde{\phi}(\xm;\thetaa)|^2\mathrm{d}\xm}{\int_\Omega \mathrm{d}\xm},\\
    \end{split}
\end{equation}
where $\Gamma_{I,II}(\thetaa)$ denotes the boundary between the structure and void {in the design domain $\Omega$}; $\Gamma$ is the boundary of the design domain $\Omega$; {$A$ denotes the surface area;} $\tilde{\phi}(\xm;\thetaa)$ is a signed distance field truncated at $\tilde{\phi}_{\min}$ and $\tilde{\phi}_{\max}$; and the parameters $w_\phi$ and $w_{\nabla\phi}$ are set at $1/\int_{\Gamma}\mathrm{d}A$ (see Barrera et al. \cite{barrera2020hole} for details). 
Note that we add the perimeter penalty to avoid the emergence of any geometric features with irregular shapes and add the level set regularization penalty to avoid any spurious oscillations in the level set field $\phi(\xm;\thetaa)$.  
The weights for the strain energy term and the penalty terms are set at $w_{\!_\Psi}=0.9$, $w_\mathrm{per}=0.025$, and $w_\mathrm{reg}=0.5$, respectively. In the constraint, we use $\gamma_\mathrm{reqd} = 0.40$, which is enforced with a penalty parameter $\kappa=1000$ in Adam as described in Section \ref{sec:sgd}. 

We describe the geometry of the structure at the macroscale using a level set and use XFEM for the solution of governing equations (see Section \ref{sec:level_set}). 
The finite element mesh of the half-domain has a total of $120\times40$ bilinear elements. {The initial geometry of the structure is shown in Figure \ref{fig:ExI_init_design}, where we add 18 holes in the global $x$-direction for each of the 6 rows in the global $y$-direction of the half domain. The level set field with these 18 holes is generated as $\phi(\xm)=\max_{i}\{\phi_i\}_{i=1}^{18}$, where 
\begin{equation}\label{eq:holes}
    \begin{split} 
        &\phi_i=1-\left(\frac{x_\mathrm{hole}}{r_\mathrm{hole}}\right)^{10}+\left(\frac{y_\mathrm{hole}}{r_\mathrm{hole}}\right)^{10}.\\ 
    \end{split}
\end{equation}
Here, $(x_\mathrm{hole},y_\mathrm{hole})$ is the local coordinate of the hole, and we use $r_\mathrm{hole}=1/15$. Note that a positive level set value represents void. 
We summarize the specifications used in this example in Table \ref{tab:ExI_specs}.}  

\begin{table}[!htb] 
\centering
\caption{Summary of specifications used to formulate and solve the optimization problem in Example I. The values are in consistent units. } \label{tab:ExI_specs} 
  \begin{tabular}{c|c|l|l}
    \hline
    \multicolumn{2}{c|}{\Tstrut Category} & Parameter & Value \Bstrut \\ \hline \multicolumn{2}{c|}{\Tstrut \multirow{4}{*}{Problem formulation}} & Weight for strain energy, $w_\Psi$ & 0.90 \\ 
    \multicolumn{2}{c|}{} & Weight for perimeter penalty, $w_\mathrm{per}$ & 0.025 \\
    \multicolumn{2}{c|}{} & Weight for regularization penalty, $w_\mathrm{reg}$ & 0.50 \\
    \multicolumn{2}{c|}{} & Point load, $P$ & 1.00 \\
    \multicolumn{2}{c|}{} & Mass constraint, $\gamma_\mathrm{reqd}$ & 0.40 \Bstrut \\ \hline 
    \multicolumn{2}{c|}{\Tstrut \multirow{3}{*}{Mesh (half domain)}} & Length, $L/2$ & 3.0 \\
    \multicolumn{2}{c|}{} & Height, $L/6$ & 1.0 \\ 
    \multicolumn{2}{c|}{} & Discretization & $120\times40$ \\ \hline 
    \Tstrut\multirow{6}{*}{Case Ia} & \multirow{5}{*}{Microstructure} & Period, $T$ & $4\pi$ \\ 
    & & Max. wavenumber, $K$ & 25 \\
    & & Elastic modulus (blue phase), $E_1$ & 1.0\\
    & & Elastic modulus (red phase), $E_2$ & 10.0\\
    & & Poisson's ratio (both phases), $\nu_1,\nu_2$ & 0.3\\ \cline{2-4}
    \Tstrut & Optimization & Step size, $\eta$ & 0.05 \\ \hline 
    \Tstrut\multirow{6}{*}{Case Ib} & \multirow{5}{*}{Microstructure} & Period, $T$ & $4\pi$ \\ 
    & & Max. wavenumber, $K$ & 25 \\
    & & Elastic modulus (blue phase), $E_1$ & 0.1\\
    & & Elastic modulus (red phase), $E_2$ & 10.0\\
    & & Poisson's ratio (both phases), $\nu_1,\nu_2$ & 0.3\\ \cline{2-4}
    \Tstrut & Optimization & Step size, $\eta$ & 0.05 \\ \hline 
    \Tstrut\multirow{6}{*}{Case IIa} & \multirow{5}{*}{Microstructure} & Period, $T$ & $2\pi$ \\ 
    & & Max. wavenumber, $K$ & 50 \\
    & & Elastic modulus (blue phase), $E_1$ & 1.0\\
    & & Elastic modulus (red phase), $E_2$ & 10.0\\
    & & Poisson's ratio (both phases), $\nu_1,\nu_2$ & 0.3\\ \cline{2-4}
    \Tstrut & Optimization & Step size, $\eta$ & 0.025 \\ \hline 
    \Tstrut\multirow{6}{*}{Case IIb} & \multirow{5}{*}{Microstructure} & Period, $T$ & $2\pi$ \\ 
    & & Max. wavenumber, $K$ & 50 \\
    & & Elastic modulus (blue phase), $E_1$ & 0.1\\
    & & Elastic modulus (red phase), $E_2$ & 10.0\\
    & & Poisson's ratio (both phases), $\nu_1,\nu_2$ & 0.3\\ \cline{2-4}
    \Tstrut & Optimization & Step size, $\eta$ & 0.025 \\ 
    \hline 
    \multicolumn{2}{c|}{\Tstrut \multirow{3}{*}{Solution strategy}} & No. of possible microstructures per elem. & 200 \\ 
    \multicolumn{2}{c|}{} & Random config. per iter., $n_s$ & 4 \\
    \multicolumn{2}{c|}{} & No. of optimization variables, $n_{\thetaa}$ & 5289 \\ 
    \multicolumn{2}{c|}{} & Penalty to implement $C(\thetaa)$, $\kappa$ & $1000$ \\ 
    \hline
  \end{tabular} 
\end{table}

\subsubsection{Microstructure Scenarios} 
In this example, each of the elements in the finite element mesh is assumed to have any microstructure from the $200$ possible scenarios. {Hence, the microstructure in each element of the finite element mesh can be thought of as a uniformly distributed discrete random variable with 200 possible realizations. 
These realizations are generated using the Gaussian random field model described in Section \ref{sec:micro_gen}. 
We use the finite element method to compute the first-order homogenized constitutive tensor for each of these random microstructures discretized in $100\times100$ mesh with bilinear elements. Note that each of the two-dimensional microstructures requires three linear finite element analyses to estimate three independent components of the strain and the homogenized constitutive tensor $\Chom$.} 

We study the design problem under two cases distinguished by the anisotropy in the two-phase microstructures. These cases and the corresponding observations are discussed below. In each of these cases, we approximate the objective, constraints, and their design sensitivities with four configurations of the microstructure layout per optimization iteration. In every configuration, each of the elements in the mesh is randomly assigned a microstructure from the 200 possible instances.

\begin{figure}[!htb]
    \centering
    \begin{subfigure}[b]{0.45\textwidth}
    \centering
    \includegraphics[scale=0.25]{Figures/T4pi_K25_1.pdf}~~~~
    \includegraphics[scale=0.248]{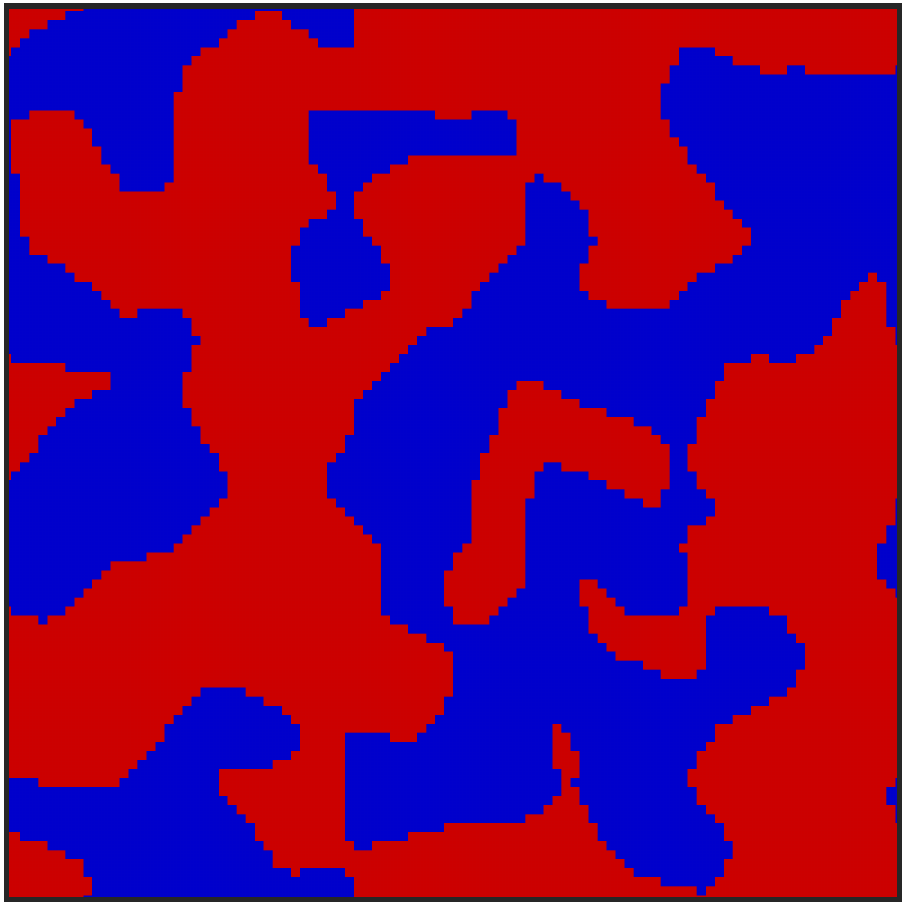}\\~\\
    \includegraphics[scale=0.25]{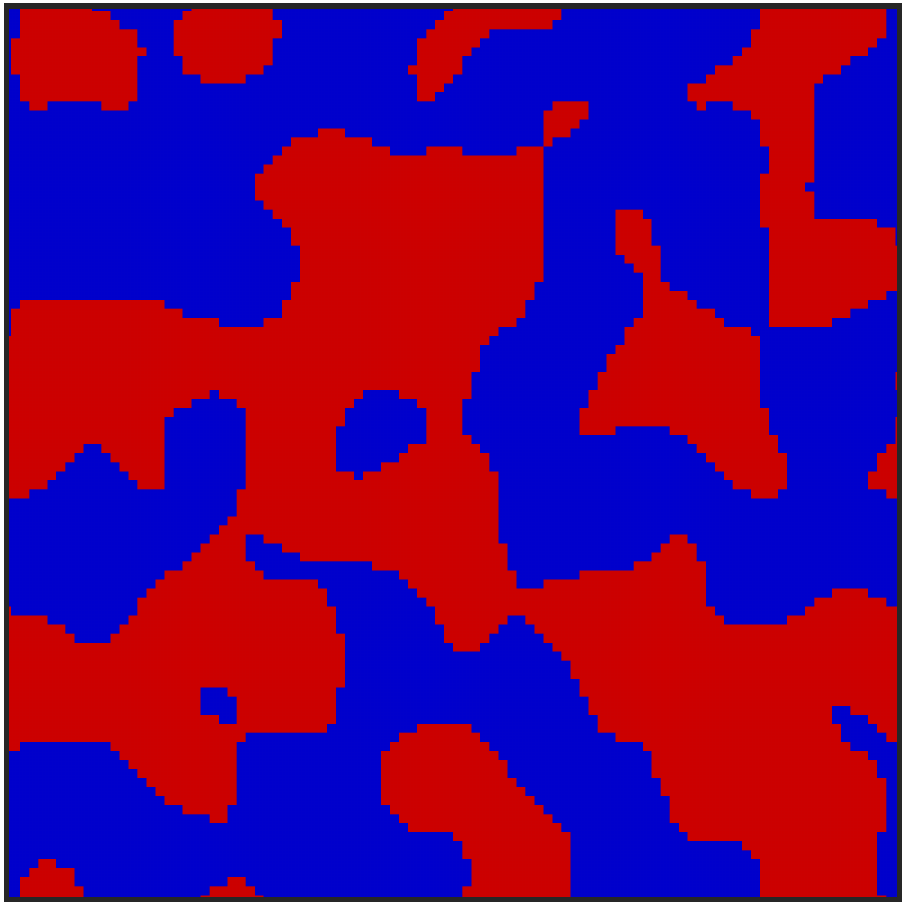}~~~~
    \includegraphics[scale=0.25]{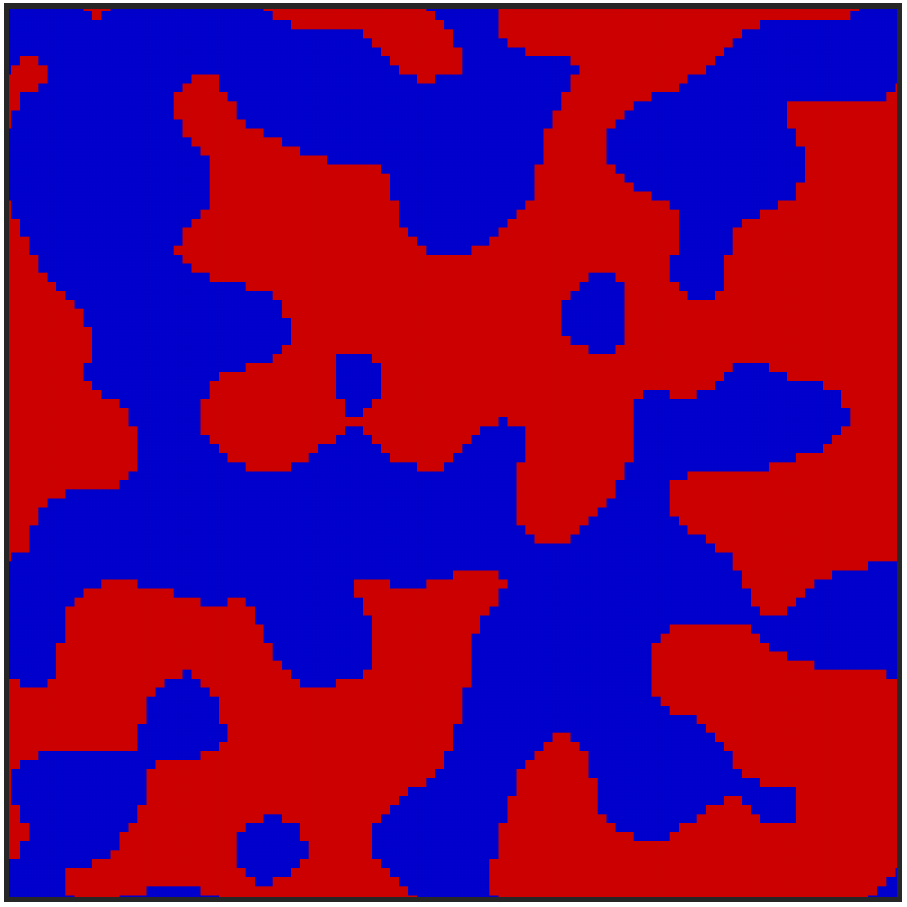}
    \caption{Case I: Four random microstructures generated using $T=4\pi$ and $K=25$. } \label{fig:exI_caseI_micro}
    \end{subfigure}~~~~
    \begin{subfigure}[b]{0.45\textwidth}
    \centering
    \includegraphics[scale=0.25]{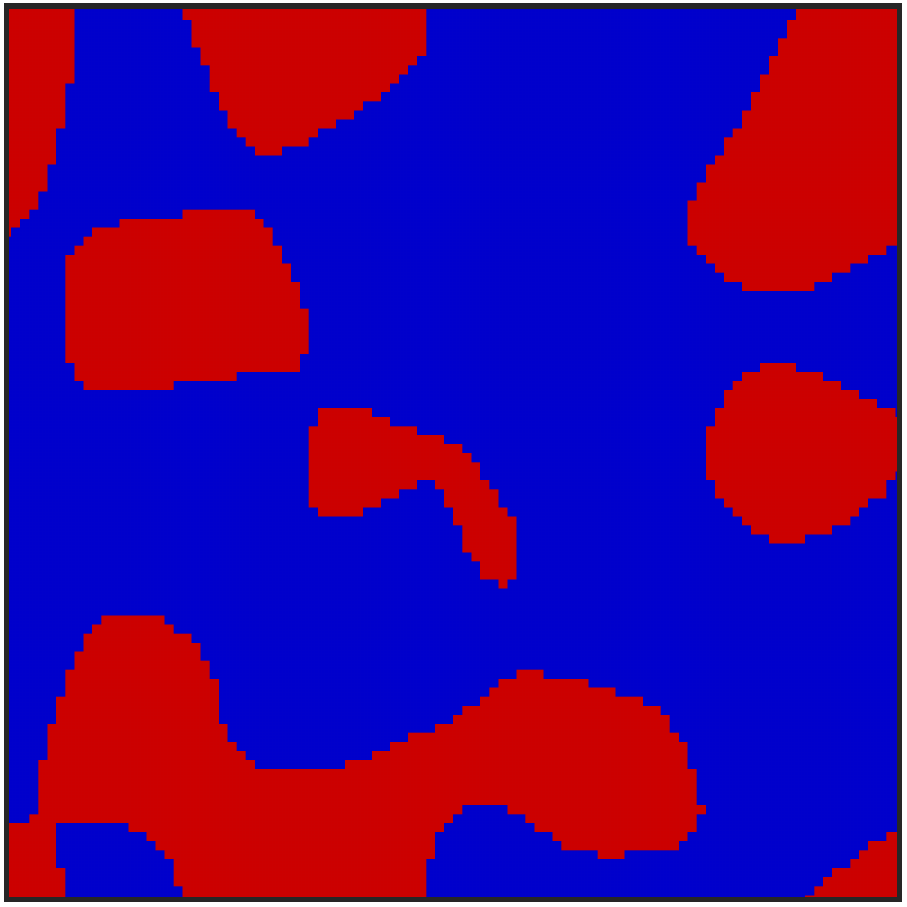}~~~~
    \includegraphics[scale=0.248]{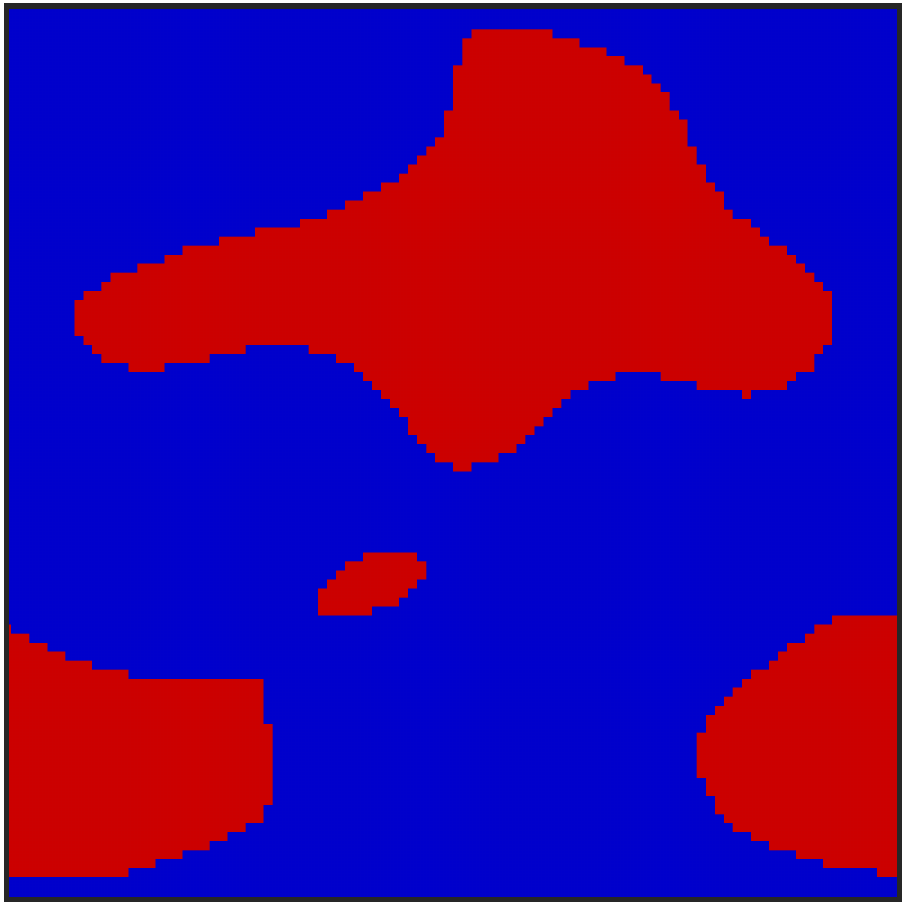}\\~\\
    \includegraphics[scale=0.25]{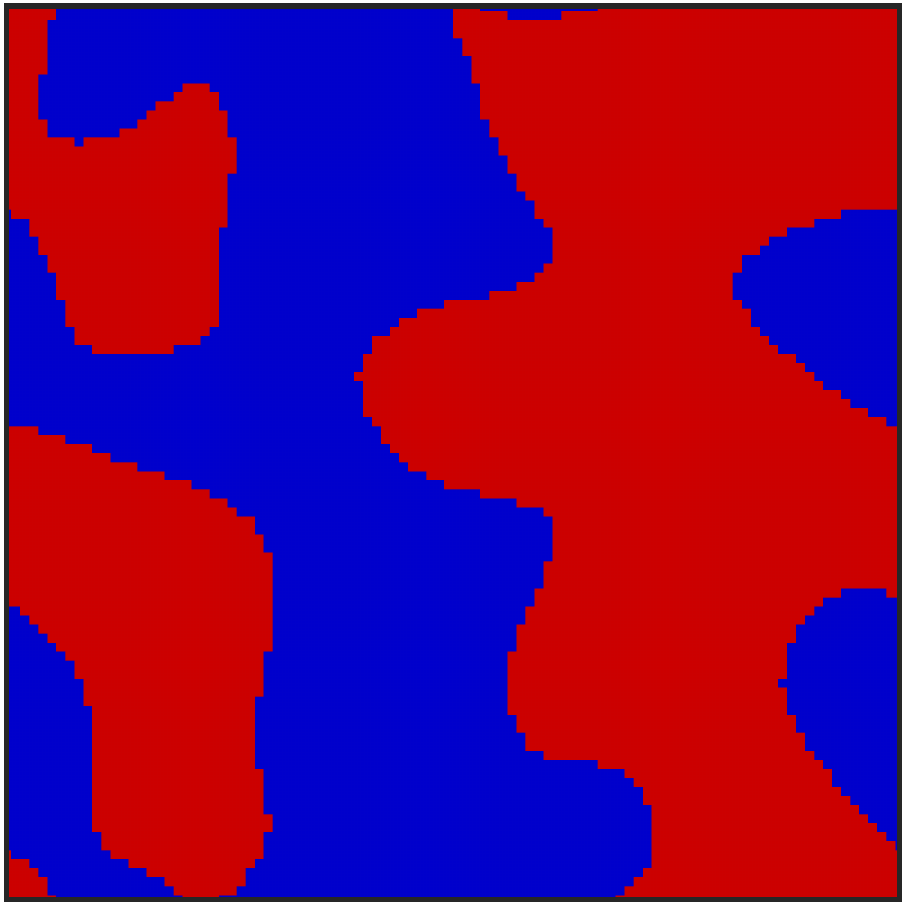}~~~~
    \includegraphics[scale=0.25]{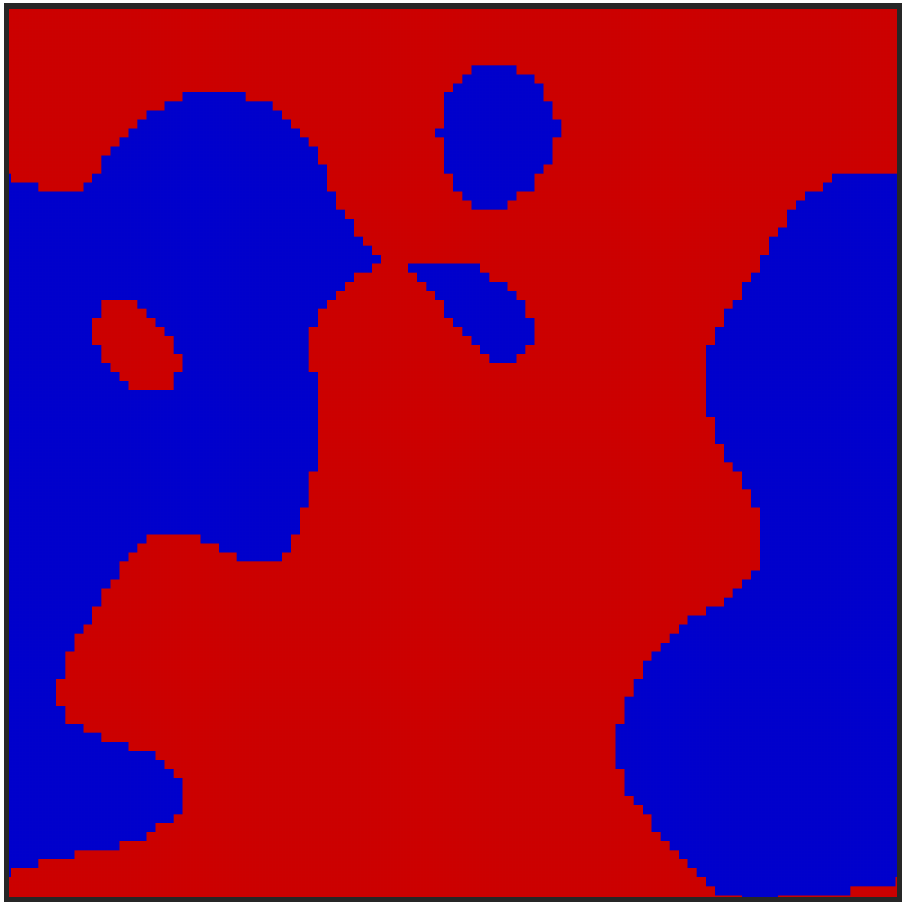}
    \caption{Case II: Four random microstructures generated using $T=2\pi$ and $K=50$. } \label{fig:exI_caseII_micro}
    \end{subfigure}
    \caption{Four random microstructures from each of the two cases studied with Example I. Note that we select $K$ such that $N=50$, 
    {which results in the same number of terms in \eqref{eq:random_field} and the same finite element discretization during the first-order computational homogenization for both of these cases.} }
    \label{fig:exI_micro} 
\end{figure}

\textbf{Case Ia: } In this case, we generate 200 random microstructures using $T=4\pi$ and $K=25$ in \eqref{eq:random_field} to give $N=50$, which results in a finely dispersed distribution of the two phases. We set the ratio of the elastic moduli of the phases to $E_1/E_2=10$, where $E_1$ is the elastic modulus of the stiff (red) phase and $E_2$ is the elastic modulus of the compliant (blue) phase. The fine dispersion and the chosen ratio of elastic moduli leads to a moderate level of anisotropy. Four such microstructures are shown in Figure \ref{fig:exI_caseI_micro}. 

\begin{figure}[!htb]
    \centering
    \begin{subfigure}[b]{\textwidth}
    \centering 
    \begin{tikzpicture}
        \node[inner sep=0pt] (structure) at (0,0){
    \includegraphics[scale=0.275]{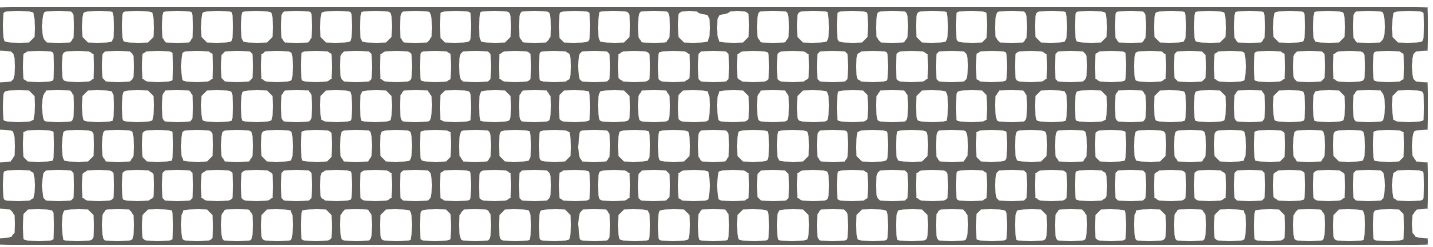}}; 
    \draw[-latex] (-6.35,-0.9) -- (-5.6,-0.9); 
    \draw[-latex] (-6.35,-0.9) -- (-6.35,-0.15); 
    \node[] at (-5.45,-0.9) {$x$}; 
    \node[] at (-6.35,0) {$y$}; 
    \end{tikzpicture} 
    \caption{Initial design} \label{fig:ExI_init_design}
    \end{subfigure}
    \begin{subfigure}[b]{0.7\textwidth}
    \begin{tikzpicture}
        \node[inner sep=0pt] (structure) at (0,0){\includegraphics[scale=0.275]{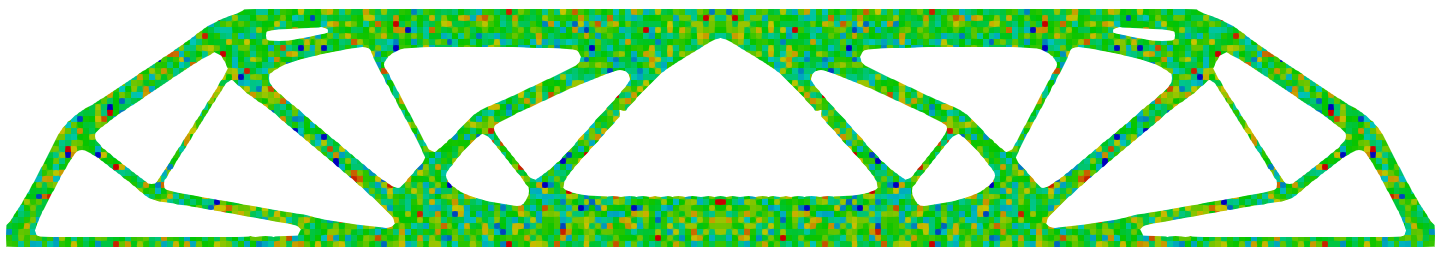}}; 
        \node[inner sep=0pt] (colorbar) at (5.5,0)
    {\includegraphics[scale=0.375]{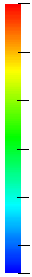}}; 
    \node at (5.75,1.8) {$C_{1111}/C_{2222}$}; 
    \node at (5.95,1.35) {$1.3$}; 
    \node at (5.95,0.4) {$1.1$}; 
    \node at (5.95,-0.55) {$0.9$}; 
    \node at (6,-1.3) {$0.74$}; 
    \draw[-latex] (-6.35,-0.9) -- (-5.6,-0.9); 
    \draw[-latex] (-6.35,-0.9) -- (-6.35,-0.15); 
    \node[] at (-5.45,-0.9) {$x$}; 
    \node[] at (-6.35,0) {$y$}; 
        \end{tikzpicture}
        \caption{Adam design}
    \end{subfigure}\\
    \begin{subfigure}[b]{0.7\textwidth}
    \begin{tikzpicture}
        \node[inner sep=0pt] (structure) at (0,0){\includegraphics[scale=0.275]{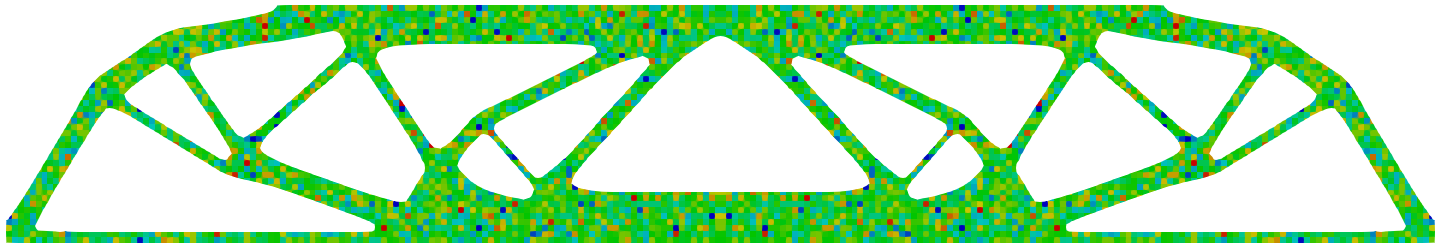}};
        \node[inner sep=0pt] (colorbar) at (5.5,0)
    {\includegraphics[scale=0.375]{Figures/colorbar_T4pi_E10.png}}; 
    \node at (5.75,1.8) {$C_{1111}/C_{2222}$}; 
    \node at (5.95,1.35) {$1.3$}; 
    \node at (5.95,0.4) {$1.1$}; 
    \node at (5.95,-0.55) {$0.9$}; 
    \node at (6,-1.3) {$0.74$}; 
    \draw[-latex] (-6.35,-0.9) -- (-5.6,-0.9); 
    \draw[-latex] (-6.35,-0.9) -- (-6.35,-0.15); 
    \node[] at (-5.45,-0.9) {$x$}; 
    \node[] at (-6.35,0) {$y$}; 
        \end{tikzpicture}
        \caption{GCMMA design}
    \end{subfigure}\\
    \begin{subfigure}[b]{0.55\textwidth}
    \centering
    \includegraphics[scale=0.275]{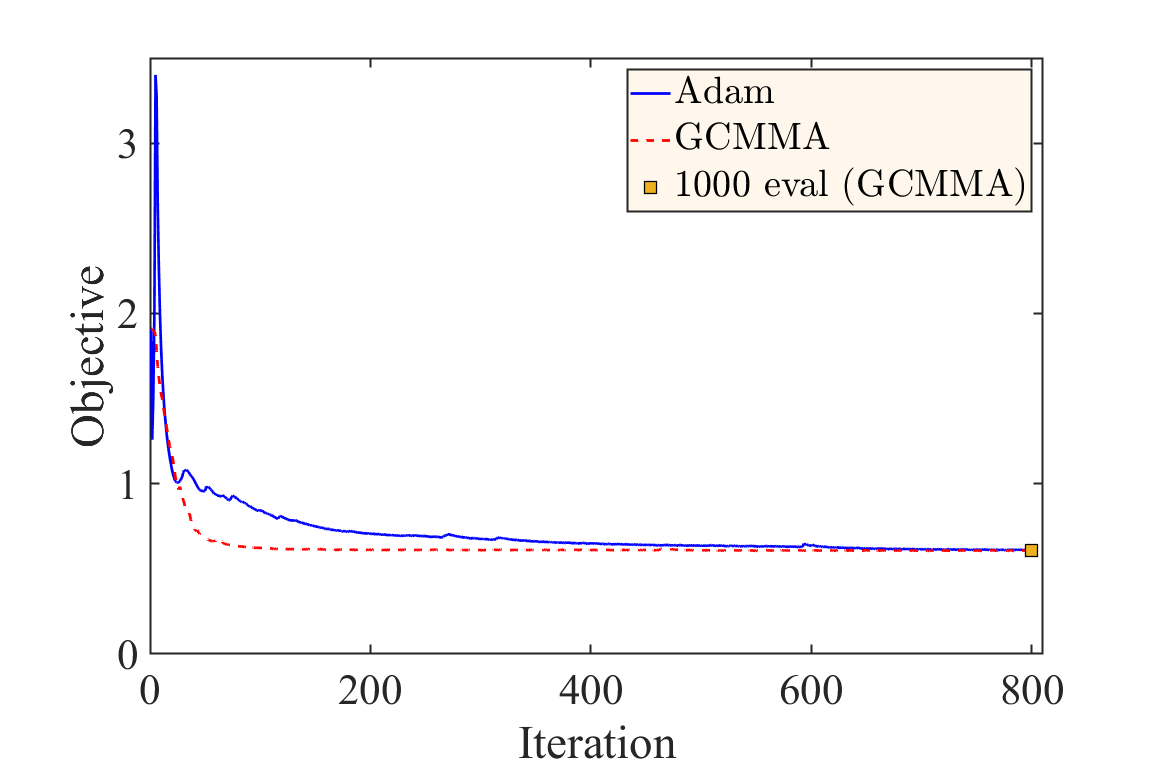}
    \caption{Objective vs. iteration} \label{fig:exI_caseI_Eratio10_obj}
    \end{subfigure}
    \caption{Comparison of the initial design (a) and designs obtained from Adam (b) and GCMMA (c) for Case Ia of Example I using $E_1/E_2=10$. The color shading of the elements in (b) and (c) corresponds to the ratio $C_{1111}/C_{2222}$  for one random microstructure layout. 
    The evolution of the objectives is shown in (d) in addition to the objective of the final design from GCMMA evaluated with 1000 random microstructure scenarios shown with a yellow square. }
    \label{fig:exI_caseI_Eratio10}
\end{figure}

We use GCMMA and Adam with a step size $\eta=0.05$. 
Figure \ref{fig:exI_caseI_Eratio10} shows the initial and optimized designs and objectives obtained from Adam and GCMMA. The color shading in the structure corresponds to the ratio $C_{1111}/C_{2222}$ of the first two diagonal elements of the constitutive matrix. This ratio represents the level of anisotropy of the two-phase composite. In this case, the final designs from Adam and GCMMA both achieve similar objectives, but the convergence of GCMMA is faster. To confirm the accuracy of estimating the expected values for the objective and constraint of the proposed approach, the final design from GCMMA is further evaluated for 1000 random configurations of the microstructure layout, {a reasonably large number. The objective value from these 1000 random configurations is shown with a (yellow) square in Figure \ref{fig:exI_caseI_Eratio10_obj}, which coincides with the (red) dashed curve for GCMMA at the end of the optimization verifying the convergence of the proposed stochastic gradient based approach using only four such configurations per iteration.} 

\begin{figure}[!htb]
    \centering
    \begin{subfigure}[b]{0.7\textwidth}
    \begin{tikzpicture}
        \node[inner sep=0pt] (structure) at (0,0){\includegraphics[scale=0.275]{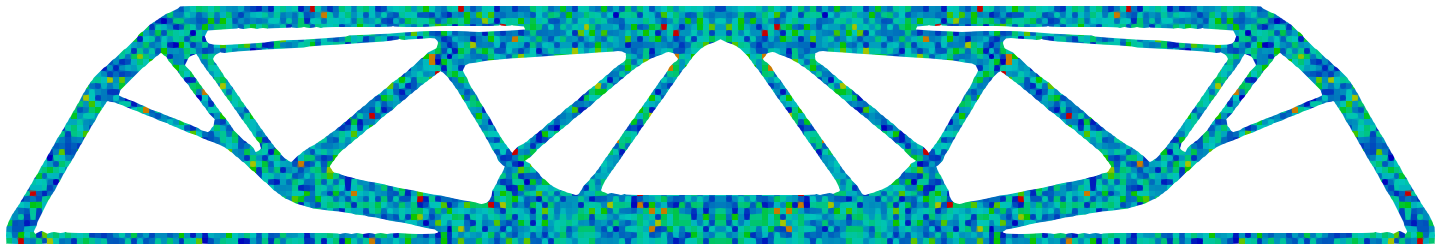}}; 
        \node[inner sep=0pt] (colorbar) at (5.5,0)
    {\includegraphics[scale=0.375]{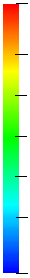}}; 
    \node at (5.75,1.8) {$C_{1111}/C_{2222}$}; 
    \node at (5.95,1.35) {$3.6$}; 
    \node at (5.95,0.45) {$2.5$}; 
    \node at (5.95,-0.35) {$1.5$}; 
    \node at (6,-1.3) {$0.29$}; 
    \draw[-latex] (-6.35,-0.9) -- (-5.6,-0.9); 
    \draw[-latex] (-6.35,-0.9) -- (-6.35,-0.15); 
    \node[] at (-5.45,-0.9) {$x$}; 
    \node[] at (-6.35,0) {$y$}; 
        \end{tikzpicture}
        \caption{Adam design}
    \end{subfigure}\\ 
    \begin{subfigure}[b]{0.7\textwidth}
    \begin{tikzpicture}
        \node[inner sep=0pt] (structure) at (0,0){\includegraphics[scale=0.275]{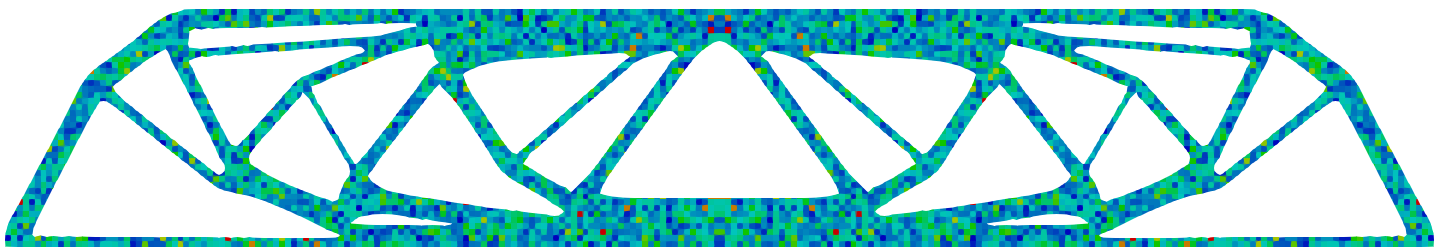}};
        \node[inner sep=0pt] (colorbar) at (5.5,0)
    {\includegraphics[scale=0.375]{Figures/colorbar_T4pi_E100.png}}; 
    \node at (5.75,1.8) {$C_{1111}/C_{2222}$}; 
    \node at (5.95,1.35) {$3.6$}; 
    \node at (5.95,0.45) {$2.5$}; 
    \node at (5.95,-0.35) {$1.5$}; 
    \node at (6,-1.3) {$0.29$}; 
    \draw[-latex] (-6.35,-0.9) -- (-5.6,-0.9); 
    \draw[-latex] (-6.35,-0.9) -- (-6.35,-0.15); 
    \node[] at (-5.45,-0.9) {$x$}; 
    \node[] at (-6.35,0) {$y$}; 
        \end{tikzpicture}
        \caption{GCMMA design}
    \end{subfigure}\\
    \begin{subfigure}[b]{0.55\textwidth} 
    \centering 
    \includegraphics[scale=0.275]{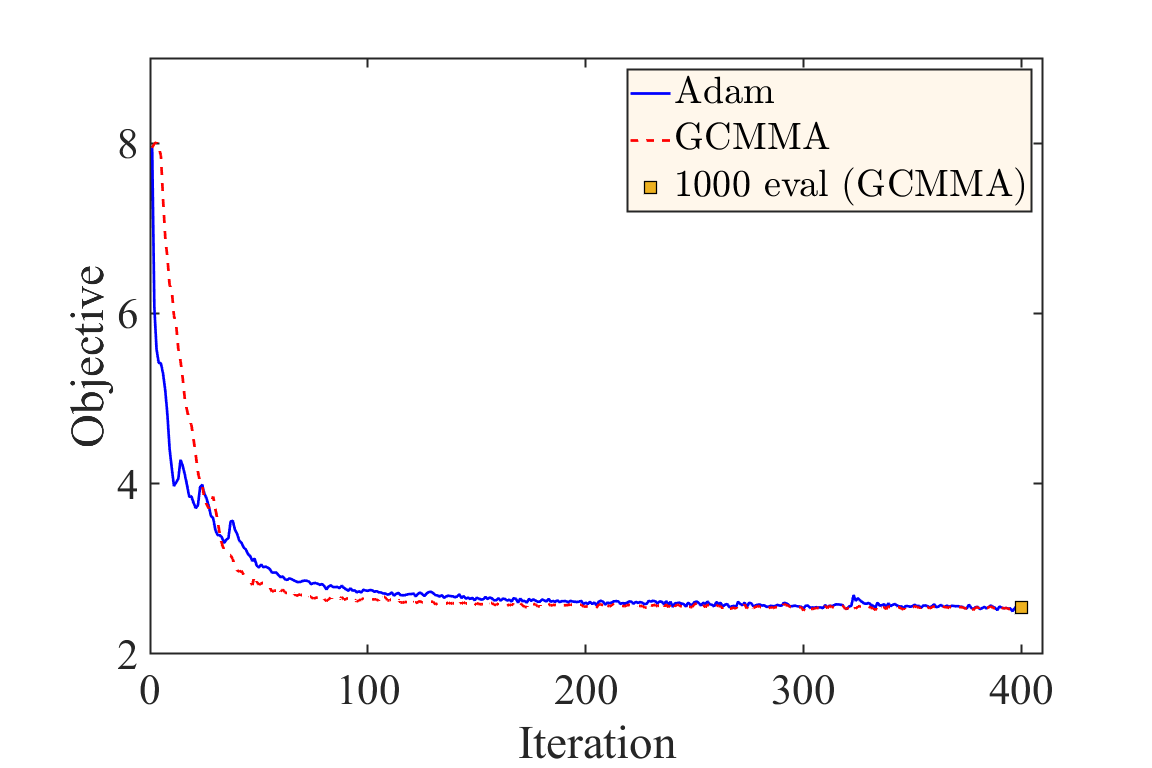}
    \caption{Objective vs. iteration} \label{fig:exI_caseI_Eratio100_obj}
    \end{subfigure}
    \caption{
    Comparison of the designs obtained from Adam (a) and GCMMA (b) for Case I of Example Ib using $E_1/E_2=100$. The color shading of the elements in (a) and (b) corresponds to the ratio $C_{1111}/C_{2222}$  for one random microstructure layout. 
    The evolution of the objectives is shown in (c) in addition to the objective of the final design from GCMMA evaluated with 1000 random microstructure scenarios shown with a yellow square. 
    }
    \label{fig:exI_caseI_Eratio100}
\end{figure}

\textbf{Case Ib: }
 To increase the level of material anisotropy we set the ratio $E_1/E_2=100$, but keep $T$ and $K$ the same as in Case Ia. 
 Figure \ref{fig:exI_caseI_Eratio100} shows the designs and objectives obtained from Adam and GCMMA with a step size $\eta = 0.05$. 
 The color shading shows that in this case the variation in $C_{1111}/C_{2222}$ is much larger compared to Case Ia with $E_1/E_2=10$. As a result, the final designs obtained from Adam and GCMMA both have more bars. The initial convergence of Adam is faster during the initial stages of the optimization. However, the final designs from Adam and GCMMA both have similar objective values at the end of the optimization. We further evaluate the final design from GCMMA by using 1000 random configurations of the microstructure layout. The (yellow) square in Figure \ref{fig:exI_caseI_Eratio100_obj} shows that it coincides with the objective values at the end of the optimization.

\begin{figure}[!htb]
    \centering
    \begin{subfigure}[b]{0.7\textwidth}
    \begin{tikzpicture}
        \node[inner sep=0pt] (structure) at (0,0){\includegraphics[scale=0.275]{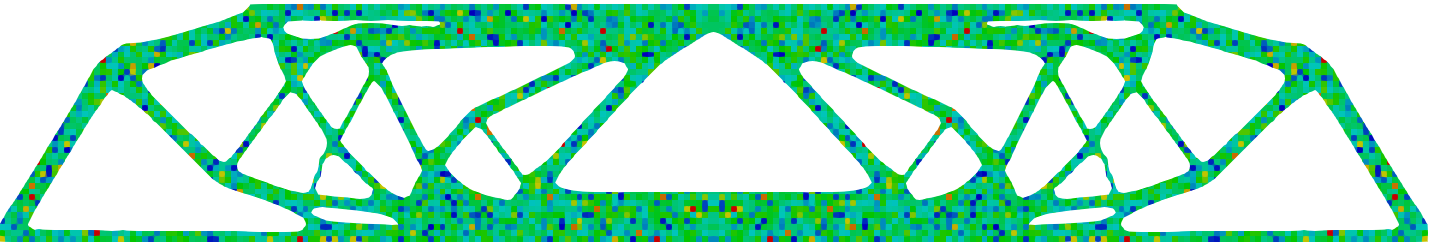}}; 
        \node[inner sep=0pt] (colorbar) at (5.5,0)
    {\includegraphics[scale=0.375]{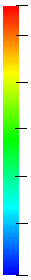}}; 
    \node at (5.75,1.8) {$C_{1111}/C_{2222}$}; 
    \node at (5.95,1.35) {$1.7$}; 
    \node at (5.95,0.6) {$1.4$}; 
    \node at (5.95,-0.35) {$1.0$}; 
    \node at (6,-1.3) {$0.57$}; 
    \draw[-latex] (-6.35,-0.9) -- (-5.6,-0.9); 
    \draw[-latex] (-6.35,-0.9) -- (-6.35,-0.15); 
    \node[] at (-5.45,-0.9) {$x$}; 
    \node[] at (-6.35,0) {$y$}; 
        \end{tikzpicture}
        \caption{Adam design}
    \end{subfigure}\\
    \begin{subfigure}[b]{0.7\textwidth}
    \begin{tikzpicture}
        \node[inner sep=0pt] (structure) at (0,0){\includegraphics[scale=0.275]{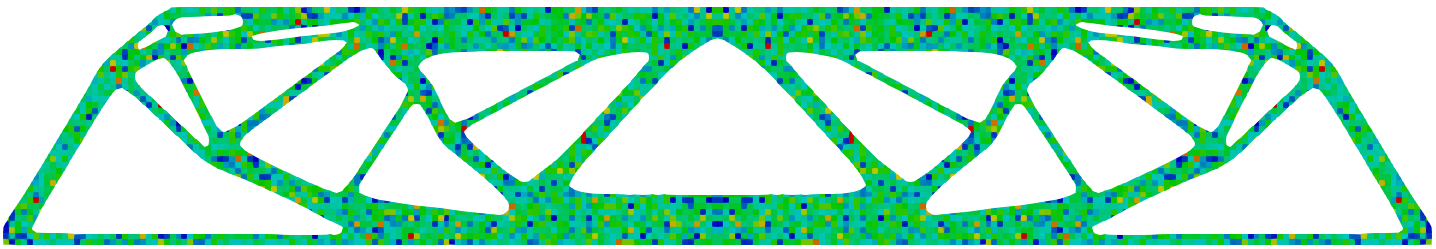}};
        \node[inner sep=0pt] (colorbar) at (5.5,0)
    {\includegraphics[scale=0.375]{Figures/colorbar_T2pi_E10.png}}; 
    \node at (5.75,1.8) {$C_{1111}/C_{2222}$}; 
    \node at (5.95,1.35) {$1.7$}; 
    \node at (5.95,0.6) {$1.4$}; 
    \node at (5.95,-0.35) {$1.0$}; 
    \node at (6,-1.3) {$0.57$}; 
    \draw[-latex] (-6.35,-0.9) -- (-5.6,-0.9); 
    \draw[-latex] (-6.35,-0.9) -- (-6.35,-0.15); 
    \node[] at (-5.45,-0.9) {$x$}; 
    \node[] at (-6.35,0) {$y$}; 
        \end{tikzpicture}
        \caption{GCMMA design}
    \end{subfigure}\\
    \begin{subfigure}[b]{0.55\textwidth} 
    \centering 
    \includegraphics[scale=0.275]{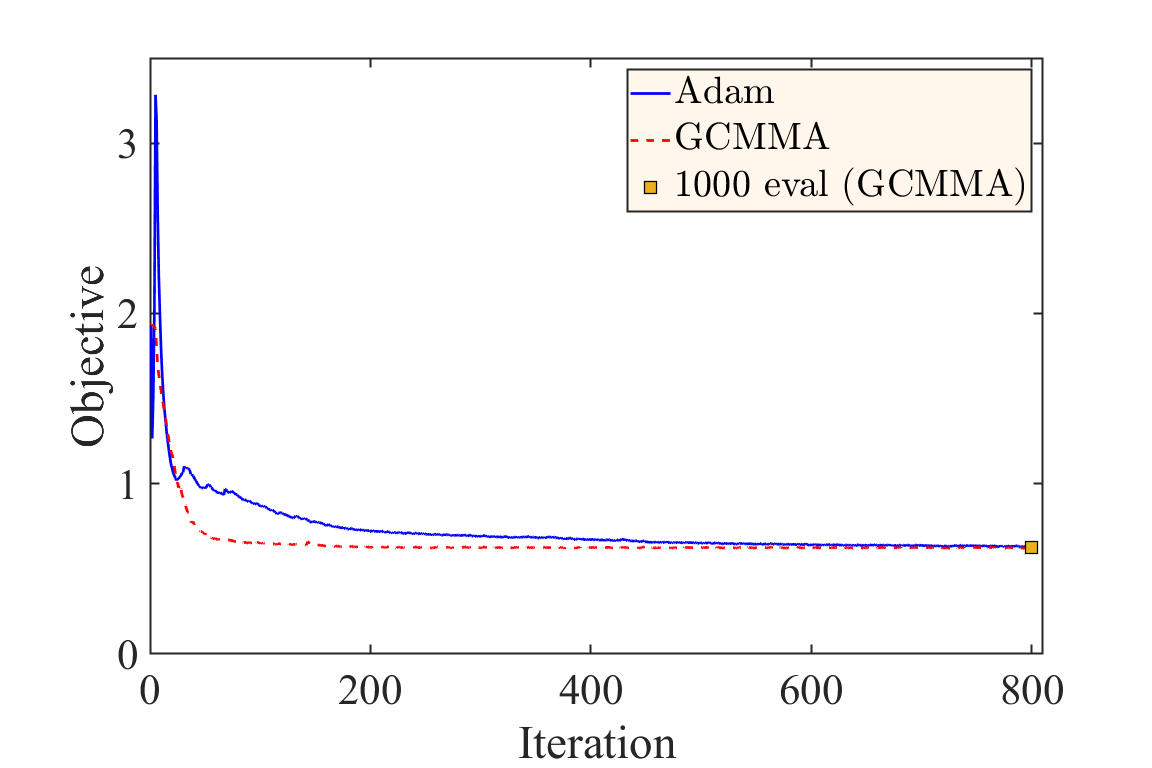}
    \caption{Objective vs. iteration}
    \label{fig:exI_caseII_Eratio10_obj}
    \end{subfigure}
    \caption{Comparison of the designs obtained from Adam (a) and GCMMA (b) for Case IIa of Example I using $E_1/E_2=10$. The color shading of the elements in (a) and (b) corresponds to the ratio $C_{1111}/C_{2222}$  for one random microstructure layout. 
    The evolution of the objectives is shown in (c) in addition to the objective for the final design from GCMMA evaluated with 1000 random microstructure scenarios shown with a yellow square. }
    \label{fig:exI_caseII_Eratio10}
\end{figure}

\textbf{Case IIa: } In this case, we generate 200 random microstructures using $T=2\pi$ and $K=50$ in \eqref{eq:random_field} to give $N=50$. The ratio of elastic moduli of the two phases are set as $E_1/E_2=10$. Figure \ref{fig:exI_caseII_micro} depicts four such microstructures, {which shows that in this case the dispersion of the two phases is coarser than Case Ia.} 
Adam and GCMMA with a step size $\eta=0.025$ 
are used for the optimization. 
Figure \ref{fig:exI_caseII_Eratio10} shows designs and objectives obtained from Adam and GCMMA for this case. As in Case Ia, final designs from Adam and GCMMA both achieve similar objectives, and again the convergence of GCMMA is faster. The final design from GCMMA, when evaluated for 1000 random configurations of the microstructure layout, produces an objective, shown with a (yellow) square in Figure \ref{fig:exI_caseII_Eratio10_obj}, that coincides with the (red) dashed curve for GCMMA. 

\begin{figure}[!htb]
    \centering
    \begin{subfigure}[b]{0.7\textwidth}
    \begin{tikzpicture}
        \node[inner sep=0pt] (structure) at (0,0){\includegraphics[scale=0.275]{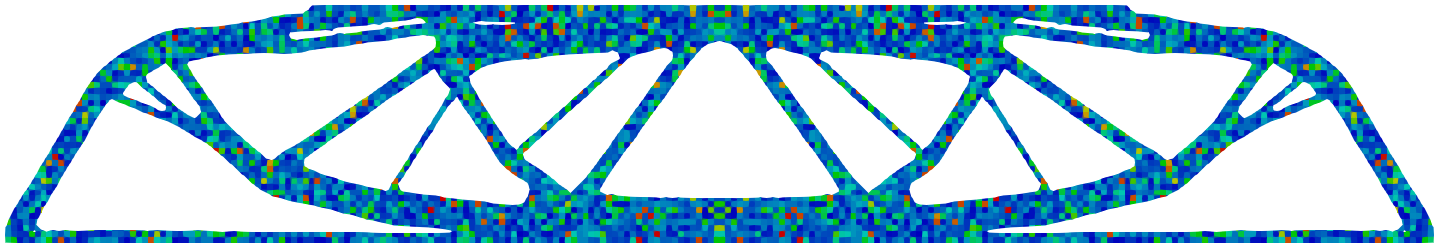}}; 
        \node[inner sep=0pt] (colorbar) at (5.5,0)
    {\includegraphics[scale=0.375]{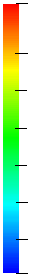}}; 
    \node at (5.75,1.8) {$C_{1111}/C_{2222}$}; 
    \node at (5.95,1.35) {$7.3$}; 
    \node at (5.95,0.9) {$6$}; 
    \node at (5.95,0.15) {$4$}; 
    \node at (5.95,-0.6) {$2$}; 
    \node at (6,-1.3) {$0.79$}; 
    \draw[-latex] (-6.35,-0.9) -- (-5.6,-0.9); 
    \draw[-latex] (-6.35,-0.9) -- (-6.35,-0.15); 
    \node[] at (-5.45,-0.9) {$x$}; 
    \node[] at (-6.35,0) {$y$}; 
        \end{tikzpicture}
        \caption{Adam design}
    \end{subfigure}\\ 
    \begin{subfigure}[b]{0.7\textwidth}
    \begin{tikzpicture}
        \node[inner sep=0pt] (structure) at (0,0){\includegraphics[scale=0.275]{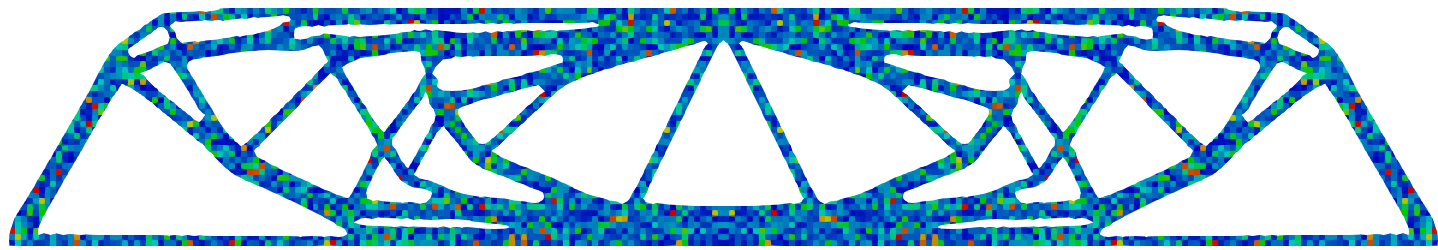}};
        \node[inner sep=0pt] (colorbar) at (5.5,0)
    {\includegraphics[scale=0.375]{Figures/colorbar_T2pi_E100.png}}; 
    \node at (5.75,1.8) {$C_{1111}/C_{2222}$}; 
    \node at (5.95,1.35) {$7.3$}; 
    \node at (5.95,0.9) {$6$}; 
    \node at (5.95,0.15) {$4$}; 
    \node at (5.95,-0.6) {$2$}; 
    \node at (6,-1.3) {$0.79$}; 
    \draw[-latex] (-6.35,-0.9) -- (-5.6,-0.9); 
    \draw[-latex] (-6.35,-0.9) -- (-6.35,-0.15); 
    \node[] at (-5.45,-0.9) {$x$}; 
    \node[] at (-6.35,0) {$y$}; 
        \end{tikzpicture}
        \caption{GCMMA design}
    \end{subfigure}\\
    \begin{subfigure}[b]{0.55\textwidth} 
    \centering 
    \includegraphics[scale=0.275]{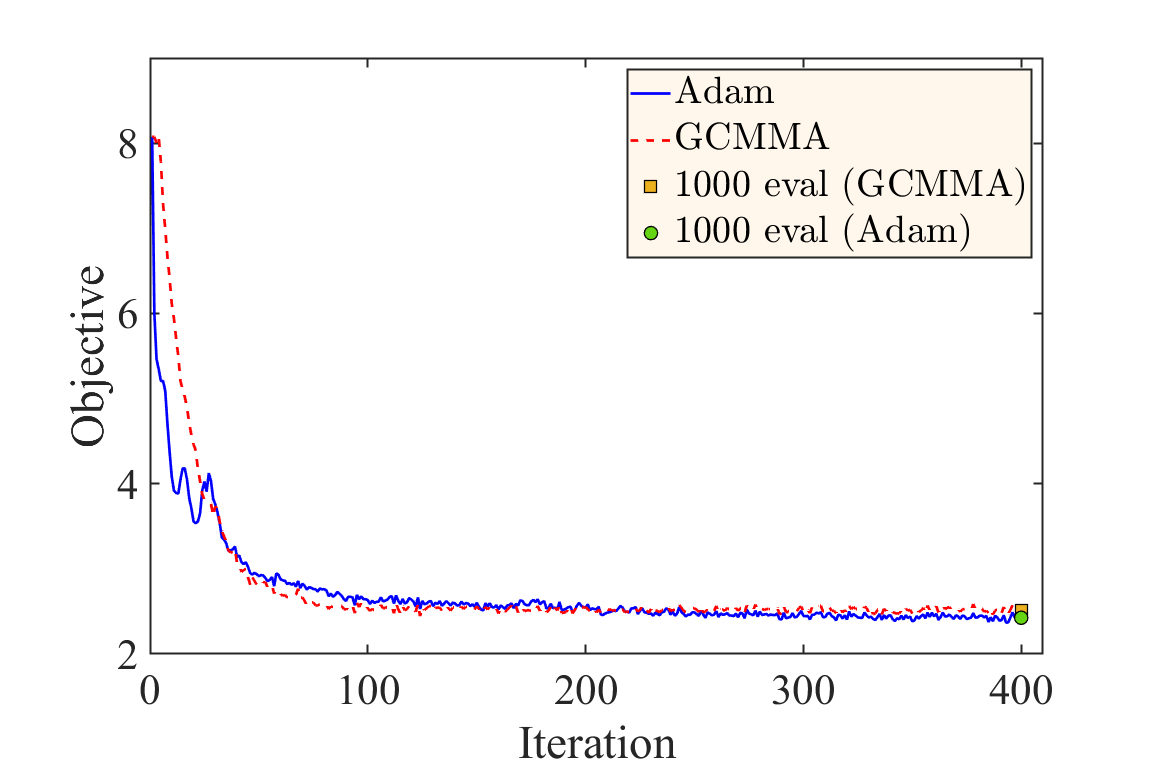}
    \caption{Objective vs. iteration}
    \label{fig:exI_caseII_Eratio100_obj}
    \end{subfigure}
    \caption{Comparison of the designs obtained from Adam (a) and GCMMA (b) for Case IIb of Example I using $E_1/E_2=100$. The color shading of the elements in (a) and (b) corresponds to the ratio $C_{1111}/C_{2222}$ for one random microstructure layout. 
    The objectives are shown in (c) in addition to the objectives for the final designs from GCMMA and Adam evaluated with 1000 random microstructure scenarios shown with a yellow square and a green circle, respectively. }
    \label{fig:exI_caseII_Eratio100}
\end{figure}

\textbf{Case IIb: }
Next, we use $E_1/E_2=100$ with $T$ and $K$ same as in Case IIa {to further increase the level of anisotropy.} 
Figure \ref{fig:exI_caseII_Eratio100} shows the designs and objectives obtained from Adam and GCMMA with step size $\eta=0.025$. 
Compared to Cases I and IIa, here the level of anisotropy is larger as can be seen from the $C_{1111}/C_{2222}$ ratios. To account for the increased anisotropy, the final designs obtained from Adam and GCMMA both have more bars. Figure \ref{fig:exI_caseII_Eratio100_obj}, however, shows that the final design from GCMMA has a slightly higher objective value compared to the Adam design. When evaluated with 1000 random configurations of the microstructure layout, the objective value of the GCMMA design is similar to the (red) dashed curve at the right end. Similarly, the (green) circle for the Adam design coincides with the (blue) solid line at the end of the optimization when evaluated with 1000 random samples. 
Note that the penalty approach used in Adam produces $C(\thetaa)\sim\mathcal{O}(10^{-6})$ to $\mathcal{O}(10^{-7})$ for all the cases in this example. However, GCMMA satisfies the mass constraint exactly. 

{This example shows that, in the presence of microstructural uncertainty, the use of stochastic gradients effectively reduce the computational cost of topology optimization at the macroscale. 
The designs obtained for different levels of anisotropy and different ratios of the elastic moduli of the two phases in the microstructure, however, have different features. The designs obtained from Adam are different from GCMMA, but their objective values remain similar. In the next two examples, we extend the proposed approach to design three-dimensional structures. 
}

\FloatBarrier
\subsection{Example II: Design of a Three-dimensional Beam} 

\subsubsection{Problem Formulation}
\begin{figure}[!htb]
    \centering
    \begin{tikzpicture}
    \node[inner sep=0pt] (beam) at (0,0)
    {\includegraphics[scale=0.5]{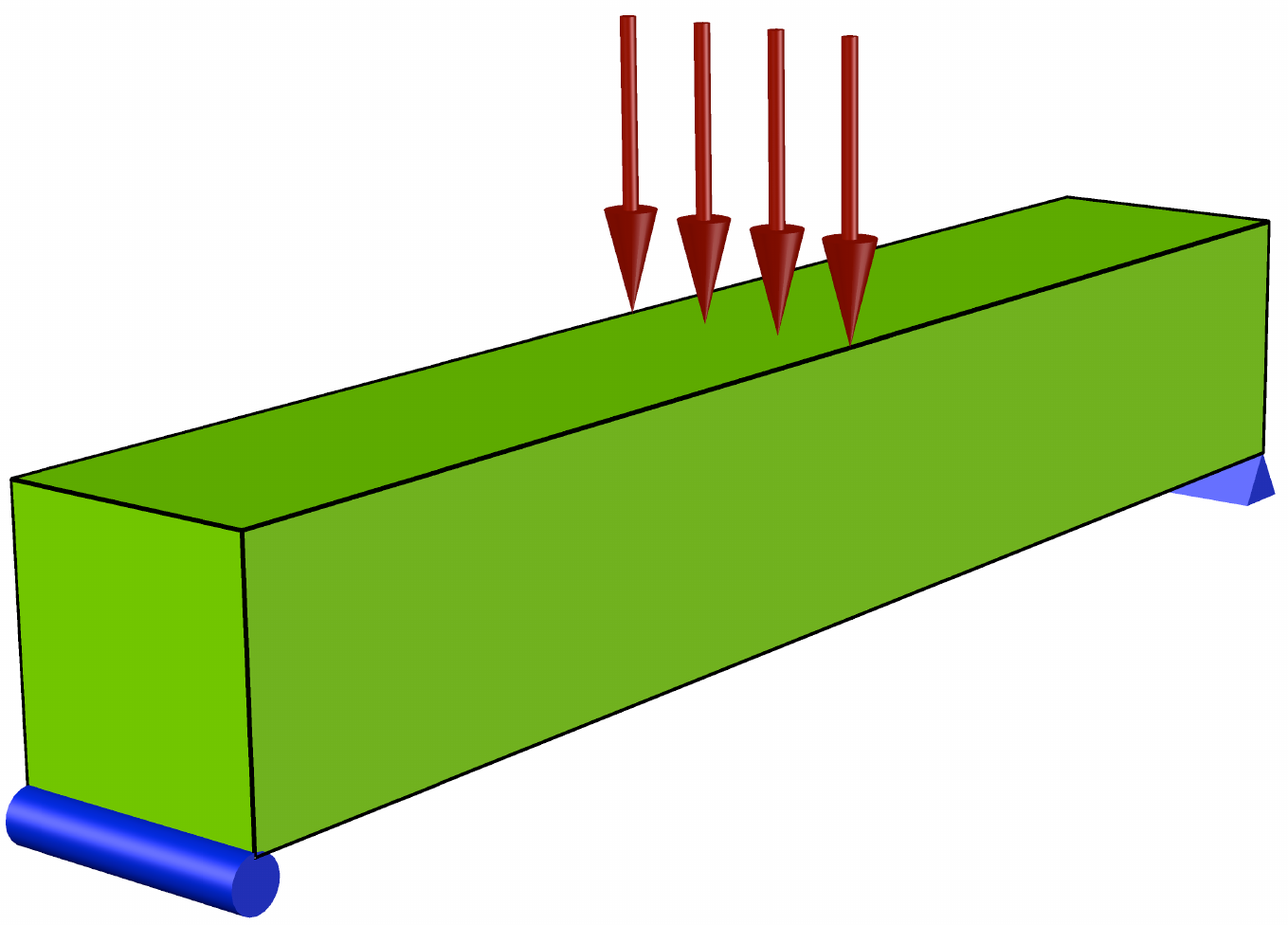}};
        \node[inner sep=0pt] (inset) at (4.5,0)
    {\includegraphics[scale=0.25]{Figures/micro-structure_boxed.png}}; 
    \draw[thick,latex-latex,white] (-2.1,-1.9) -- (3.3,0.25); 
\node[white,rotate=25] at (1,-0.35) {$L$}; 
\draw[thick,latex-latex,white] (-2.3,-2) -- (-2.35,-0.35); 
\node[white,rotate=-15] at (-2.7,-1.2) {$L/6$}; 
\draw[thick,latex-latex,white] (-2,-0.25) -- (-3.1,-0); 
\node[white,rotate=20] at (-2,0.1) {$L/6$}; 
    \draw[draw=gray!50,very thick] (2,0) circle (0.1cm); 
\draw[draw=gray!50,very thick] (4.5,0) circle (1cm); 
\draw[draw=gray!50,very thick] (2,-0.1) -- (4.1,-0.925); 
\draw[draw=gray!50,very thick] (2,0.1) -- (4.2,0.95); 
\node at (0.8,2.75) {$4P$}; 
    \draw[-latex] (-4.5,-2) -- (-4,-1.75); 
    \draw[-latex] (-4.5,-2) -- (-4.525,-1.5); 
    \draw[-latex] (-4.5,-2) -- (-4,-2.125); 
    \node[] at (-3.9,-1.65) {$x$}; 
    \node[] at (-4.525,-1.4) {$y$}; 
    \node[] at (-3.9,-2.25) {$z$}; 
    \end{tikzpicture}
    \caption{A three-dimensional beam composed of chopped-fiber composite is designed in Example II. }
    \label{fig:3d_mbb_schematic}
\end{figure}

\begin{figure}[!htb]
    \centering
    \begin{subfigure}[b]{0.75\textwidth}
    \centering 
    \begin{tikzpicture}
    \node[inner sep=0pt] (beam) at (0,0)
    {
    \includegraphics[scale=0.275]{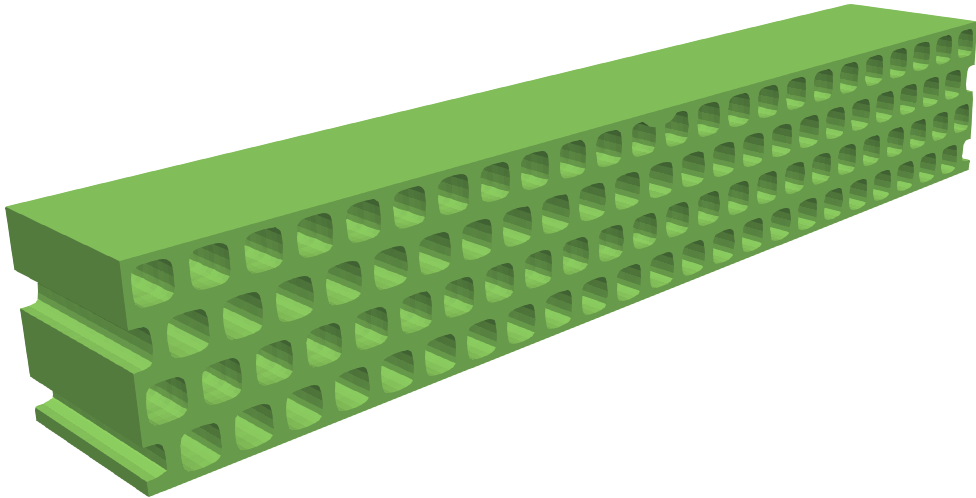}}; 
    \draw[-latex] (-4,-2) -- (-3.5,-1.75); 
    \draw[-latex] (-4,-2) -- (-4,-1.5); 
    \draw[-latex] (-4,-2) -- (-3.6,-2.25); 
    \node[] at (-3.4,-1.65) {$x$}; 
    \node[] at (-4.025,-1.4) {$y$}; 
    \node[] at (-3.4,-2.25) {$z$}; 
    \end{tikzpicture}
    \caption{Initial design } \label{fig:exII_init}
    \end{subfigure}\\
    \begin{subfigure}[b]{0.45\textwidth}
    \centering 
    \begin{tikzpicture}
    \node[inner sep=0pt] (beam) at (0,0)
    {
    \includegraphics[scale=0.275]{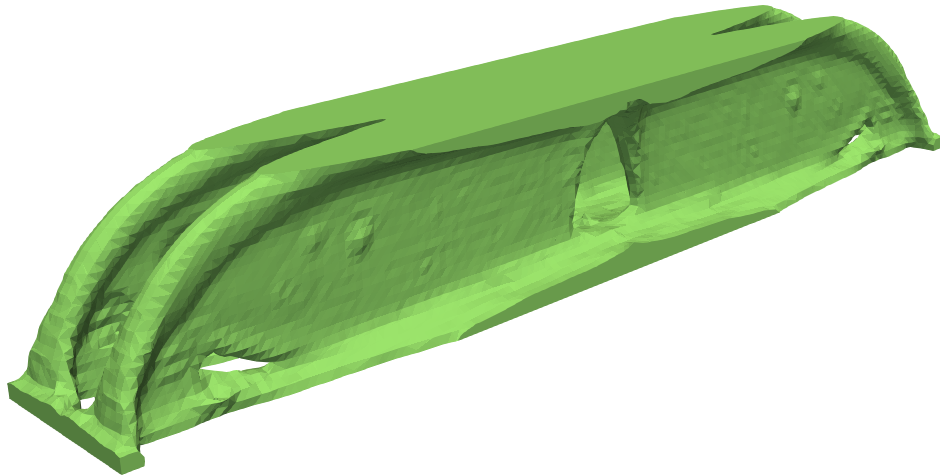}}; 
    \draw[-latex] (-4,-2) -- (-3.5,-1.75); 
    \draw[-latex] (-4,-2) -- (-4,-1.5); 
    \draw[-latex] (-4,-2) -- (-3.6,-2.25); 
    \node[] at (-3.4,-1.65) {$x$}; 
    \node[] at (-4.025,-1.4) {$y$}; 
    \node[] at (-3.4,-2.25) {$z$}; 
    \end{tikzpicture}
    \caption{Adam design } \label{fig:exII_adam}
    \end{subfigure}~~~~
    \begin{subfigure}[b]{0.45\textwidth}
    \centering 
    \begin{tikzpicture}
    \node[inner sep=0pt] (beam) at (0,0)
    {
    \includegraphics[scale=0.275]{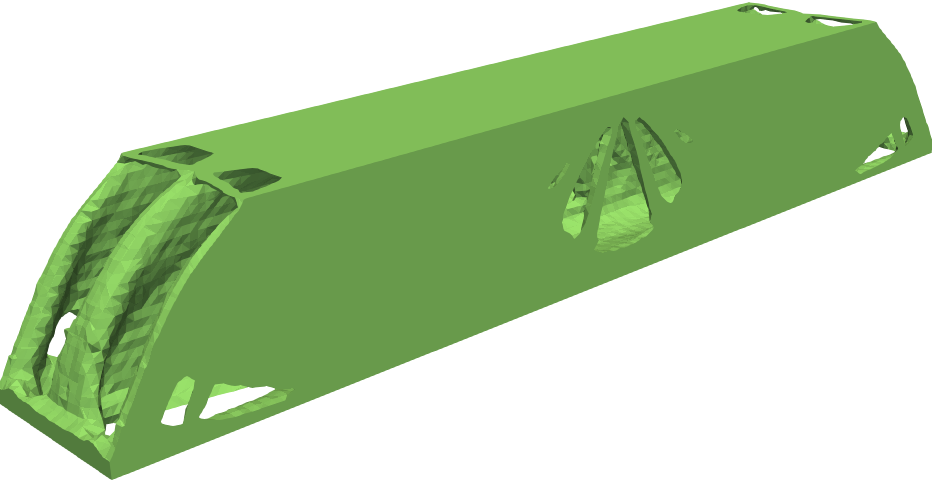}}; 
    \draw[-latex] (-4,-2) -- (-3.5,-1.75); 
    \draw[-latex] (-4,-2) -- (-4,-1.5); 
    \draw[-latex] (-4,-2) -- (-3.6,-2.25); 
    \node[] at (-3.4,-1.65) {$x$}; 
    \node[] at (-4.025,-1.4) {$y$}; 
    \node[] at (-3.4,-2.25) {$z$}; 
    \end{tikzpicture}
    \caption{GCMMA design } \label{fig:exII_gcmma}
    \end{subfigure}
    \caption{
    Comparison of the initial design (a) and designs obtained from Adam (b) and GCMMA (c) in Example II. }
    \label{fig:exII_struc} 
\end{figure}
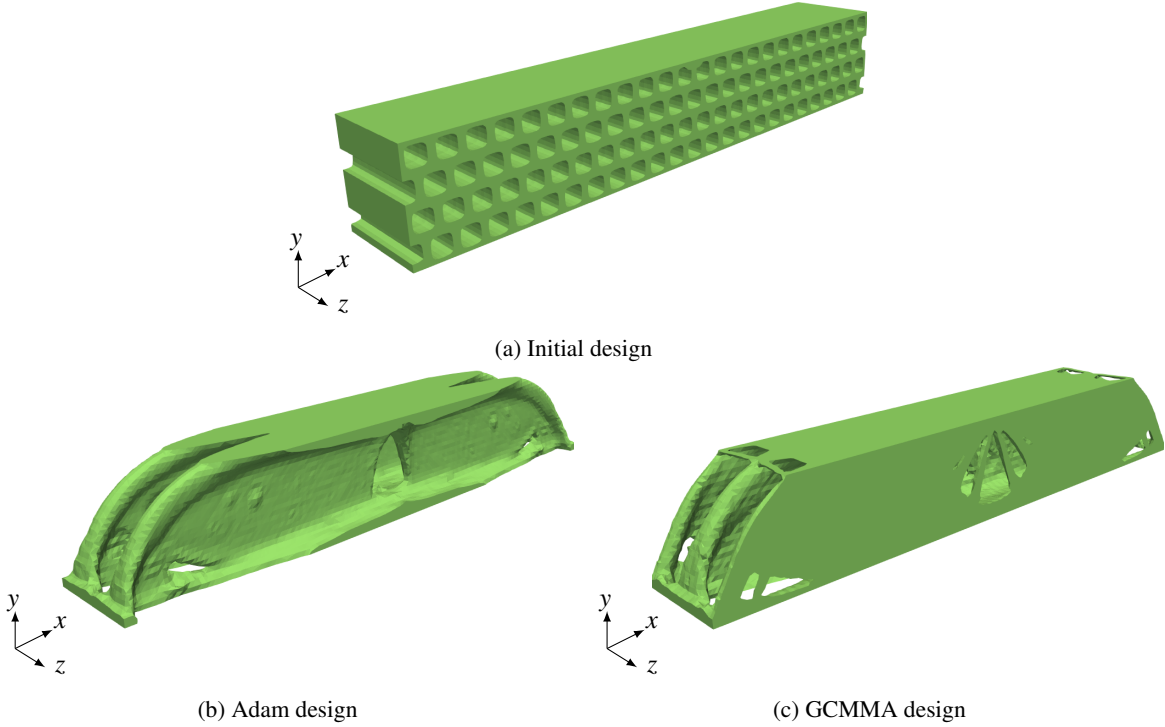

In our second example, we 
consider the design of a beam with a line load at the midspan (see Figure \ref{fig:3d_mbb_schematic}) subject to a 20\% mass constraint. Assuming symmetry, only one-fourth of the domain is used for optimization. 
In this example, we define the optimization problem as follows 
\begin{equation}
\begin{split}
    &\min_{\thetaa} ~~ R(\thetaa) = \Exp_{\xii}\left[{w_{\!_\Psi}}\frac{ {\Psi(\thetaa;\xii)} }{\Psi_0} + w_{\mathrm{m}}\frac{\mathcal{M}(\thetaa;\xii)}{\mathcal{M}_0} \right] + w_{\mathrm{per}} P_\mathrm{per}(\thetaa) + w_{\mathrm{reg}} P_\mathrm{reg}(\thetaa);\\
    &\text{subject to } C(\ppm) = g(\ppm) = 
    \frac{\int_{\Omega} \rho(\ppm)\mathrm{d}\xm}{\int_{\Omega}\mathrm{d}\xm}
    - \gamma_\mathrm{req}\leq 0,\\
    \end{split}
\end{equation}
where a term involving the total mass of the structure $\mathcal{M}(\thetaa;\xii)$ is also added to objective with $\mathcal{M}_0=\int_{\Omega} \mathrm{d}\xm$ and $w_\mathrm{m}=1$; and $\gamma_\mathrm{reqd}=0.20$ is used in the mass constraint. The two penalty terms are the same as in the previous example. For the weights $w_{\!_\Psi}$, $w_\mathrm{per}$, and $w_\mathrm{reg}$, we also use the same values as in the previous example. 
{We start the optimization from a structure with many square like holes as shown in Figure \ref{fig:exII_init} to facilitate hole seeding. 
In the initial design, we make 12 holes in the global $x$-direction for each of the 4 rows in the global $y$-direction, according to \eqref{eq:holes}, and with $r_\mathrm{hole}=0.09$. These holes are then protruded in the global $z$ direction along the width. 
We summarize the specifications used in this example in Table \ref{tab:ExII_specs} for convenient reproduction of the results presented herein. }

\begin{table}[!htb] 
\centering
\caption{Summary of specifications used to formulate and solve the optimization problem in Example II. The values are in consistent units. } \label{tab:ExII_specs} 
  \begin{tabular}{c|c|l|l}
    \hline
    \multicolumn{2}{c|}{\Tstrut Category} & Parameter & Value \Bstrut \\ \hline \multicolumn{2}{c|}{\Tstrut \multirow{6}{*}{Problem formulation}} & Weight for strain energy, $w_\Psi$ & 0.90 \\ 
    \multicolumn{2}{c|}{} & Weight for mass, $w_\mathrm{m}$ & 1.00 \\
    \multicolumn{2}{c|}{} & Weight for perimeter penalty, $w_\mathrm{per}$ & 0.025 \\
    \multicolumn{2}{c|}{} & Weight for regularization penalty, $w_\mathrm{reg}$ & 0.50 \\
    \multicolumn{2}{c|}{} & Line load, $P$ & 1.00 \\
    \multicolumn{2}{c|}{} & Mass constraint, $\gamma_\mathrm{reqd}$ & 0.20 \Bstrut \\ \hline 
    \multicolumn{2}{c|}{\Tstrut \multirow{4}{*}{Mesh (quarter domain)}} & Length, $L/2$ & 3.0 \\
    \multicolumn{2}{c|}{} & Height, $L/6$ & 1.0 \\ 
    \multicolumn{2}{c|}{} & Width, $L/12$ & 0.5 \\ 
    \multicolumn{2}{c|}{} & Discretization & $72\times24\times12$ \\ \hline 
    \multicolumn{2}{c|}{\Tstrut \multirow{5}{*}{Solution strategy}} & No. of possible microstructures per elem. & 200 \\ 
    \multicolumn{2}{c|}{} & Random config. per iter., $n_s$ & 4 \\
    \multicolumn{2}{c|}{} & No. of optimization variables, $n_{\thetaa}$ & $3.04\times10^4$ \\ 
    \multicolumn{2}{c|}{} & Step size, $\eta$ & $0.1$ \\ 
    \multicolumn{2}{c|}{} & Penalty to implement $C(\thetaa)$, $\kappa$ & $2000$ \\ 
    \hline
  \end{tabular} 
\end{table}

\begin{table}[!htb]
\caption{Uncertain parameters of the microstructure used in Example II. These parameters are assumed uniformly distributed between a lower and an upper limit. } 
\centering 
\begin{tabular}{c | c | c} 
\hline 
\Tstrut Parameter & Lower limit & Upper limit \\ [0.5ex] 
\hline  
\Tstrut Elastic modulus of fiber, $E_\mathrm{f}$ & 0.95 & 1.05  \\ 
Elastic modulus of matrix, $E_\mathrm{m}$ & 0.0095 & 0.0105  \\
Aspect ratio, $l/d$ & 10 & 100  \\
In-plane angle, $\theta_i$ & 0$^\circ$ & 180$^\circ$  \\  
Out-of-plane angle, $\theta_o$ & 0$^\circ$ & 180$^\circ$  \Bstrut \\  
\hline 
\end{tabular}
\label{tab:exII_fiber} 
\end{table} 

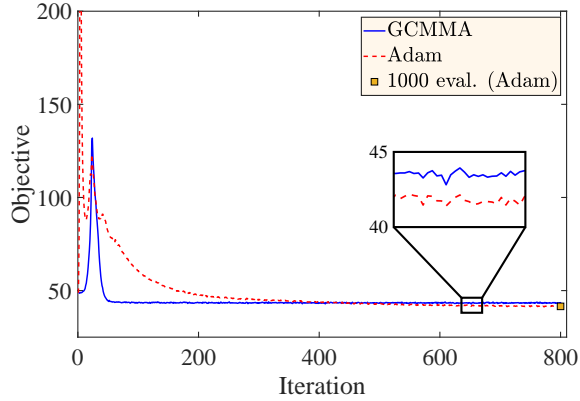
\begin{figure}[!htb]
    \centering
    \begin{tikzpicture}
    \node[inner sep=0pt] (plot) at (0,0)
    {\includegraphics[scale=0.275]{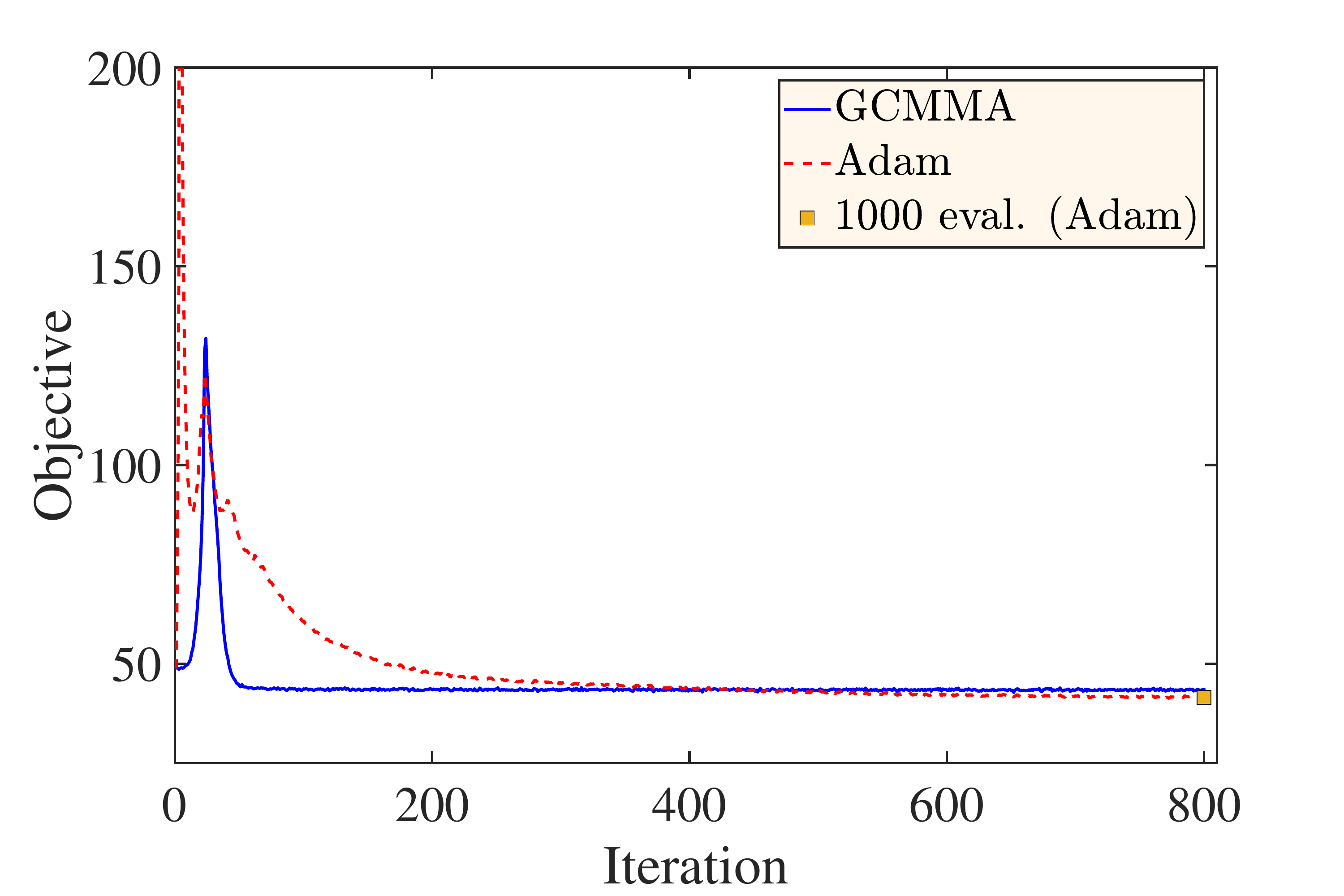}};
    \node[inner sep=0pt] (plot) at (1.98,0)
    {\includegraphics[scale=0.0775]{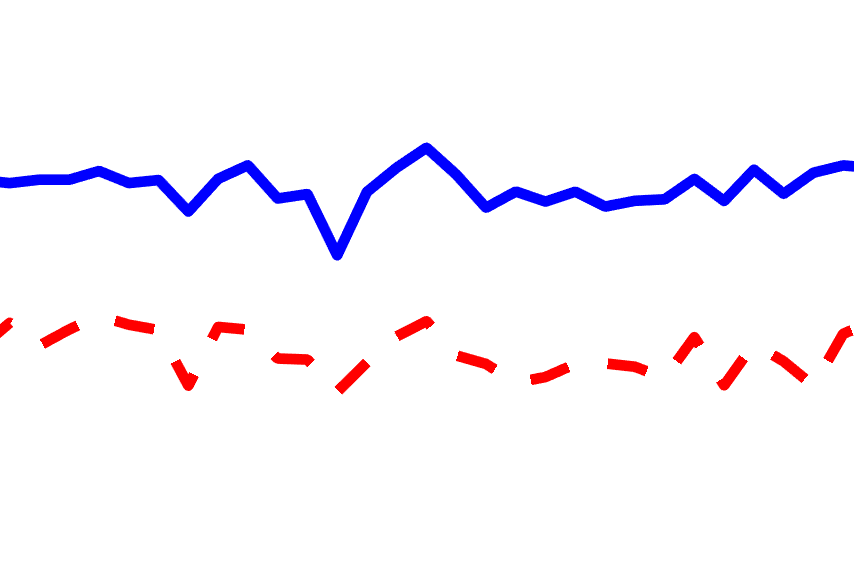}}; 
    \draw[draw=black,thick] (1.1,-0.5) rectangle ++(1.75,1); 
    \draw[draw=black,thick] (2,-1.64) rectangle ++(0.275,0.2); 
    \draw[thick] (2,-1.44) -- (1.1,-0.5); 
    \draw[thick] (2.275,-1.44) -- (2.85,-0.5); 
    \node[] at (0.9,0.5) {\scriptsize{45}}; 
    \node[] at (0.9,-0.5) {\scriptsize{40}}; 
    \end{tikzpicture}
    \caption{Comparison of the evolution of the objectives in Example II with a zoomed-in portion showing some oscillations at the end of the optimization due to the stochastic nature of gradients. The yellow square at the right end shows the objective from the final Adam design when evaluated using 1000 random microstructure layouts. }
    \label{fig:exII_obj_con} 
\end{figure} 

\begin{figure}[!htb]
    \centering
    \begin{subfigure}[b]{0.7\textwidth}
    \centering
    \begin{tikzpicture}
        \node[inner sep=0pt] (structure) at (0,0){
    \includegraphics[scale=0.275]{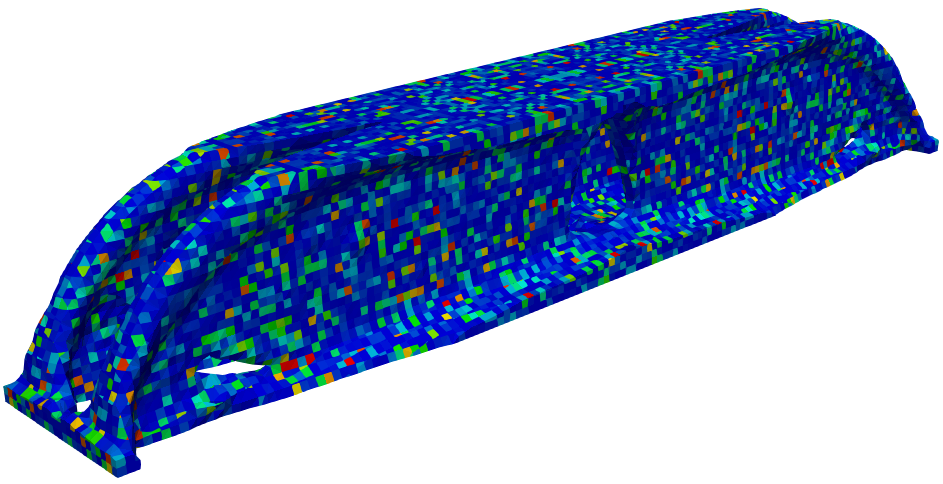}}; 
    \node[inner sep=0pt] (colorbar) at (4.5,0)
    {\includegraphics[scale=0.375]{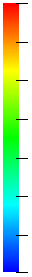}}; 
    \node at (4.75,1.8) {$C_{1111}/C_{2222}$}; 
    \node at (4.95,1.35) {$14.0$}; 
    \node at (4.95,0.625) {$10.0$}; 
    \node at (4.95,-0.2) {$6.0$}; 
    \node at (4.95,-0.95) {$2.0$}; 
    \node at (5,-1.3) {$0.07$}; 
    \draw[-latex] (-4,-2) -- (-3.5,-1.75); 
    \draw[-latex] (-4,-2) -- (-4,-1.5); 
    \draw[-latex] (-4,-2) -- (-3.6,-2.25); 
    \node[] at (-3.4,-1.65) {$x$}; 
    \node[] at (-4.025,-1.4) {$y$}; 
    \node[] at (-3.4,-2.25) {$z$}; 
        \end{tikzpicture} 
    \caption{Adam design } \label{fig:exII_adam_C11_C22}
    \end{subfigure}\\
    \begin{subfigure}[b]{0.7\textwidth}
    \centering
    \begin{tikzpicture}
        \node[inner sep=0pt] (structure) at (0,0){\includegraphics[scale=0.275]{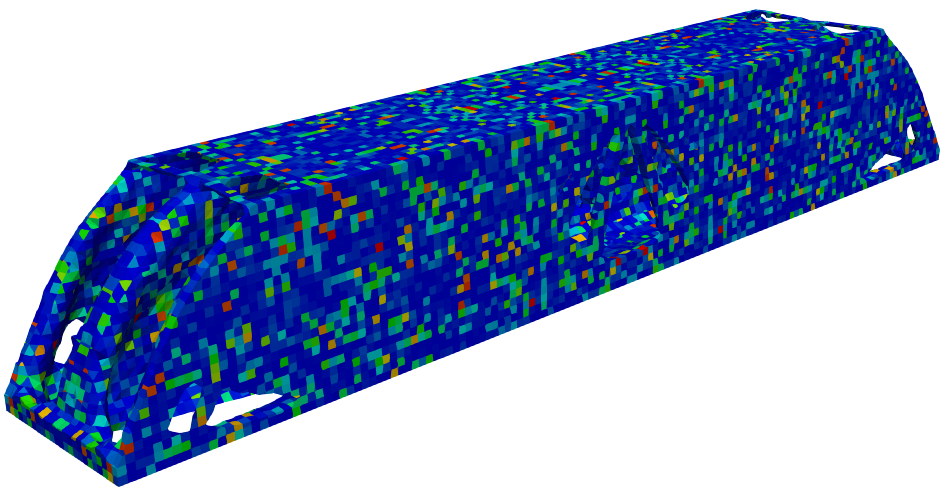}};
        \node[inner sep=0pt] (colorbar) at (4.5,0)
    {\includegraphics[scale=0.375]{Figures/colorbar_MBB3D.png}}; 
    \node at (4.75,1.8) {$C_{1111}/C_{2222}$}; 
    \node at (4.95,1.35) {$14.0$}; 
    \node at (4.95,0.625) {$10.0$}; 
    \node at (4.95,-0.2) {$6.0$}; 
    \node at (4.95,-0.95) {$2.0$}; 
    \node at (5,-1.3) {$0.07$}; 
    \draw[-latex] (-4,-2) -- (-3.5,-1.75); 
    \draw[-latex] (-4,-2) -- (-4,-1.5); 
    \draw[-latex] (-4,-2) -- (-3.6,-2.25); 
    \node[] at (-3.4,-1.65) {$x$}; 
    \node[] at (-4.025,-1.4) {$y$}; 
    \node[] at (-3.4,-2.25) {$z$}; 
        \end{tikzpicture}
    \caption{GCMMA design } \label{fig:exII_gcmma_C11_C22}
    \end{subfigure}
    \caption{The color shading shows the ratio $C_{1111}/C_{2222}$ of the first two diagonal elements of the constitutive tensor $\Chom$ for the designed structure using Adam and GCMMA in Example II for one random layout of the microstructure. }
    \label{fig:exII_C11_C22} 
\end{figure}

\subsubsection{Microstructure Scenarios} 
The microstructure in the beam is assumed to be chopped-fiber composite with uncertain properties listed in Table \ref{tab:exII_fiber}. 
We generate 200 realizations of the uncertain parameters of the microstructure from their respective uniform probability distributions as stated in Table \ref{tab:exII_fiber}. The homogenized constitutive properties are estimated using the Mori-Tanaka method described in Section \ref{sec:mori_tanaka}. The finite element model of one-quarter of the domain has a total $72\times24\times12=20,736$ elements. For one realization of the random microstructure layout, we assign each element of the finite element mesh to one randomly selected microstructure out of 200 possible ones, similar to the previous example. 

\subsubsection{Optimization Results} 
We use Adam and GCMMA with a step size $\eta=0.1$ and penalty parameter $\kappa=2000$ to implement the constraint and perform the optimization with four configurations of the microstructure layout per iteration for gradient calculations. 
Note that an assumption of symmetry holds here since every element in the finite element mesh is equally likely to have any of the 200 possible microstructures. 

The optimized structures obtained from Adam and GCMMA are shown in Figures \ref{fig:exII_adam} and \ref{fig:exII_gcmma}. Interestingly, Adam produces a design with two webs, but GCMMA produces a design with three webs. {However, both of these designs use the same amount of mass.} Figure \ref{fig:exII_obj_con} depicts the objective values for these two methods, where the Adam design with two webs can achieve a smaller objective when compared to a three-web design obtained using GCMMA. The objective of the final design from Adam is further verified with 1000 microstructural layout configurations and is shown in Figure \ref{fig:exII_obj_con} with a yellow square. Note that both GCMMA and Adam satisfy the constraint. The ratio $C_{1111}/C_{2222}$ of the first two diagonal elements in the homogenized constitutive tensor $\Chom$ for all the elements in the finite element model for one random layout of the microstructure for designs obtained from Adam and GCMMA shows the variability in material stiffness and anisotropy. As a result, the designs show a less smooth surface with kinks compared to a deterministic design. 
{Further constraints can be added to produce smoother surfaces more amenable to 3D printing \cite{schmitt2016formulation}, but that is beyond the scope of the current paper as our focus here is on the optimization process itself.} 

{The last two examples considered a common geometry used in the TO literature in two and three dimensions. The results show that the proposed approach is capable of producing average designs for two different types of random microstructures. In the next example, we use the stochastic gradient based approach to design a bracket for supporting a payload box to show the usefulness of the proposed approach for more practical design exercises. 
}

\subsection{Example III: Design of a Bracket}

\subsubsection{Problem Formulation}
\begin{figure}[!htb]
    \centering
    \begin{subfigure}[b]{0.45\textwidth}
    \centering
        \begin{tikzpicture}
    \node[inner sep=0pt] (beam) at (0,0)
    {\includegraphics[scale=0.5]{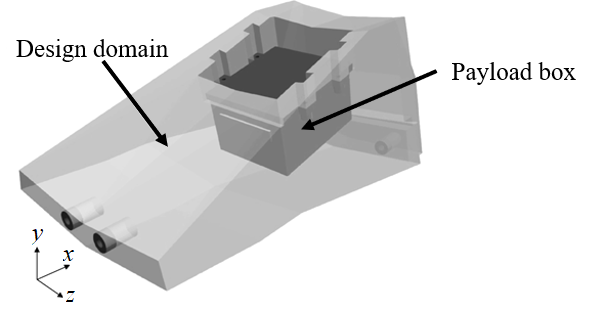}};
        \node[inner sep=0pt] (inset) at (2.75,-0.5)
    {\includegraphics[scale=0.075]{Figures/3d_mat_v2_1.png}};
    \draw[draw=orange!50,very thick] (0,-0.5) circle (0.1cm); 
\draw[draw=orange!50,very thick] (2.75,-0.5) circle (1cm); 
\draw[draw=orange!50,very thick] (0,-0.59) -- (2.35,-1.425); 
\draw[draw=orange!50,very thick] (0,-0.41) -- (2.45,0.45); 
    \end{tikzpicture}
    \caption{Isometric view} 
    \end{subfigure}\\
    \begin{subfigure}[b]{0.45\textwidth}
    \centering 
    \begin{tikzpicture}
    \node[inner sep=0pt] (beam) at (0,0)
    {
    \includegraphics[scale=0.5]{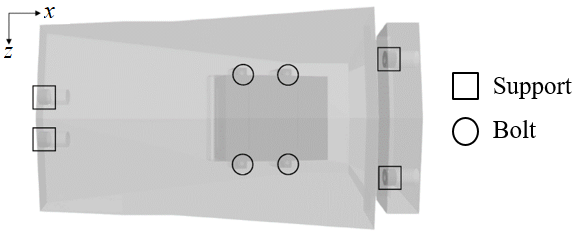}}; 
    \draw[thick,latex-latex] (-4,1.45) -- (-4,-1.5); 
    \draw[thick] (-3.75,-1.5) -- (-4.25,-1.5); 
    \draw[thick] (-3.75,1.45) -- (-4.25,1.45); 
    \node[rotate=90] at (-4.25,0) {10.16 cm}; 
    \draw[thick,latex-latex] (-3.45,-1.75) -- (1.75,-1.75); 
    \draw[thick] (-3.45,-1.5) -- (-3.45,-2); 
    \draw[thick] (1.75,-1.5) -- (1.75,-2); 
    \node[] at (-0.85,-2) {15.24 cm}; 
    \end{tikzpicture}
    \caption{Bottom view} 
    \end{subfigure}~~~~~~~~
    \begin{subfigure}[b]{0.45\textwidth}
    \centering
    \begin{tikzpicture}
    \node[inner sep=0pt] (beam) at (0,0)
    {
    \includegraphics[scale=0.35]{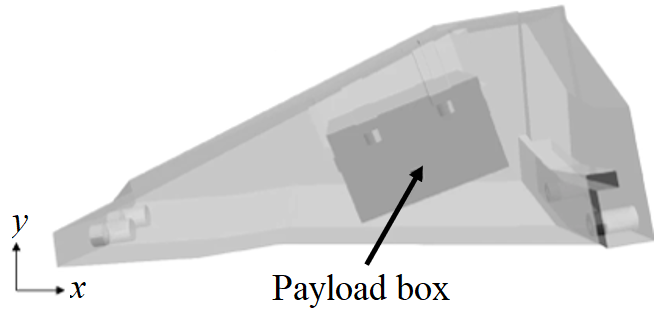}}; 
    \draw[thick,latex-latex] (3.4,1.45) -- (3.4,-1); 
    \draw[thick] (3.15,-1) -- (3.65,-1); 
    \draw[thick] (3.15,1.45) -- (3.65,1.45); 
    \node[rotate=90] at (3.65,0.225) {10.16 cm}; 
    \end{tikzpicture}
    \caption{Side view} 
    \end{subfigure}
    \caption{Schematic of the bracket design problem used in Example III with a random microstructure shown in the inset figure. The structure is supported at the left and right ends and designed to carry a payload box, which is connected to the structure using four bolts. }
    \label{fig:exIII_bracket_schem} 
\end{figure}

In the third example, we consider a structure to support a payload box given a set of supports at the left and right ends as shown in Figure \ref{fig:exIII_bracket_schem}. 
Note that the design domain is non-trivial, \textit{i.e.}, non box-shaped. This example is adapted from one of the challenge problems under DARPA TRADES program (\texttt{\href{http://solidmodeling.org/trades-cp/}{http://solidmodeling.org/trades-cp/}}) and Barrera et al. \cite{barrera2020hole,barrera2020ambiguous}. The payload is attached to the structure using four bolts. Uniform pressure of $1.2\times10^4$ N/cm$^2$ is applied on top of the payload, and an equivalent shock loading is applied to the entire assembly by subjecting it to a body force in the $y$-direction. 


\begin{figure}[!htb]
    \centering
    \begin{tikzpicture}
    \node[inner sep=0pt] (beam) at (0,0)
    {
    \includegraphics[scale=0.3]{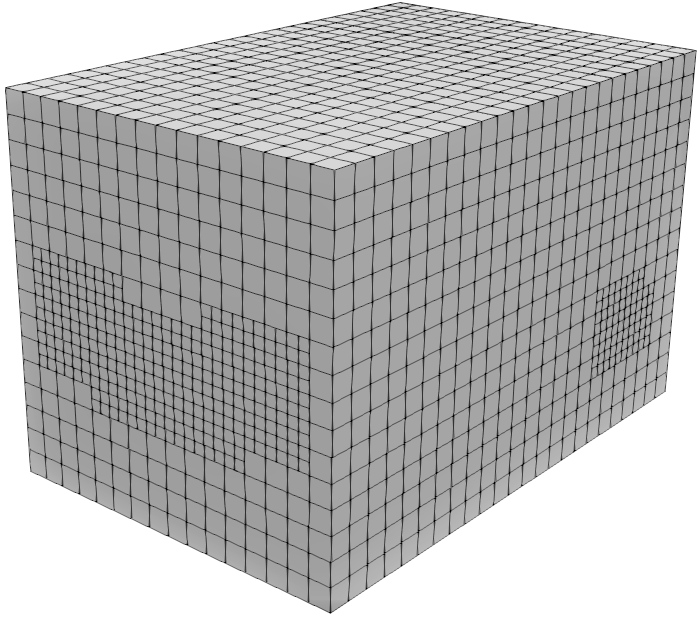}}; 
    \draw[ultra thick,black] (-1.5,2.5) -- (1.6,1.8); 
    \draw[ultra thick,black] (-1.5,2.5) -- (-1.5,2); 
    \draw[ultra thick,black] (1.6,1.8) -- (1.55,-1.9); 
    \draw[ultra thick,black] (1.55,-1.9) -- (0.95,-1.7); 
    \draw[ultra thick,black] (-3.2,0.2) -- (0.2,-1.05); 
    \draw[ultra thick,black] (-3.2,0.2) -- (-2.65,0.4); 
    \draw[ultra thick,black] (0.2,-1.05) -- (3.25,0.7); 
    \draw[ultra thick,black] (3.25,0.7) -- (2.65,0.8); 
    \draw[ultra thick,black] (-1.75,1.7) -- (1.75,2.65); 
    \draw[ultra thick,black] (-1.75,1.7) -- (-1.68,-2.3); 
    \draw[ultra thick,black] (-1.68,-2.3) -- (-1.05,-2); 
    \draw[ultra thick,black] (1.75,2.65) -- (1.75,2.2); 
    \node[inner sep=0pt,rotate=180] (beam) at (5.5,-2){
    \includegraphics[scale=0.175]{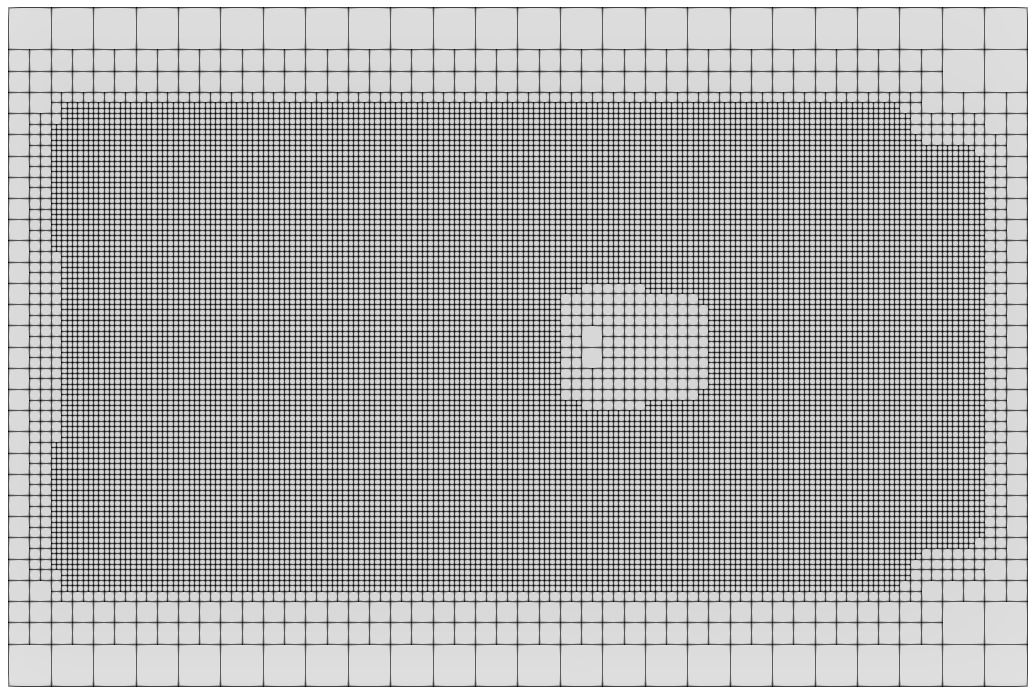}}; 
    \node[inner sep=0pt] (beam) at (-5,2.5){
    \includegraphics[scale=0.225]{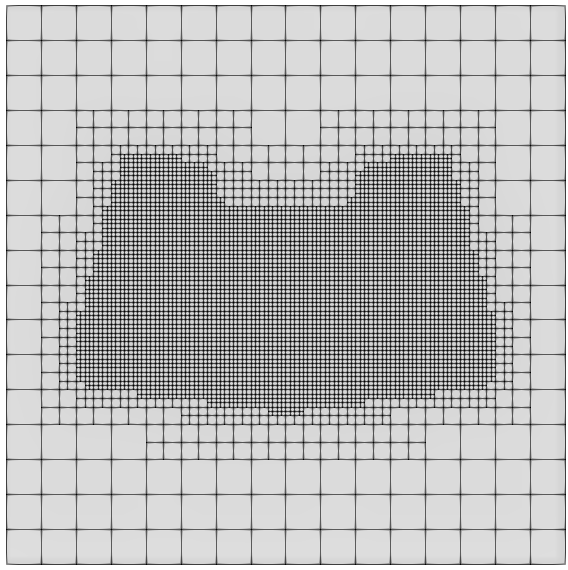}}; 
    \node[inner sep=0pt] (beam) at (6,2.5){
    \includegraphics[scale=0.15]{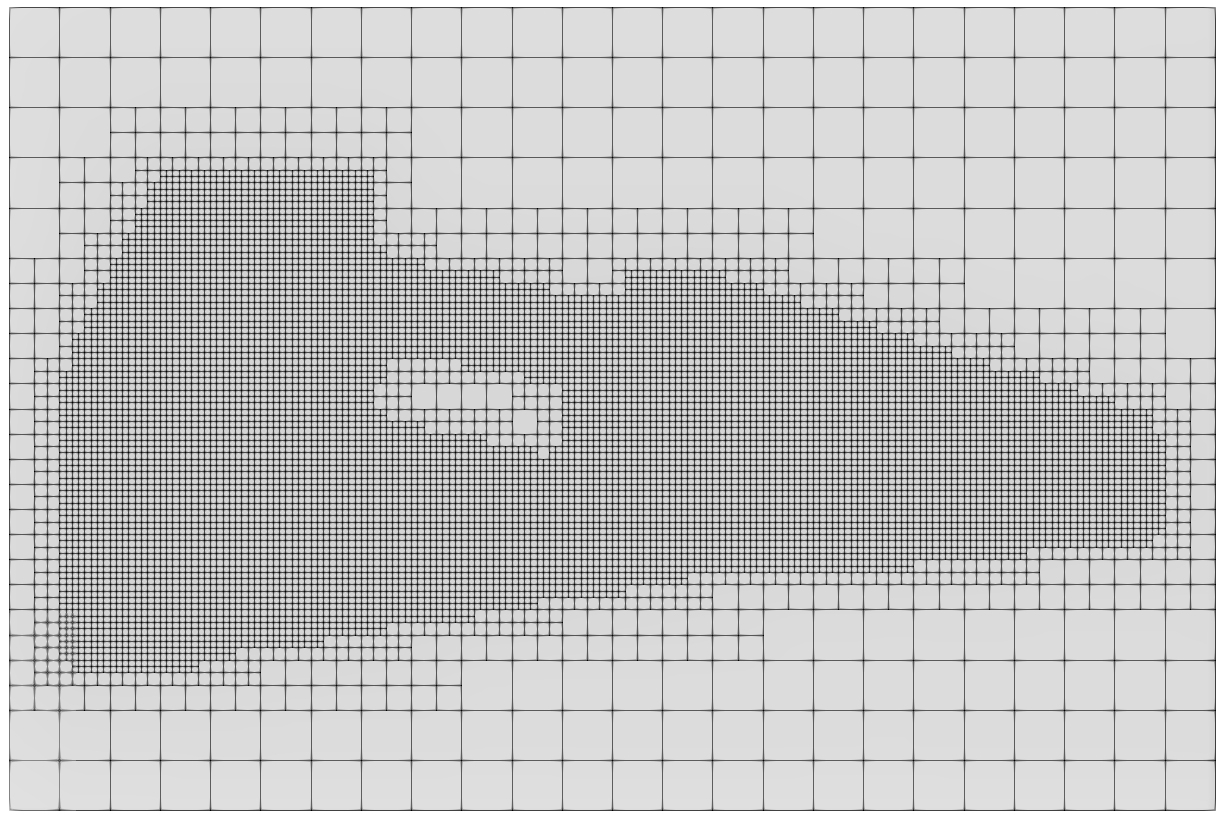}}; 
    \draw[-{Triangle[width=12pt,length=10pt]}, line width=5pt,color=black!60] (1.9,2.4) -- (3.6,2.5); 
    \draw[-{Triangle[width=12pt,length=10pt]}, line width=5pt,color=black!60] (3,0.4) -- (4,-0.35); 
    \draw[-{Triangle[width=12pt,length=10pt]}, line width=5pt,color=black!60] (-1.6,2.2) -- (-3.25,2.4); 
    \draw[-latex] (-3,-2) -- (-2.5,-1.7); 
    \draw[-latex] (-3,-2) -- (-3.025,-1.5); 
    \draw[-latex] (-3,-2) -- (-2.5,-2.25); 
    \node[] at (-2.4,-1.65) {$x$}; 
    \node[] at (-3.025,-1.4) {$y$}; 
    \node[] at (-2.4,-2.25) {$z$}; 
    \draw[-latex] (8.8,2) -- (9.5,2); 
    \draw[-latex] (8.8,2) -- (8.8,2.7); 
    \node[] at (9.7,2) {$x$}; 
    \node[] at (8.8,2.9) {$y$}; 
    \draw[-latex] (8.25,-2) -- (8.95,-2); 
    \draw[-latex] (8.25,-2) -- (8.25,-2.7); 
    \node[] at (9.15,-2) {$x$}; 
    \node[] at (8.25,-2.9) {$z$}; 
    \draw[-latex] (-2.9,3.2) -- (-2.2,3.2); 
    \draw[-latex] (-2.9,3.2) -- (-2.9,3.9); 
    \node[] at (-2.1,3.2) {$z$}; 
    \node[] at (-2.9,4.1) {$y$}; 
    \end{tikzpicture}
    \caption{Locally refined mesh with three cross-sectional views used in Example III. }
    \label{fig:exIII_bracket_mesh} 
\end{figure}
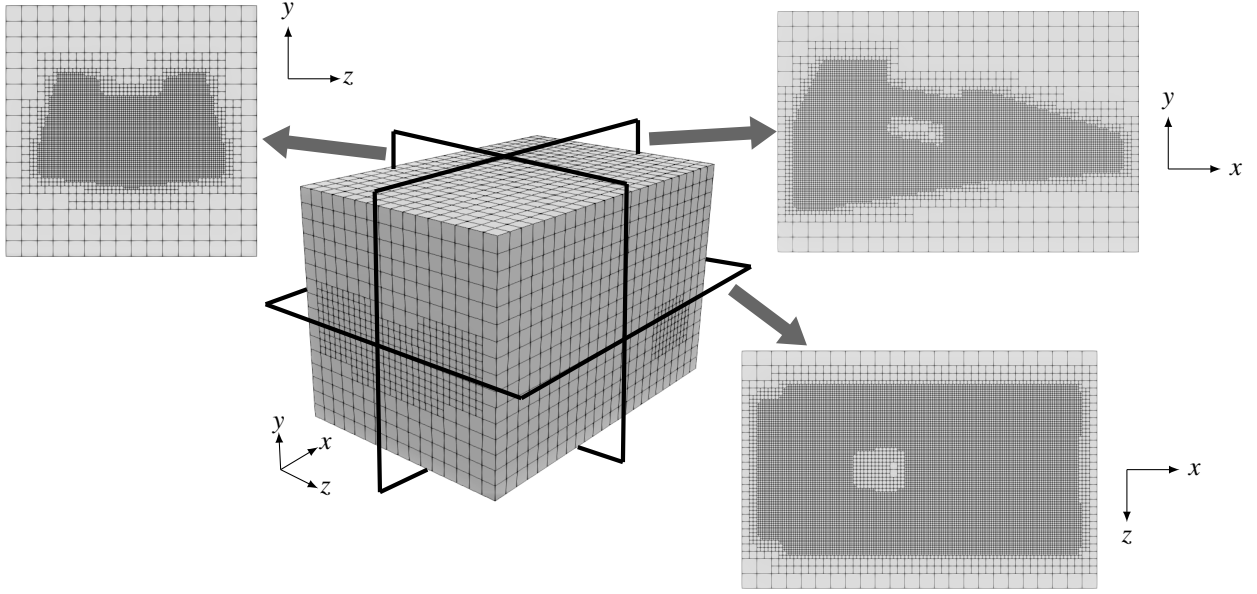 

\begin{figure}[!htb]
    \centering
    \begin{subfigure}[b]{0.475\textwidth}
    \centering 
    \begin{tikzpicture}
    \node[inner sep=0pt] (beam) at (0,0)
    {
    \includegraphics[scale=0.32]{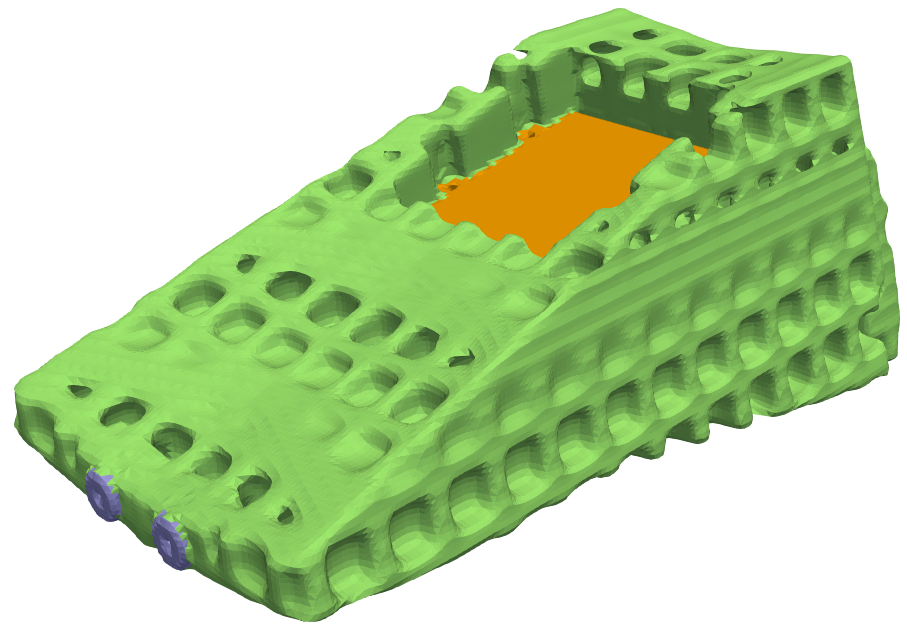}}; 
    \draw[-latex] (-4,-2) -- (-3.5,-1.75); 
    \draw[-latex] (-4,-2) -- (-4,-1.5); 
    \draw[-latex] (-4,-2) -- (-3.6,-2.25); 
    \node[] at (-3.4,-1.65) {$x$}; 
    \node[] at (-4.025,-1.4) {$y$}; 
    \node[] at (-3.4,-2.25) {$z$}; 
        \end{tikzpicture}
    \caption{Initial design } \label{fig:exIII_initial}
    \end{subfigure}~~~~
    \begin{subfigure}[b]{0.475\textwidth}
    \centering 
    \begin{tikzpicture}
    \node[inner sep=0pt] (beam) at (0.5,0)
    {
    \includegraphics[scale=0.275]{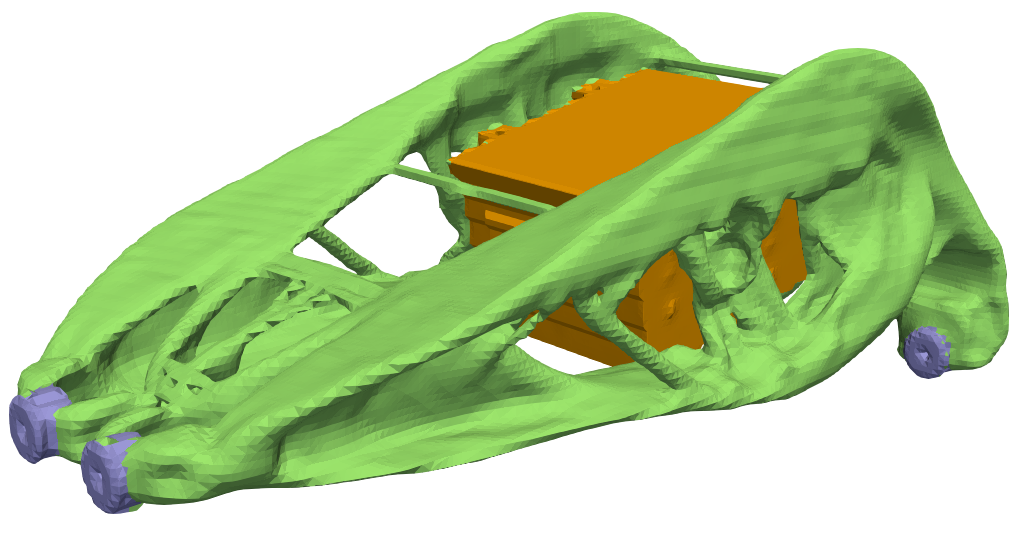}}; 
    \draw[-latex] (-4,-2) -- (-3.5,-1.75); 
    \draw[-latex] (-4,-2) -- (-4,-1.5); 
    \draw[-latex] (-4,-2) -- (-3.6,-2.25); 
    \node[] at (-3.4,-1.65) {$x$}; 
    \node[] at (-4.025,-1.4) {$y$}; 
    \node[] at (-3.4,-2.25) {$z$}; 
        \end{tikzpicture}
    \caption{Adam design} \label{fig:exIII_adam}
    \end{subfigure}\\
    \begin{subfigure}[b]{0.45\textwidth}
    \centering 
    \begin{tikzpicture}
    \node[inner sep=0pt] (beam) at (0.6,0)
    {
    \includegraphics[scale=0.275]{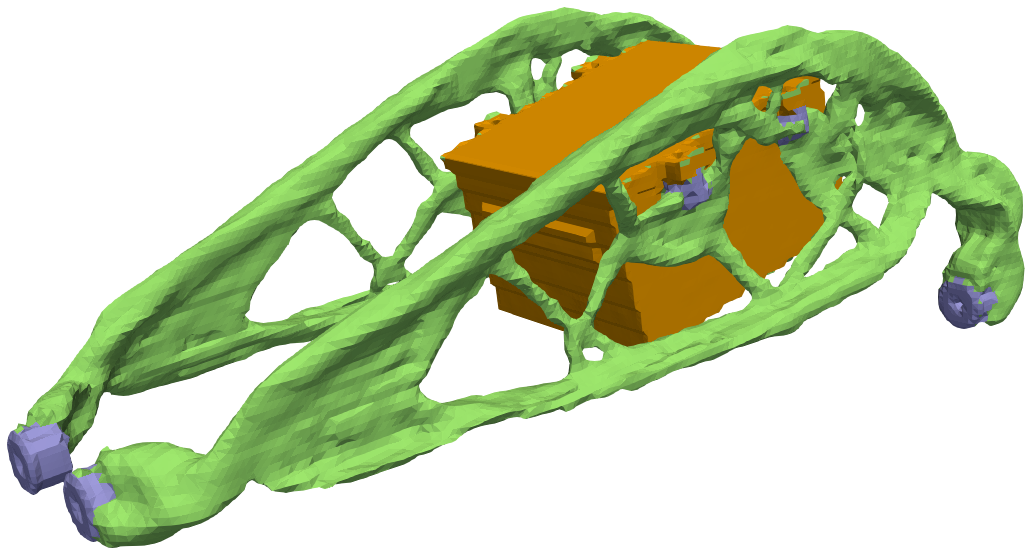}}; 
    \draw[-latex] (-4,-2) -- (-3.5,-1.75); 
    \draw[-latex] (-4,-2) -- (-4,-1.5); 
    \draw[-latex] (-4,-2) -- (-3.6,-2.25); 
    \node[] at (-3.4,-1.65) {$x$}; 
    \node[] at (-4.025,-1.4) {$y$}; 
    \node[] at (-3.4,-2.25) {$z$}; 
        \end{tikzpicture}
    \caption{Deterministic design} \label{fig:exIII_det}
    \end{subfigure}
    \caption{We start the optimization from the initial design shown in (a) and Adam produces the final design shown in (b). If we assume the Ti-6Al-4V alloy does not have any random impurities, i.e., no microstructural uncertainties, then a deterministic design process using GCMMA results in the design shown in (c), where members are slender. }
    \label{fig:exIII_design} 
\end{figure}
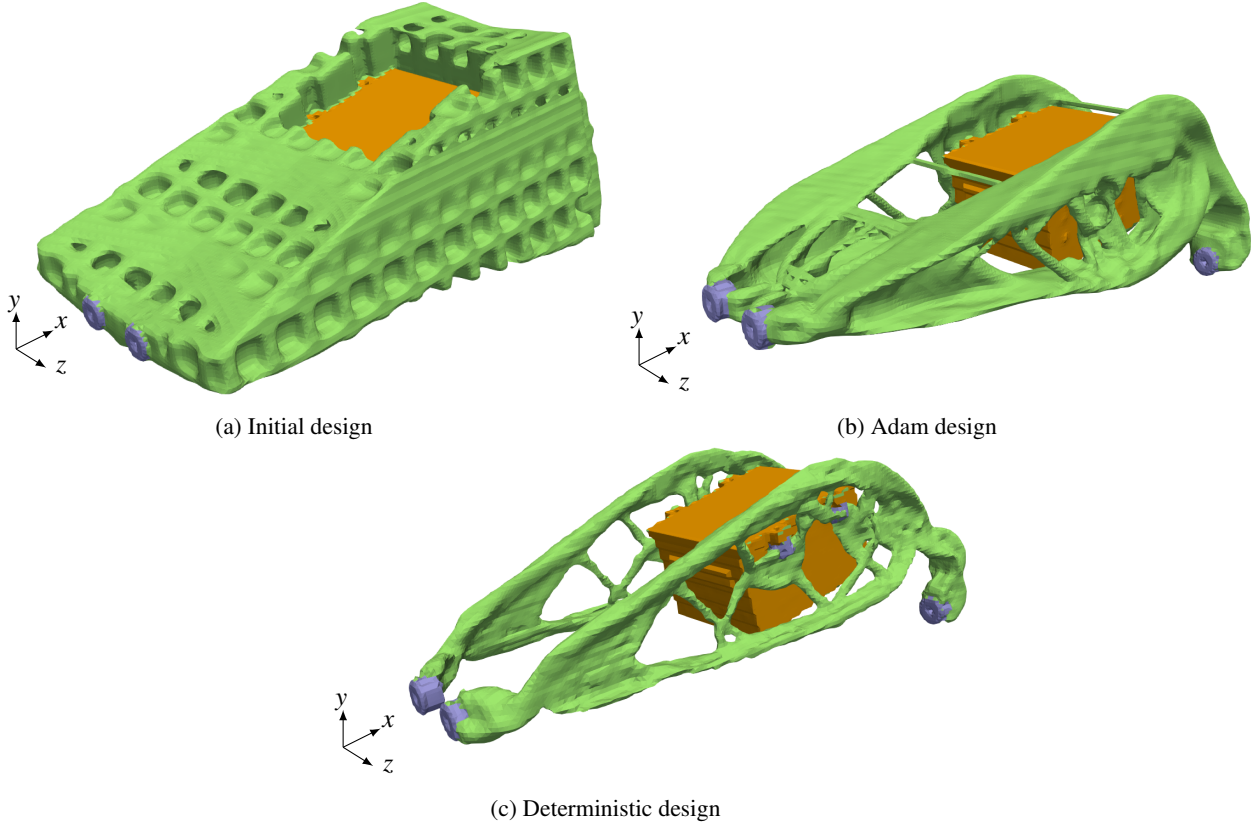

In this example, we minimize a combination of strain energy and mass of the structure subject to a mass constraint of 30\%. The optimization problem is defined as follows
\begin{equation}
\begin{split}
    &\min_{\thetaa} ~~ R(\thetaa) = \Exp_{\xii}\left[{w_{\!_\Psi}}\frac{ {\Psi(\thetaa;\xii)} }{\Psi_0} + w_{\mathrm{m}}\frac{\mathcal{M}(\thetaa;\xii)}{\mathcal{M}_0} \right] + w_{\phi}P_{\tilde{\phi}}(\thetaa) + w_{\mathrm{per}} P_\mathrm{per}(\thetaa) + w_{\mathrm{reg}} P_\mathrm{reg}(\thetaa);\\
    &\text{subject to } C_1(\ppm) = g_1(\ppm) = 
    \frac{\int_{\Omega} \rho(\ppm)\mathrm{d}\xm}{\int_{\Omega}\mathrm{d}\xm}
    - \gamma_\mathrm{req}\leq 0,\\
    &\qquad \qquad ~ C_2(\ppm) = g_2(\ppm) = P_{\tilde{\phi}}(\thetaa) = 0 
    \end{split}
\end{equation}
where the mass of the structure is $\mathcal{M}(\thetaa)=\int_{\Omega}\rho(\thetaa)\mathrm{d}\xm$; $\mathcal{M}_0=\int_{\Omega} \mathrm{d}\xm$; the strain energy of the structure is $\Psi(\thetaa;\xii) = \int_{\Omega} \overline{{\sigmaa}}(\xm;\ppm,\xii):\overline{\epsl}(\xm;\ppm,\xii) \mathrm{d}\xm$; $\Psi_0$ is the strain energy of the initial design; $P_{\tilde{\phi}}(\thetaa)=\left({\int_{\Omega}\lvert \phi(\xm;\thetaa)-\tilde\phi(\xm;\thetaa) \rvert\mathrm{d}\xm}\right)^2/{\int_{\Gamma_D}\mathrm{d}A}$ is a penalty term used to obtain a structure that connects to the payload box at the bolts only with $\tilde{\phi}$ defined as a level set field that is trimmed to the lower or upper limit at the bolts and supports; the other two penalty terms are same as before. We use $w_{\!_\Psi}=0.05$, $w_\mathrm{m}=50$, $w_\mathrm{\tilde{\phi}}=5\times10^4$, $w_\mathrm{per}=0.1$, and $w_\mathrm{reg}=1.0$ to define the optimization problem. The geometry of the structure is defined using a level set description (see Section \ref{sec:level_set}) {with locally refined mesh shown in Figure \ref{fig:exIII_bracket_mesh}.} The finite element model of the structure has a total of 833,306 elements. We start the optimization from the initial design shown in Figure \ref{fig:exIII_initial}. {We summarize the specifications used in this example in Table \ref{tab:ExIII_specs}.} 



\begin{table}[!htb] 
\centering
\caption{Summary of specifications used to formulate and solve the optimization problem in Example III. The values are in consistent units. } \label{tab:ExIII_specs} 
  \begin{tabular}{c|c|l|l}
    \hline
    \multicolumn{2}{c|}{\Tstrut Category} & Parameter & Value \Bstrut \\ \hline \multicolumn{2}{c|}{\Tstrut \multirow{7}{*}{Problem formulation}} & Weight for strain energy, $w_\Psi$ & 0.05 \\ 
    \multicolumn{2}{c|}{} & Weight for mass, $w_\mathrm{m}$ & 50.0 \\
    \multicolumn{2}{c|}{} & Weight for bolt connection penalty, $w_\mathrm{\tilde{\phi}}$ & $5\times10^4$ \\
    \multicolumn{2}{c|}{} & Weight for perimeter penalty, $w_\mathrm{per}$ & 0.10 \\
    \multicolumn{2}{c|}{} & Weight for regularization penalty, $w_\mathrm{reg}$ & 1.00 \\
    \multicolumn{2}{c|}{} & Uniform pressure on the payload & 120 MPa \Bstrut \\ 
    \multicolumn{2}{c|}{} & Mass constraint, $\gamma_\mathrm{reqd}$ & 0.30 \Bstrut \\ \hline 
    \multicolumn{2}{c|}{\Tstrut \multirow{7}{*}{Microstructure}} & Period, $T$ & $4\pi$ \\ 
    \multicolumn{2}{c|}{} & Max. wavenumber, $K$ & 25 \\
    \multicolumn{2}{c|}{} & Elastic modulus of Ti-6Al-4V alloy, $E_\mathrm{Ti-6Al-4V}$ & $1.138\times10^5$ MPa \\ 
    \multicolumn{2}{c|}{} & Elastic modulus of impurities, $E_\mathrm{imp}$ & $2.276\times10^3$ MPa \\ 
    \multicolumn{2}{c|}{} & Poisson's ratio of Ti-6Al-4V alloy, $\nu_\mathrm{Ti-6Al-4V}$ & 0.342\\ 
    \multicolumn{2}{c|}{} & Poisson's ratio of impurities, $\nu_\mathrm{imp}$ & 0.3\\ 
    \multicolumn{2}{c|}{} & Density of the impure alloy & 44.3 kg/m$^3$\\ \hline
    \multicolumn{2}{c|}{\Tstrut \multirow{5}{*}{Solution strategy}} & No. of possible microstructures per elem. & 200 \\ 
    \multicolumn{2}{c|}{} & Random config. per iter., $n_s$ & 1 \\
    \multicolumn{2}{c|}{} & No. of optimization variables, $n_{\thetaa}$ & $4.03\times10^5$ \\ 
    \multicolumn{2}{c|}{} & Step size, $\eta$ & 0.05 \\ 
    \multicolumn{2}{c|}{} & Penalties to implement $C_1(\thetaa)$ and $C_2(\thetaa)$, $\kappaa$ & $[10^4,10^4]^T$ \\ 
    \hline
  \end{tabular} 
\end{table}

\begin{figure}[!htb]
    \centering
    \begin{tikzpicture}
    \node[inner sep=0pt] (plot) at (0,0)
    {\includegraphics[scale=0.275]{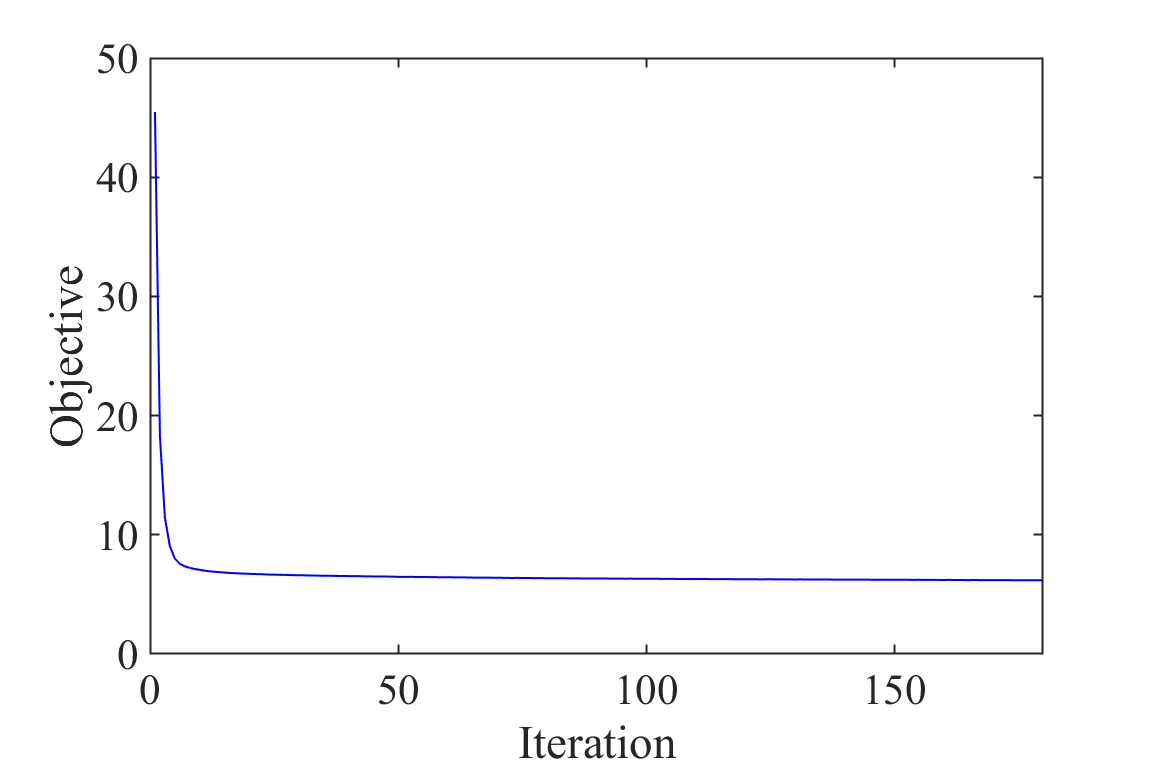}};
    \node[inner sep=0pt] (plot) at (2,0)
    {\includegraphics[scale=0.1]{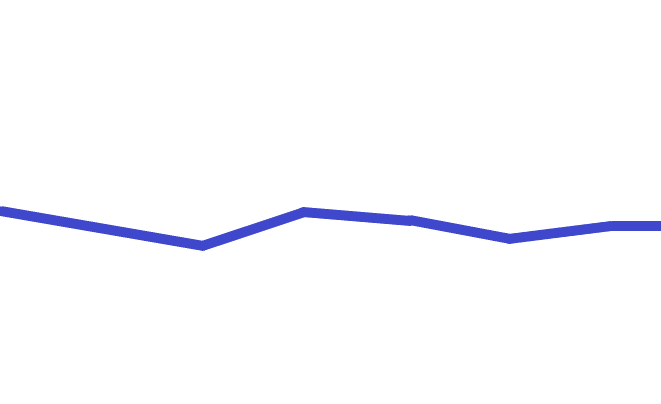}}; 
    \draw[draw=black,thick] (1.1,-0.5) rectangle ++(1.75,1); 
    \draw[draw=black,thick] (2.15,-1.5) rectangle ++(0.275,0.2); 
    \draw[thick] (2.15,-1.3) -- (1.1,-0.5); 
    \draw[thick] (2.425,-1.3) -- (2.85,-0.5); 
    \node[] at (0.8,0.5) {\scriptsize{6.45}}; 
    \node[] at (0.8,-0.5) {\scriptsize{6.10}}; 
    \end{tikzpicture}
    \caption{Reduction of objective in Example III using Adam with a zoomed-in portion showing some oscillations at the end of the optimization due to the stochastic nature of gradients. }
    \label{fig:exIII_obj} 
\end{figure}
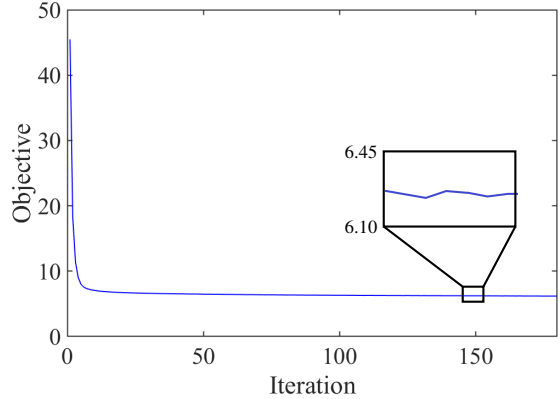 

\begin{figure}[!htb]
    \centering
    \begin{tikzpicture}
        \node[inner sep=0pt] (structure) at (0,0){\includegraphics[scale=0.275]{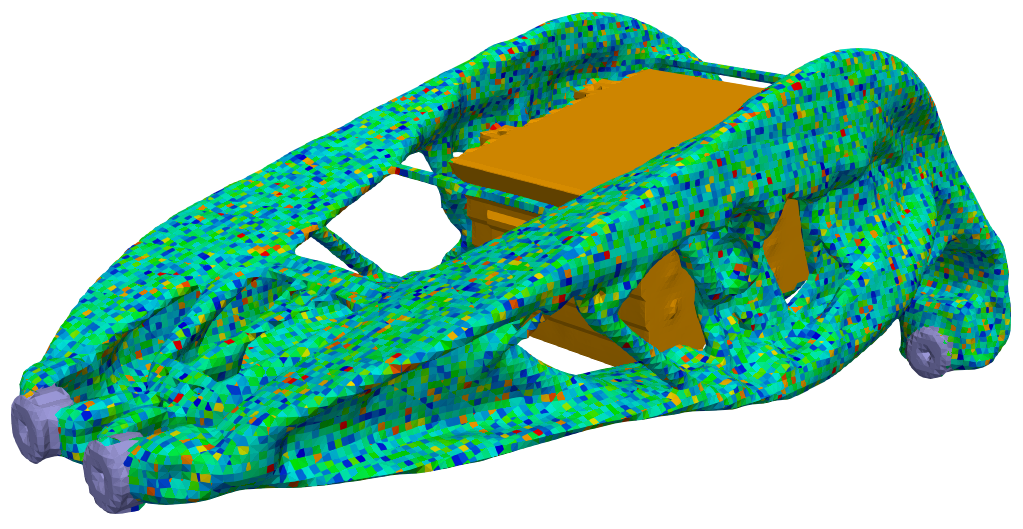}};
        \node[inner sep=0pt] (colorbar) at (4.5,0)
    {\includegraphics[scale=0.375]{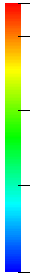}}; 
    \node at (4.75,1.8) {$C_{1111}/C_{2222}$}; 
    \node at (4.95,1.35) {$2.2$}; 
    \node at (4.95,1.0) {$2.0$}; 
    \node at (4.95,0.275) {$1.5$}; 
    \node at (4.95,-0.45) {$1.0$}; 
    \node at (4.95,-1.3) {$0.4$}; 
    \draw[-latex] (-4.5,-2) -- (-4,-1.75); 
    \draw[-latex] (-4.5,-2) -- (-4.5,-1.5); 
    \draw[-latex] (-4.5,-2) -- (-4.1,-2.25); 
    \node[] at (-3.9,-1.65) {$x$}; 
    \node[] at (-4.525,-1.4) {$y$}; 
    \node[] at (-3.9,-2.25) {$z$}; 
        \end{tikzpicture}
    \caption{The color shading shows the ratio $C_{1111}/C_{2222}$ of the first two diagonal elements of the constitutive tensor $\Chom$ for the designed structure in Example III for one random layout of the microstructure. }
    \label{fig:exIII_adam_C11_C22} 
\end{figure}
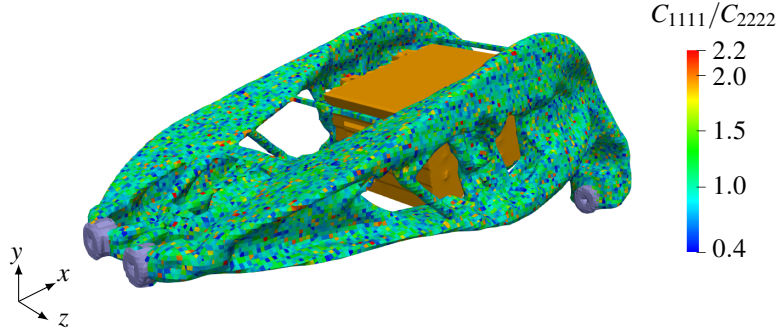

\subsubsection{Microstructure Scenarios} 
The structure is assumed to be made of Ti-6Al-4V alloy with impurities in it. The alloy has elastic modulus $E_\mathrm{Ti-6Al-4V}=1.138\times10^5$ MPa and Poisson's ratio $\nu_\mathrm{Ti-6Al-4V} = 0.342$. The impurities in the material are assumed to have an elastic modulus $E_\mathrm{imp} = 0.02E_\mathrm{Ti-6Al-4V}$ and Poisson's ratio $\nu_\mathrm{imp}=0.3$. 
We generate 200 random microstructures with $T=4\pi$ and $K=25$ 
similar to shown in Figure \ref{fig:3dmat1} (see Section \ref{sec:micro_gen}). The first-order computational homogenization uses a discretization of $100\times100\times100$ and six linear analyses for each of these microstructures to estimate the homogenized constitutive tensor. For one layout of the microstructure configuration, we randomly assign each element in the finite element model one of the 200 possible microstructures. Hence, the stochastic dimension of this problem is extremely high. Instead, in our proposed approach, we only consider one random material layout per iteration and use a stochastic gradient-based approach. 

\subsubsection{Optimization Results} 

Using Adam with a step size $\eta=0.05$ and penalties $\kappaa = [10^4,10^4]^T$ to implement the constraints, we obtain the final design shown in Figure \ref{fig:exIII_adam}. When compared to a deterministic design for Ti-6Al-4V alloy with no impurity using GCMMA (see Figure \ref{fig:exIII_det}), the Adam design under microstructural uncertainty produces a design with thicker members. Figure \ref{fig:exIII_obj} shows the progression of the objective during optimization. Note that the Adam algorithm satisfies the mass inequality constraint. The violation in the equality constraint to ensure connections at the bolts only is $\mathcal{O}(10^{-6})$. In Figure \ref{fig:exIII_adam_C11_C22}, we show the ratio of the first two diagonal elements of the constitutive tensor, illustrating the variability in the microstructural properties. 

{We also use GCMMA with the same step size as above, but
GCMMA fails to converge in this example as it removes material from the design domain abruptly based on the stochastic gradients of the current optimization iteration. Smaller step sizes for GCMMA, on the other hand, slow down the convergence significantly. Similarly, increasing $n_s$ does not result in convergence within a reasonable computational budget for GCMMA.} In one of our previous works \cite{de2019topology}, we showed that GCMMA requires more accurate gradients to converge when the variance of the stochastic gradients is significantly large. As we are using only one random microstructure layout per iteration compared to four in the previous examples, GCMMA fails here. 


\section{Conclusions} 

Uncertainties in the microstructure of composite materials are frequently encountered across engineering applications. To design structures that are robust with respect to microstructural uncertainty, the mean performance needs to be estimated. Using standard gradient-based optimization methods requires accurate estimations of objectives, constraints, and their design sensitivities. Computing this information by the standard Monte Carlo methods may require a large number of samples increasing the computational cost. 
In this paper, to significantly reduce the computational cost, we presented an approach where only {a handful ($\sim\mathcal{O}(1)$) of the possible microstructural configurations selected randomly} need to be considered per optimization iteration. This results in stochastic gradients, which we used with two algorithms, namely Adam and GCMMA. 
{The proposed approach, to the best of our knowledge, is the first to tackle such TO problems, as they are beyond the capabilities of current tools, which require a large computational budget.} 
We illustrated this approach with a two-dimensional and two three-dimensional problems. These examples show the effectiveness of the proposed approach in reducing the computational cost of the optimization. Among the two algorithms, in the presence of large uncertainty, Adam outperforms GCMMA. Further, without an accurate estimation of the gradients, GCMMA is prone to diverge. {Increasing the number of realizations of the uncertain parameters to get more accurate estimates of the gradients in GCMMA, however, results in an impossible computational budget.} In the future, we plan to incorporate additive manufacturing constraints into the proposed topology optimization process.

\section*{Acknowledgment} 
The authors acknowledge the support of the Defense Advanced Research Projects Agency (DARPA) under the TRADES
program (agreement HR0011-17-2-0022). The opinions and conclusions presented in this paper are those of
the authors and do not necessarily reflect the views of DARPA.

\bibliographystyle{unsrt}
\bibliography{Multiscale_TO.bib}

\end{document}